\documentclass[11pt]{article}
\usepackage{amsmath,amsfonts, amssymb,amsthm,multicol,graphicx,tikz, enumitem}
\usepackage{makeidx}
\usepackage{mleftright}
\usetikzlibrary{calc}

\DeclareMathOperator{\Var}{Var}

\renewcommand{\a}{\alpha}

\newcommand{\ep}{\epsilon}
\newcommand{\ZZ}{\mathbb{Z}}
\newcommand{\NN}{\mathbb{N}} 
\newcommand{\PP}{\mathbb{P}}

\newcommand{\cA}{\mathcal{A}}
\newcommand{\cB}{\mathcal{B}}
\newcommand{\cF}{\mathcal{F}}

\newcommand{\cE}{\mathcal{E}}

\newcommand{\ccI}{\mathcal{I}}
\newcommand{\cI}{\mathcal{U}}

\newcommand{\cX}{\mathcal{X}}

\newcommand{\cZ}{\mathcal{Z}}
\newcommand{\cU}{\mathcal{U}}
\newcommand{\cS}{\mathcal{S}}

\newcommand{\EE}{\mathbb{E}}

\newcommand{\hm}{\hat{m}}
\newcommand{\hr}{\hat{r}}

\newcommand{\hk}{\hat{k}}

\newcommand{\vd}{\vec{d}}
\newcommand{\vchi}{\vec{\chi}}

\newcommand{\tk}{\tilde{k}}

\newtheorem{prop}{Proposition}   
\newtheorem{lem}[prop]{Lemma} 
\newtheorem{thm}[prop]{Theorem} 
\newtheorem{define}[prop]{Definition} 

\newtheorem{claim}[prop]{Claim}

\title{A Note on Two-Point Concentration of the Independence Number of $G_{n,m}$}
\author{Tom Bohman\thanks{Supported in part by NSF grant DMS-2246907} \  and Jakob Hofstad\thanks{Supported in part by NSF grant DMS-2246907} }

\begin{document}

\maketitle

\begin{abstract}
    We show that the independence number of $ G_{n,m}$ is concentrated on two values for $ n^{5/4+ \epsilon} < m \le \binom{n}{2}$. This result establishes a distinction between $G_{n,m}$ and $G_{n,p}$ with $p = m/ \binom{n}{2}$ in the regime $ n^{5/4 + \epsilon} < m< n^{4/3}$. In this regime the independence number of $ G_{n,m}$ is concentrated on two values while the independence number of $ G_{n,p}$ is not; indeed, for $p$ in this regime variations in $ \alpha( G_{n,p})$ are determined by variations in the number of edges in $ G_{n,p}$.
\end{abstract}

\section{Introduction}
An independent set in a graph $G=(V,E)$ is a set of vertices that contains no edge. The independence number of $ G$, denoted $ \alpha(G)$, is the maximum number of vertices in an independent set in $G$. In this work we address the independence number of $ G_{n,m}$, which is a graph drawn uniformly at random from the collection of all graphs on a set of $n$ labeled vertices with $m$ edges. We also consider the binomial random graph $ G_{n,p}$ which is chosen from the space of all graphs on $n$ labeled vertices such that each potential edge appears as an edge with probability $p$, independently of all other potential edges. Given a sequence of random graphs (e.g. $ G_{n,m}$ with $m = m(n)$ or $ G_{n,p} $ with $ p = p(n)$) we say that an event occurs with high probability (whp) if the probability as $n$ tends to infinity of this event is 1. A random variable $X$ defined on a random graph is concentrated on $\ell = \ell(n)$ values if there is a sequence $k = k(n)$ such that $ X \in \{ k , k+1, \dots, k+ \ell-1\}$ whp.

Early in the study of random graphs it was observed that certain natural graph parameters can be sharply concentrated. The sharpest version of this phenomenon is so-called 2-point concentration in which the given parameter is one of two consecutive values with high probability. This was first observed in the 1970's by Matula \cite{matula} as well as Bollob\'as and Erd\H{o}s \cite{1976cliques}, who showed that $  \alpha (G_{n,p})$ is concentrated on two values when $p$ is a constant. An influential series of works starting with Shamir and Spencer \cite{ss} and culminating in the result of Alon and Krivelevich \cite{ak} establish two-point concentration of the chromatic number of $ G_{n,p}$ for $ p \le n^{-1/2 - \epsilon}$ (also see \cite{luczak} \cite{an} \cite{cps} \cite{lutz}). Heckel and Riordan recently showed the chromatic number of $ G_{n,1/2}$ is not concentration on any interval of length $ n^{1/2 - \epsilon}$ \cite{h} \cite{hr}. They also have a fascinating conjecture regarding the fine detail of the extent of concentration of the chromatic number of $ G_{n,1/2}$. The extent of concentration of the chromatic number of $ G_{n,p}$ for $ n^{-1/2} \le p \ll 1$ remains a central open question. There has also been interest in the extent of concentration of the domination number of random graphs \cite{bwz} \cite{gls} \cite{wg}, and recent work on variations in the matching number of sparse random graphs \cite{k} \cite{gkss}.

Returning to the independence number, the authors recently observed that 2-point concentration of the independence number of $ G_{n,p}$ holds for $ p \ge n^{-2/3 + \epsilon}$ \cite{bhindy}. Furthermore, an argument of Sah and Sawhney shows that this behavior breaks down for $ p < n^{-2/3}$. Thus the authors' result for 2-point concentration of the independence number of $ G_{n,p}$ is essentially best possible. 

The argument of Sah and Sawhney uses that fact that if $ p $ is close to $p'$ then $ G_{n,p}$ and  $ G_{n,p'}$ can be coupled in such a way that the two graphs are equal with some significant probability. The argument does not apply to $ G_{n,m}$ as such a coupling is not available. This leaves open the possibility that there is $ m = m(n)$ such that $ \alpha( G_{n,m})$ is concentrated on 2 values while $ \alpha( G_{n,p})$ with $ p = m/ \binom{n}{2}$ is not concentrated on 2 values. Our main result is that there is a range for the parameters where this is indeed the case.
\begin{thm} \label{thm:main}
    Let $ \epsilon > 0$ and $ m = m(n) > n^{5/4 + \epsilon}$. Then there exists $ k_0 = k_0(n)$ such that 
    \[ \lim_{n \to \infty} \PP \left( \alpha( G_{n,m}) \in \{k_0-1,k_0\} \right) =1. \]
\end{thm}
\noindent
As a corollary we pin down the extent of concentration of $ \alpha( G_{n,p})$ for $ p \ge n^{-3/4 + \epsilon}$.
\begin{thm} \label{thm:Gnp} 
    \begin{itemize}
     \item[(i)] If $p = \omega( (\log n)^{2/3} n^{-2/3} )$ then  there exists $ k_0 = k_0(n)$ such that 
    \[ \lim_{n \to \infty} \PP \left( \alpha( G_{n,p}) \in \{k_0-1,k_0\} \right) =1. \]
        \item[(ii)] Let $ \epsilon > 0$ and suppose $ n^{-3/4+ \epsilon} < p = O\left( ( \log n)^{2/3} n^{-2/3} \right)$. For any function $ \ell= \omega(1)$ 
        the independence number 
        $\alpha(G_{n,p})$ is concentrated on $\ell (\log n) p^{-3/2} n^{-1}$ values. Furthermore, $\alpha(G_{n,p})$ is not concentrated on any interval of length 
        $ C \log(n) p^{-3/2} n^{-1} $ for any constant $C$.  
    \end{itemize}
\end{thm}

We prove Theorem~\ref{thm:main} by applying the first and second moment methods to a generalized version of an independent set that we call an extended independent set. These objects are very similar to the structures that Ding, Sly and Sun used in their breakthrough work on the independence numbers of random regular graphs with large constant degree \cite{dss}. Extended independent sets are introduced and discussed in the following Subsection. Our second moment argument requires precise estimates for the probability that a pair of extended independent sets appear in $ G_{n,m}$. Our estimates make crucial use of a novel variation on Janson's inequality that was recently introduced by Bohman, Warnke and Zhu in the context of establishing two-point concentration of the domination number of $ G_{n,p}$ \cite{bwz}.

\subsection{Extended independent sets}

We establish two-point concentration of $ \alpha( G_{n,m})$ by applying the second moment method. The choice of random variable to which we apply the second moment is an important part of the argument. The natural first choice is the number of independent sets of a fixed size $k$. Let $ k_V$ be the smallest value of $k$ for which the expected number of independent sets of size $k$ is at most $ n^{\epsilon}$ (the `V' stands for vanilla). Its turns out that if $ m > n^{5/3+ \epsilon}$ one can establish two-point concentration of $ \alpha( G_{n,m})$ around $ k_V$ by working directly to the random variable that counts the number of independent sets of size $k$. For smaller values of $m$ this approach does not succeed due to variance coming from potential independent sets that have large intersection. Working with maximal independent sets allows one to extend the argument for two-point concentration to $ m > n^{3/2 + \epsilon}$. 

Below roughly $  n^{3/2} $ random edges the behavior of the independence number of $ G_{n,m}$ changes. In particular, the independence number is less then $k_V$ with high probability in this regime. In a previous work, the authors showed that if $ n^{4/3+ \epsilon} < m < n^{3/2} $ then the independence number is driven by the appearance of structures called augmented independent sets, where an augmented independent set of order $k$ in a graph $ G$ is a set $K$ of $k +r$ vertices such that $ G[K]$ is a matching with $r$ edges combined with $k-r$ isolated vertices and every vertex outside of $K$ has at least two neighbors in $K$. It turns out that the maximum order of an augmented independent set in $G$ is equal to $ \alpha(G)$. Note that an augmented independent set includes $2^r$ independent sets of size of $k$. Thus the expected number of independent sets of size $k$ can be large even though the probability that there is a single augmented independent set of order $k$ is vanishingly small.

Below roughly $ n^{4/3}$ random edges we need to consider a more general object. To motivate this definition, consider an augmented independent
set on vertex set $K$ with matching $M$. If there is a vertex $ v \not\in K$ that has exactly two neighbors
in $K$ such that one of the 
neighbors is a vertex $x$ such that the edge $xy$ is in $M$ while the other neighbor is 
$ w \in K \setminus ( \cup_{e \in M}e )$ then there are potentially two augmented independent sets in play. These are the
vertex set $K$ with matching $M$ and the vertex set $ (K \cup \{v\}) \setminus \{x\}$ with matching $(M \cup \{wv\}) \setminus \{xy\}$. 
This has the potential to generate more variance in the count of augmented independent sets. It turns out this situation
is unlikely when $ m > n^{4/3}$, but becomes more prominent for smaller $m$. In order to extend our previous two-point concentration result to smaller values of $m$, we further generalize the definition of an augmented independent set. In the situation presented here we would move the vertex $v$ into the vertex set $K$ and add the edge $ vw$ to the matching. Note
that this then allows edges that connect pairs of edges in $M$. If $v$ is a vertex in a graph $G$ then we use $ {\mathcal N}(v)$ to
denote the neighborhood of $v$.
\begin{define}
Let $G$ be a graph.
An {\bf extended independent set of order $k$} in $G$ is a set $K$ of $k + r$ vertices for some $ 0 \le r \leq k $ such that 
\begin{enumerate}
    \item $G[K]$ has no cycles, exactly $2r$ non-isolated vertices, and a perfect matching among those $2r$ vertices, and
    \item For every vertex $v \not\in K$, the intersection $ \mathcal{N}(v) \cap K $ contains 
    \begin{itemize}
        \item At least two isolated vertices of $G[K]$ or
        \item At least two vertices of the same connected component of $G[K]$.
    \end{itemize}
\end{enumerate}
\end{define}

\begin{lem} \label{Lem:MNI connection}
For any graph $G$, $\a(G) = k$ if and only if $G$ has an extended independent set of order $k$, but no extended independent
set of any larger order.
\end{lem}
\begin{proof}
Let $ \hat{\alpha}(G)$ be the largest $k$ for which there is an extended independent set of order $k$. 

Let $K$ be the vertex set of an extended independent set of maximum order.
Each non-singleton component of $G[K]$ is a tree and has a perfect matching, hence $G[K]$ contains an
independent set on $k$ vertices, so we have
\[  \hat{\alpha}(G) \le \alpha(G).   \]
Now suppose that $G$ has a maximum independent set $S$ of size $k$. Let $T$ be a maximum set of vertices that contains $S$ and has the additional property that the connected components on $G[T]$ are all isolated vertices or trees with perfect matchings. We claim that $T$ is an extended independent set of order $k$. Assume not. Then some $v \not\in T$ must have 0 or 1 neighbors which are isolated in $G[T]$, and at most 1 neighbor in each tree in $G[T]$. If $\mathcal{N}(v) \cap T$ includes no isolated vertices of $G[T]$, then an independent set of size $k+1$ can be found in $G[T + \{v\}]$ by choosing the partite set of each tree in $G[T]$ which does not include any vertices of $\mathcal{N}(v)$, $v$ itself, and the isolated vertices of $G[T]$. This contradicts the choice of $S$ as a maximum-size independent set of $G$. If $\mathcal{N}(v) \cap T$ includes just one isolated vertex in $G[T]$ (call this vertex $w$), then $T \cup \{w\}$ has the property that each component of $ G[ T \cup \{w\}]$ is an isolated vertex or a tree with a perfect matching. (The edge $vw$ is a matching edge in a tree that includes some set - possibly empty - of the trees in $G[T]$). This contradicts the maximality of $T$. Therefore,
\[ \alpha(G) \le \hat{\alpha}(G). \]
\end{proof}

\begin{define}
    Let the random variable $ Y_{k,r} (G) $ be the number of extended independent sets of order $k$ with matchings of size $r$ in $G$. 
\end{define}
\noindent
While our main interest is in the variable $ Y_{k,r}$, we mostly work with a variable that counts sets with slightly different requirements; specifically, these sets allow for cycles in the interior of $K$ and do not allow vertices outside of $K$ to satisfy their degree requirement by having two neighbors within some component of $ G[K]$.
\begin{define}
For $ 0 \le r \le k$
let the random variable $ X_{k,r}$ be the number of sets $K$ of $k + r$ vertices such that 
\begin{enumerate}
    \item The induced graph $G[K]$ has exactly $2r$ non-isolated vertices and a perfect matching among those $2r$ vertices, and
    \item Every vertex $v \not\in K$ has at least neighbors among the collection of $k-r$ vertices in $K$ that are
    isolated in $G[K]$.
\end{enumerate}
\end{define}
\noindent
Note that the variable $ X_{k,r}$ has the advantage that in order to determine whether or not a vertex set $K$ and matching $M$ is counted by $ X_{k,r}$ we do not need to observe the potential edges joining the edges in $M$. 

We need to set some further definitions and conventions in order to state our main technical results. 
In order to make comparisons (often simply heuristic, but also occasionally rigorous) with $ G_{n,p}$ we set $ p = m/\binom{n}{2}$ throughout
this work. We restrict our attention to $k$ that satisfy
\begin{equation} 
k = \frac{2}{p} \left[ \log(np) - \log\log(np) + \log(e/2) \pm o(1) \right]. \label{equ: k used in this range}
\end{equation}
Note that a celebrated result of Frieze establishes that $ \alpha( G_{n,p}) $ satisfies this condition \cite{frieze}. 
Set
\[  r_1 = 8 ( \log np)^3/( n p^2) ,\]
and define
$ k_0 = k_0(n,m)$ to be the smallest $k$ such that $ \EE[ X_{k,r}] < n^{\epsilon}$ for all $r \leq r_1$. For any $k$ that satisfies (\ref{equ: k used in this range}) 
we set $ r_0(k)$ to be the value of $r \le r_1 $ that maximizes $ \EE[ X_{k,r}] $.
We are now ready to state the two main technical results in this paper.
\begin{thm} \label{thm:techy1}
Let $ \epsilon > 0$ and $  n^{5/4 + \epsilon} < m = m(n) < n^{4/3 + \epsilon}$.  If $ k$ satisfies (\ref{equ: k used in this range}),
$ r = O( ( \log n)^3 / (np^2))$ and $ \EE[ X_{k,r}] \to \infty$ then
\[ \Var[ X_{k,r}] = o ( \EE[ X_{k,r}]^2).\]
\end{thm}
\begin{thm} \label{thm:techy2}
Let $ \epsilon > 0$ and $  n^{5/4 + \epsilon} < m = m(n) < n^{4/3 + \epsilon}$.  If $ k$ satisfies (\ref{equ: k used in this range}) then
\[ \sum_{r=0}^{ r_1} \EE[ X_{k, r}] = \Theta \left( \sqrt{ r_0(k) }  \cdot \EE[ X_{k, r_0(k)}] \right). \]
    \end{thm}
\noindent
While these two results comprize the main technical contributions of this work, note that
some work needs to be done in order to deduce Theorem~\ref{thm:main} from Theorem~\ref{thm:techy1} and 
Theorem~\ref{thm:techy2}. In particular, we need to examine the
relationship between $ Y_{k,r}$ and $ X_{k,r}$ and to extend the estimate given in Theorem~\ref{thm:techy2} to
a sum over $ 0 < r < k$ and $ k_0+1 \le k \le k_V$ 

The remainder of the paper is organized as follows. We close the Introduction with some commentary
on the extent of concentration of $ \alpha( G_{n,m})$ for $ m = O( n^{5/4})$ and potential obstacles to
extending Theorem~\ref{thm:main} to that regime. Section~2 contains a number of facts and 
observations that are used throughout the rest of the paper. The proof of Theorem~\ref{thm:techy2} is given in 
Section~3, and the proof of Theorem~\ref{thm:techy1} is given in Section~4. The proof of Theorem~\ref{thm:techy1} requires
two applications of the Poisson paradigm. These are crucial parts of the argument, and they are presented in Section~5. 
Theorem~\ref{thm:main} is derived from Theorems~\ref{thm:techy1}~and~\ref{thm:techy2} in Section~6. Section~7 contains the proof of Theorem~\ref{thm:Gnp}.

\subsection{Discussion}

The extent of concentration of $ \alpha(G_{n,m})$ for $  n (\log n)/2 < m < n^{5/4} $ 
and $ \alpha(G_{n,p})$ for $ (\log n)/n < p < n^{-3/4}$ remain interesting open questions. 
(In a previous work \cite{bhindy}, the authors showed that $\alpha( G_{n,p})$ with $ p = c/n$ where $c$ is a constant has variations of size $ \Theta( \sqrt{n})$. This argument exploits vertices of degree 1, and we focus this discussion on values of the parameter where such arguments do not apply.) In this Section we make some remarks about this problem and note some obstacles to applying the arguments given in this work to random graphs with fewer edges. There are reasons to believe that 2-point concentration of $ G_{n,m}$ breaks down at about $ m = n^{5/4}$. There are also reasons to believe that 2-point concentration of $ G_{n,m}$ persists even for $m$ close to $ n (\log n)$. We attempt to present both perspectives in the section.

We begin by giving a rough heuristic explanation for the change of behavior of the variation in $ \alpha( G_{n,p})$ at about $ p = n^{-2/3}$. We focus on the variance in the number of independent sets of size $k$ where $k$ satisfies (\ref{equ: k used in this range}). Throughout this discussion we work with counts of independent sets of size $k$, but we assume that we can ignore the contribution to the variance coming from independent sets with large intersection. (N.b. This is a significant assumption as the large intersection case is a substantial part of both this work and its predecessor \cite{bhindy}).)  The rationale for this assumption is that we should be able to eliminate these contributions by working with extended independent sets (or some related structure). So we focus on the contribution to the variance coming from pairs of potential independent sets of size $k$ that have typical intersection size. So consider potential independent sets $A$ and $B$ of size $k$ such that $ |A \cap B| = k^2/n$.  Let $ \ccI_A$ be the event that $A$ is an independent set. We have
\begin{multline*}
\PP( \ccI _A \wedge \ccI_B) = (1-p)^{2\binom{k}{2} - \binom{ k^2/n}{2}}
\approx \PP(\ccI_A) \PP (\ccI_B) \exp\{ p k^4/( 2n^2) \} \\
=  \PP(\ccI_A) \PP( \ccI_B) \exp \left\{ \tilde\Theta \left(  \frac{ 1 }{ p^3 n^2} \right) \right\}. \end{multline*} 
Thus, we expect such terms to start to contribute to the variance at around $ p = n^{-2/3}$. This is where two point concentration of $ \alpha( G_{n,p})$ breaks down.

Now we apply the same line of reasoning to $ G_{n,m}$. Again let $A$ and $B$ be potential independent sets of size $k$ where $k$ satisfies (\ref{equ: k used in this range}), but here we will allow some variation of the size of the intersection of $A$ and $B$. So we set $ |A \cap B| = i $. Here we have
\begin{multline*}
\frac{ \PP( \ccI_A \wedge \ccI_B)}{ \PP( \ccI_A ) \PP( \ccI_B)} 
= \frac{ \binom{ \binom{n}{2} - 2 \binom{k}{2} + \binom{ i }{2}}{ m} \binom{ \binom{n}{2}}{m}}{ \binom{ \binom{n}{2} - \binom{k}{2}}{m}  ^2}
= \prod_{a=0}^{m-1} \frac{  \left( \binom{n}{2} - 2 \binom{k}{2} + \binom{ i }{2} -a  \right) \left( \binom{n}{2} -a  \right)  }{ \left( \binom{n}{2} - \binom{k}{2} - a \right)^2} \\
= \prod_{a=0}^{m-1} \left( 1 - \frac{ \binom{k}{2}^2 - \binom{i}{2}\left( \binom{n}{2} - a \right)}{ \left( \binom{n}{2} - \binom{k}{2} - a \right)^2} \right)
. 
\end{multline*} 
Note that if we set $ i = k^2/n$ then we have 
\[  \binom{k}{2}^2 - \binom{i}{2}\left( \binom{n}{2} - a \right) =  \frac{ k^2n}{4} \pm O \left( (\log n) k^3 \right),    \]
and there is a slight negative correlation between the events $ \ccI_A$ and $\ccI_B$. 
(To put this in other terms: If we condition on the event $ \ccI_A$ then the probability that each potential edge among the remaining $ \binom{n}{2} - \binom{k}{2}$ vertices appears as an edge goes up slightly. This makes the event $ \ccI_B$ less likely. On the other hand, conditioning on $ \ccI_A$ assures that certain edges within $B$ do not appear. It turns out that these two effects balance out almost perfectly.) Now we consider what happens to this computation when we consider slightly larger intersections; in particular, consider $i = k^2/n + k / n^{1/2}$. (Of course, there are many pairs $A,B$ that have such intersection sizes.) Now we have
\[  \binom{k}{2}^2 - \binom{i}{2}\left( \binom{n}{2} - a \right) = 
\frac{ -k^3 n^{1/2}}{2} - O( k^2n).  \]
So we conclude that in this situation we have
\[\frac{ \PP( \ccI_A \wedge \ccI_B)}{ \PP( \ccI_A ) \PP( \ccI_B)} \approx \exp \left\{ \frac{ p k^3}{ n^{3/2}}  \right\}  =   \exp \left\{ \tilde{\Theta}\left( \frac{ 1}{ p^2 n^{3/2}} \right)  \right\}.  \]
This term becomes significant at about $ p = n^{-3/4}$. If we draw an analogy with $ G_{n,p}$, we might expect that this is the point where two-point concentration of $ \alpha( G_{n,m})$ breaks down. Furthermore, the arguments presented in this work to do not handle pairs of potential extended independent sets that have intersections of this size when $ m < n^{5/4}$. This is one of the most significant `gaps' if we attempt to apply the methods to developed in this work to smaller $m$.

The argument in favor of two-point concentration holding for $m$ smaller than
$ n^{5/4}$ rests on two considerations. The first is the breakthrough work of
Ding, Sly and Sun on the independence number of the random regular graph with large constant degree \cite{dss}. They determine the expected size of this independence number and show that the variations around this expectations are of constant order. Furthermore, their argument makes use of structures - inspired by the notion of replica symmetry breaking from statistical physics - that are nearly identical to the extended independent sets that we work with here. It is natural to suspect that the maximum size of such structures (and hence the independence number) should be sharply concentrated for a wide range of sparse or relatively sparse random graph models with a fixed number of edges. The second consideration is the fact that we are not aware of {\em any} anti-concentration argument for the independence number of sparse random graphs with a fixed number of edges that does not require a significant number of vertices of degree 1. (The authors presented an anti-concentration argument that exploits vertices of degree 1 in their prior work \cite{bhindy}.) For example, it would be interesting to show that $ \alpha( G_{n,n(\log n)})$ is not concentrated on any interval of constant length.

We close this discussion by mentioning a second obstacle to extending the arguments given here to $m$ smaller than $n^{5/4}$. Consider a potential
extended independent set $K$ with matching $M$. In order to calculate the probability that every vertex outside of $K$ has two neighbors
in $ K \setminus ( \cup_{e \in M} e) $ we establish an estimate for the number of unbalanced bipartite graphs with a certain number of edges and
minimum degree 2 for all vertices in one of the parts in the bipartition. The lower bound for this estimate is given in Theorem~\ref{thm: enumeration} below, which we prove using
methods introduced by Pittel and Wormald for counting graphs with minimum degree 2 \cite{pittelwormald}. There is a second order term in the bound that becomes relevant at $ m = n^{5/4}$, and this term is not present in our upper bound which we achieve by simply passing to $ G_{n,p}$. 
Furthermore, when we apply Theorem~\ref{thm: enumeration} we condition on the number edges between $ K \setminus ( \cup_{e \in M} e)$ and $ V \setminus K$,
and when $ m < n^{5/4}$ the estimate itself varies significantly when we change this number of edges by a less than a standard deviation. Accounting for such variations introduces yet another difficulty.
So applications of the proof given in this work to smaller values of $m$ would require more sophisticated methods for 
counting bipartite graphs with a fixed number of edges and fixed minimum degree conditions.

\section{Preliminaries}

Recall that it is useful to couch our calculations for $ G_{n,m}$ in the parameters of 
the associated binomial model. For that purpose we set $p = m \binom{n}{2}^{-1}$
throughout. We use standard asymptotic notation; for example, if $f$ and $g$ are non-negative functions of $n$ then we 
write $ f \sim g$ if $ \lim_{n \to \infty} f/g =1$ , we write $ f = o(g)$ is this limit is $ 0$ and we write 
$ f = \omega(g)$ if this limit is $ \infty$. We also follow the standard notational conventions regarding 
$ O(g) $, $ \Omega(g)$ and $\Theta(g)$, and we introduce tildes into this notation when we suppress factors that are logarithmic in $n$.
%

\subsection{Properties of $k$ that satisfy (\ref{equ: k used in this range})}

\begin{lem} \label{lem: first preliminary equation}
    For $k$ that is a function of $n$ and $p$ (which itself is a positive function of $n$ such that $p = \omega(n^{-1})$), the following two statements are equivalent:
    \begin{align}
        &k = \frac{2}{p} \left[ \log(np) - \log\log(np) + \log(e/2) \pm o(1) \right] \\
        &\frac{n}{k}e^{-kp/2 + 1} \sim 1. \label{equ: initial k pre-approx}
    \end{align}
\end{lem}
\begin{proof}
    Say that, for some function $a = a(n)$, we have
    \begin{equation*}
    k = \frac{2}{p} \left[ \log(np) - \log\log(np) + \log(e/2) + a \right]. 
    \end{equation*}
Then we have
\begin{align*}
\frac{n}{k} e^{-kp/2+1} &= \frac{n}k \exp\{ - \log(np) + \log( \log(np)) - \log( e/2) +1 - a\} 
\\ &= e^{-a} \left(\frac{ \log(np) }{ \log(np) - \log\log(np)  + \log(e/2) + a} \right).
\end{align*}
The final expression is $1+o(1)$ if and only if $a = o(1)$.
\end{proof}
Note that it follows from Lemma~\ref{lem: first preliminary equation} that if $k$ satisfies
(\ref{equ: k used in this range}) then we have

\begin{equation} 
 (1-p)^{k/2} \sim e^{-kp/2} \sim \frac{k}{e n}. \label{equ: k pre-approx}
\end{equation}
and

\begin{equation}
\label{eq:maxxy}
 \left[ 1 - (1-p)^k - pk(1-p)^{k-1} \right] ^{n-k} = e^{ - (1 +o(1)) \frac{ p k^3}{ e^2 n} }.
 \end{equation}

Recall that we define $k_{V}$ (V for ``Vanilla") to the smallest $k$ for which the expected number of
independent sets of size $k$ is at most $ n^{\epsilon}$. We now observe that $ k_V$ satisfies (\ref{equ: k used in this range}) and
that the expected
number of independent sets of size $k_V+1$ is $o(1)$. 
Hence $ k_V$ serves as both a natural
reference point and trivial upper bound on the independence number. 

\begin{lem} \label{lem: k vanilla}
    $k_V = \frac{2}{p} \left[ \log(np) - \log\log(np) + \log(e/2) \pm o(1) \right]$ where $ p = m / \binom{n}{2}$. In other words, $ k_V$ satisfies (\ref{equ: k used in this range}).
\end{lem}
\begin{proof}
The expected number of independent sets of size $k$ equals $\binom{n}{k} \frac{\binom{\binom{n}{2} - \binom{k}{2}}{m}}{\binom{\binom{n}{2}}{m}}$. 
Recalling a standard bound on binomial coefficients (see Lemma~21.2 in \cite{alanmichal}) we have 
\begin{align*}
    \binom{n}{k} \frac{\binom{\binom{n}{2} - \binom{k}{2}}{m}}{\binom{\binom{n}{2}}{m}} &\sim \frac{1}{\sqrt{2\pi k} e^{k^2/(2n)+k^3/(6n^2)}}\left(\frac{ne}{k}\right)^{k} \frac{\left[\binom{n}{2} - \binom{k}{2}\right]_m}{\left[\binom{n}{2}\right]_m} \\&=
    \frac{1}{\sqrt{2\pi k} e^{k^2/(2n)+k^3/(6n^2)}}\left(\frac{ne}{k}\right)^{k} \prod_{i=0}^{m-1}\left(1 - \frac{\binom{k}{2}}{\binom{n}{2} - i}\right) \\&=
    \left((1+o(1)) \frac{ne}{k}\right)^{k}\exp\left\{ -\frac{mk^2}{n^2}(1 + e^{-\Omega(\log(n))})\right\} \\&=
    \left((1+o(1)) \frac{ne}{k}\right)^{k}\exp\left\{ -k^2 p / 2(1 + e^{-\Omega(\log(n))})\right\} \\&=
    \left((1+o(1)) \frac{n}{k} e^{-kp/2+1}\right)^{k}.
\end{align*}
By Lemma \ref{lem: first preliminary equation} (and because the expression inside the parenthesis above is strictly decreasing over $k > 0$), the final expression shifts from $\omega(n)$ to $o(n^{-1})$ when $k$ satisfies (\ref{equ: k used in this range}). The Lemma follows.
\end{proof}

\begin{lem} \label{lem: k_v markov}
The expected number of independent sets in $ G_{n,m}$ of cardinality $ k_V+1$ is $o(1)$.
\end{lem}
\begin{proof}
By considering the ratio of the expected number independent sets of size $ k_V+1$ and the expected
number of independent sets of size $k_V$, we see that
expected number of independent sets of size $ k_V+1$ is at most
\begin{multline*} n^\epsilon \frac{ n -k}{k+1} \cdot \frac{ \binom{ \binom{n}{2} - \binom{k+1}{2}}{m}} { \binom{ \binom{n}{2} - \binom{k}{2}}{m}} 
\le O( n^{1 +\epsilon} k^{-1} ) \frac{n-k}{k+1} \prod_{i=0}^{m-1} \left( 1 -\frac{k}{ \binom{n}{2} - \binom{k}{2} -i}   \right) \\
= O( n^{1 +\epsilon} k^{-1} ) \exp \left\{ - \frac{ mk}{ \binom{n}{2}} + O \left( \frac{ m k^2}{ \binom{n}{2}^2}  \right) \right\}
= O( n^{1 +\epsilon} k^{-1} ) e^{-pk} \\ = O ( n^{1 +\epsilon} k^{-1} \cdot k^2 n^{-2}) = o(1),
\end{multline*}
where $ p = m/ \binom{n}{2}$ and we apply Lemma~\ref{lem: first preliminary equation} in the penultimate estimate.
\end{proof}

\subsection{Passing from $G_{n,m}$ to $G_{n,p}$.} \label{sec:passthru}

While we focus on $ G_{n,m}$ in the bulk of this work, there
are some situations in which we will be able to get useful
approximations by appealing to $ G_{n,p}$. In this subsection we
set forth a general framework for approximations of this kind. We emphasize that
the estimates we establish here only give upper bounds
on the probabilities of events in $ G_{n,m}$. 

In situations where
we are working with both random graphs we use the subscript to
indicate which graph we are working with. So $ \PP_p( \cE)$ is the
probability of the event $ \cE$ in $ G_{n,p}$ and $ \EE_m( X)$ is
the expectation of random variable $X$ in $ G_{n,m}$.

Let $m$ be fixed and set $ p = m/\binom{n}{2}$. 
Also, let $ \cE$ be some increasing event. This event will generally be of the
form ``every vertex in some set of vertices has at least two neighbors in some other set of vertices."
Set $ p' = (m + \xi)/ \binom{n}{2}$ where $ \xi$ will be chosen later. 
We compare $ G_{n,p'}$ and $ G_{n,m}$. Let the random variable $e( G_{n,p'})$ be the number of edges that appear in $ G_{n,p'}$.

We
have
\begin{gather*}
\begin{split}
\PP_{p'} ( \cE) & = \PP_{p'} ( \cE | e(G_{n,p'})  \ge m ) \PP_{p'} ( e(G_{n,p'})  \ge m)
+ \PP_{p'} ( \cE | e(G_{n,p'}) < m ) \PP_{p'} ( e(G_{n,p'}) < m) \\
&\ge \PP_{p'} ( \cE| e( G_{n,p'}) \ge m) \PP_{p'}( e(G_{n,p'}) \ge m ) \\
& \ge  \PP_{m} ( \cE ) \PP_{p'}( e( G_{n,p'}) \ge m ).
\end{split}
\end{gather*}
Therefore, we have
\begin{equation}
\label{eq:passing to Gnp}
 \PP_m( \cE) \le  \frac{ \PP_{p'}( \cE) }{ \PP_{p'}( e(G_{n,p'}) \ge m) } .
\end{equation}

Now, we are often interested in using such a bound 
in situations in which we condition on
some small collection of edges in $ G_{n,m}$. Furthermore, we will want to
have good approximations for $p'$ and 
$ (1-p')^{k}$. This motivates the following Lemma. 
\begin{lem} \label{lem: passing to Gnp}
Let $V$ be a vertex set such that $ |V|=n$ and let $ F,G \subset \binom{V}{2} $ be disjoint. Let $ \cA$ be the event in the space of graphs on $V$ 
that the pairs in $F$ appear as edges and
the pairs in $G$ do not appear. Let $ \cE$ be an increasing event that is independent of the status of the potential edges in $ F \cup G$ 
and suppose $ F = O(k)$ and $ G =O(k^2)$. If we set $p' =p + ( \log n) \max\{ k^2 p n^{-2}, p^{1/2} n^{-1} \}$ then  we have
\[ \PP_m ( \cE \mid \cA ) \le \left(1 + o(1) \right) \PP_{p'}( \cE) \]
and 
\[ (1 -p')^k = (1 -p)^k \left( 1 - O \left( \frac{p^2 k^4}{n^2} \right) - O \left(  \frac{ p^{3/2} k^2}{n} \right)  \right)\]
\end{lem}
\begin{proof}
Set $E' = e( G_{n,p'})$.
Note that we have
\[ p'\left[ \binom{n}{2} - |F| - |G| \right] = m + \Omega \left( (\log n) p^{1/2} n \right).
\]
So the Chernoff bound
implies $ \PP( E' \ge m \mid \cA) = 1 - o(1)$. Following the argument that gave (\ref{eq:passing to Gnp}) above, we have 
\[ \PP_{p'}( \cE) = \PP_{p'} (\cE \mid \cA) \ge \PP_{p'}( \cE \mid \cA \wedge  \{E' \ge m \}) \PP_{p'}( E' \ge m \mid \cA) \ge \PP_m( \cE \mid \cA) (1 - o(1)).\]

For the additional assertion note that we have
\[  \frac{ (1 - p')^k}{ (1 -p)^k} = \left( 1 - \frac{ ( \log n) \max \{ p k^2/n^2, p^{1/2}/n \}}{ (1-p)} \right)^k.  \]

\end{proof}

\subsection{A Chernoff bound for hypergeometric random variables}

Let $ C \subset A $ be finite sets such that $ |A| =a, |C|=c$, 
and let $ 0< b < a $ be an integer. Throughout this work
we encounter situations in which we choose a $b$-element subset of $A$ uniformly at random and
consider the distribution of the number of elements of this subset that fall in the set $C$. The
random variable $H = H(a,b,c)$ that gives the count of elements of the random set that fall 
in $C$ has a hypergeometric distribution;
that is,
\begin{equation}
    \label{eq:hyperdef}
\PP(H = i) = \frac{  \binom{ c}{i} \binom{a-c}{b-i} }{ \binom{a}{b}}.
\end{equation}
Note that we can view $H$ as the sum of $c$ negatively correlated indicator random
variables. It follows that we have $ \EE[H] = \frac{bc}{a} $, and the Chernoff bound applies. 
\begin{thm} \label{thm:hyperchernoff}
Let $ H = H(a,b,c)$ be a hypergeometric random variable as defined in (\ref{eq:hyperdef}) and let
$ \mu = bc/a$. Then for $ t \ge 0$ we have
\[ \PP( H > \mu + t) \le \exp \left\{ - \frac{t^2}{2( \mu + t/3)}  \right\}\]
and for $ 0 \le t \le \mu$ we have
\[  \PP( H < \mu - t) \le \exp \left\{ - \frac{t^2}{2( \mu - t/3)}  \right\}\]
\end{thm}
See Sections 21.4 and 21.5 of \cite{alanmichal} for proofs and discussion.

\subsection{Preliminary computations for $ \EE[ X_{k,r}]$}
\label{sec:expected prelims}

In this Section we gather some simple facts regarding expected value
computations that are used throughout the paper. 

The expectation $\EE[X_{k,r}]$ 
has three natural parts: one counting the total number of ways to choose the possible locations 
of the set $K$ and matching $M$, one giving the probability of choosing the specified edges and non-edges in the interior of $K$, and the last giving the probability that the two-neighbors requirement for outside vertices is satisfied {\it given} that the interior requirements are satisfied. 
In the interest of simplifying computations we treat these three components separately.
The computation of the first two factors is straightforward. Let $ N(k,r)$ denote the number of ways to choose $ K$ and $M$. We have

\begin{equation}
N(k,r) = \binom{n}{k+r} \frac{[k+r]_{2r}}{2^r r!}. \label{equ: N definition}
\end{equation}
Next let $ \cI$ be the event that the specified edges and non-edges appear within $K$. The probability of $\cI$ can be expressed by considering the number of ways one can place the remaining $m-r$ edges in the remaining $\binom{n}{2} - \binom{k+r}{2} + \binom{2r}{2} -r$ positions; call this probability $U(m,k,r)$.  Setting
\[ M_1 = M_1(k,r) = \binom{k+r}{2} - \binom{2r}{2} + r, \]
we have
\begin{equation}
    U(m,k,r) = \PP( \cU) = \frac{\binom{\binom{n}{2} - M_1}{m-r}}{\binom{\binom{n}{2}}{m}}. \label{equ: E definition}
\end{equation}
Finally, let $ \Phi = \Phi_{n,m,k,r}$ be the probability that every vertex outside of $K$ 
has two neighbors in $ K \setminus ( \cup_{e \in M} e)$ conditioned on the event $\cI$. We have
\begin{equation}
    \EE[X_{k,r}] = N(k,r)U(m,k,r) \Phi. \label{equ: expectation of X}
\end{equation}

The probability $ \Phi$ is the most challenging part of the expression for $ \EE[ X_{k,r}]$ to analyze. For each vertex $ v \not\in K$, the probability, conditioned on $ \cI$, that $v$ has 
fewer than two neighbors
in $ K \setminus ( \cup_{e \in M} e)$ is 
\begin{multline*}  
\frac{ \binom{  \binom{n}{2} - M_1 - (k-r) }{ m-r} }{ \binom{ \binom{n}{2} - M_1}{ m-r}  } 
+ (k-r) \frac{ \binom{  \binom{n}{2} - M_1 - (k-r) }{ m-r -1 } }{ \binom{\binom{n}{2} - M_1}{m-r}  } \\
\approx \left( 1 - \frac{ m-r}{ \binom{n}{2} - M_1} \right)^{k-r} + \frac{ (k-r)(m-r)}{  \binom{n}{2} - M_1 }\left( 1 - \frac{ m-r}{ \binom{n}{2} - M_1} \right)^{k-r-1} \\
\approx (1-p)^{k-r} + (k-r)p( 1 - p)^{k-r-1}.
\end{multline*}
Thus, recalling (\ref{eq:maxxy}) and that we restrict our attention to $ r = O( pk^3/n) = o(1/p)$ when we work with $X_{k,r}$, we anticipate that we should have
\[ \Phi \approx \left( 1 -  (1-p)^{k-r} - (k-r)p( 1 - p)^{k-r-1} \right)^{n-k-r}
\approx \exp \left\{  - ( 1 +o(1))\frac{ pk^3}{ e^2n} \right\}. \] 
Note that such an estimate fits in the framework of the Poisson paradigm.

Establishing estimates of this type in the settings of the first and second moment is a central technical 
component of this work. We
get a lower bound on $ \Phi$ by appealing to direct counting of the number of
ways to place edges that satisfy the given minimum degree conditions, following the work
of Pittel and Wormald \cite{pittelwormald}. We establish an upper bound on $ \Phi$ by comparison with $ G_{n,p}$ (n.b. such an upper bound can also be derived from the calculation above, using negative correlation). In
the context of the second moment calculation we only need an upper bound, and there we also pass to $ G_{n,p}$ and use a novel variation on Janson's
inequality that was recently established by Bohman, Warnke and Zhu \cite{bwz}. These arguments are presented in
Section~\ref{sec:poisson}.

In order to replace these loose approximations for $ \Phi = \Phi_{n,m,k,r}$ with proper estimates we 
define 
\[ c = c(m,k,r) = p(k-r) = \frac{ (k-r)m}{ \binom{n}{2}} ,\]
and introduce
the function
\begin{equation}
\varphi(c) : = \exp\left\{ -(c+1)e^{-c} \right\}. \label{eq:varphi definition} 
\end{equation}
\begin{lem} \label{lem:Phi} If $k$ which satisfies (\ref{equ: k used in this range}) and $ r = o(1/p)$ then we have
\[ \Phi \sim \varphi(c(m,k,r))^{n-k-r} \sim \left(1- (1-p)^{k-r} - (k-r)p( 1 - p)^{k-r-1} \right)^{n-k-r}. \]
\end{lem}
\noindent The proof of Lemma~\ref{lem:Phi} is given in Section~\ref{sec:poisson}. 

 Define
\begin{equation}
    X'(m,k,r) := N(k,r) U(m,k,r) \varphi( c(m,k,r) )^{n-k-r}. \label{equ: X' definition}
\end{equation}
Note that Lemma~\ref{lem:Phi} and (\ref{equ: expectation of X}) imply 
\begin{equation} \label{equ: X' asymptotic formula}
    \EE[ X_{k,r}] \sim X'(m,k,r).
\end{equation}
Throughout this work we make use of estimates for the changes in the factors that comprise the expected value of $ \EE[X_{k,r}]$ that
result from changes in the parameters $k$ and $r$.
\begin{lem} \label{lem: many ratios} The following equations hold for $k$ that satisfy (\ref{equ: k used in this range}) and
all $r = O( (\log n)^3/ (np^2))$; moreover, (\ref{equ: N ratio over r}) and (\ref{equ: U ratio over r})  hold for all $k$ that satisfy (\ref{equ: k used in this range}) and all $r < k$:

\begin{align}
    \frac{N(k+1,r)}{N(k,r)} &= (1\pm o(1)) \frac{n}{k} \label{equ: N ratio over k} \\
    \frac{N(k,r+1)}{N(k,r)} &= \frac{(n-k-r)(k-r)}{2(r+1)} \label{equ: N ratio over r}\\
    \frac{U(m,k+1,r)}{U(m,k,r)} &= (1 \pm o(1))\frac{k^2}{e^2n^2} \label{equ: U ratio over k}    \\
    \frac{U(m,k,r+1)}{U(m,k,r)} &= (1 \pm O(p^2k+ p k^3 n^{-2})) e^{(3r-k)p}p \label{equ: U ratio over r} \\
    \frac{\varphi( c(m,k+1,r))^{n-k-r-1} }{\varphi( c(m,k,r))^{n-k-r}} &= 1 \pm O(p^2k^3 n^{-1}) \label{equ: phi ratio over k} \\
    \frac{\varphi( c(m,k,r+1))^{n-k-r-1}}{\varphi( c(m,k,r))^{n-k-r}} &= 1 \pm O(p^2k^3n^{-1}) \label{equ: phi ratio over r} 
\end{align}
\end{lem}

\begin{proof} 

We leave the verification of (\ref{equ: N ratio over k}) and (\ref{equ: N ratio over r}) to the reader, note that $ (\log n)^3/( np^2) \ll k \ll n$
and note that (\ref{equ: N ratio over r}) holds for all $r < k$, as required. 

Next we prove (\ref{equ: U ratio over k}) and (\ref{equ: U ratio over r}) together. 
The following observation will be used for verifying both equations: for any $a$ (positive {\it or} negative) such that 
$a = O(k)$ (note that, here, $r$ {\it could} be as large as $k-1$), we have 
\begin{align}
    \frac{\binom{ \binom{n}{2} - M_1 + a}{m-r}}{\binom{ \binom{n}{2} - M_1 }{m-r}} = \frac{[ \binom{n}{2} - M_1 + a]_{m-r}}{[ \binom{n}{2} - M_1]_{m-r}} &= \prod_{i=0}^{m-r-1} \left(1 + \frac{a}{\binom{n}{2} - M_1 -i} \right) \nonumber \\&= \exp\left\{\frac{a (m-r)}{\binom{n}{2}}  \pm O\left(m \left(\frac{a(m + k^2)}{n^4}+\frac{a^2}{n^4}\right)\right)\right\} \nonumber \\&=
    \exp\left\{a p  - \frac{ar}{\binom{n}{2}} \pm O\left(\frac{a m (m+k^2)}{n^4}\right)\right\} \nonumber \\&=
    \exp\{ap \pm O\left( kp^2 +  k^3 p n^{-2}\right)\}. \label{equ: E ratio grunt work}
\end{align}
As $ M_1(k+1,r) = M_1(k,r)+ k+r$, it follows that for $r \leq r_1$ we have:
\begin{equation*}
    \frac{U(m,k+1,r)}{U(m,k,r)} = \exp\{(-k-r)p \pm O(kp^2 + k^3 p n^{-2})\} \sim e^{-kp} \sim \frac{k^2}{e^2n^2}.
\end{equation*}

Equation (\ref{equ: U ratio over r}) follows similarly, using $a = M_1(k,r+1) - M_1(k,r) = k-3r$, except here there is one additional ratio needed
(note that, again, we must allow $r$ to be as large as $k-1$):
\begin{align*}
        \frac{\binom{ \binom{n}{2} - M_1(k,r+1) }{m-r-1}}{ \binom{\binom{n}{2} - M_1(k,r)}{m-r}} &= \frac{\binom{ \binom{n}{2} - M_1(k,r+1) }{m-r-1}}{\binom{\binom{n}{2} - M_1(k,r)}{m-r-1}} \cdot\frac{\binom{\binom{n}{2} - M_1(k,r)}{m-r-1}}{\binom{\binom{n}{2} - M_1(k,r)}{m-r}} \nonumber \\&=
        \exp\left\{(3r-k)p \pm O\left( kp^2 + k^3 p n^{-2} \right)\right\} \left(\frac{m-r}{ \binom{n}{2} - M_1(k,r) - m +r +1}\right) \nonumber \\&=
        \exp\left\{(3r-k)p \pm O\left( kp^2 + k^3 p n^{-2} \right)\right\} \left(1 \pm O\left(\frac{k}{m}  +\frac{k^2 + m}{n^2} \right)\right)\left(\frac{m}{\binom{n}{2}}\right) \nonumber \\&=
        \left(1 \pm O\left( kp^2 +  k^3pn^{-2} + k m^{-1} + k^2 n^{-2}+p\right)\right)e^{(3r-k) p} p \nonumber \\&=
    (1 \pm O( kp^2 + k^3 p n^{-2})) e^{(3r-k)p} p.
\end{align*}

Verifying (\ref{equ: phi ratio over k}) and (\ref{equ: phi ratio over r}) 
involves estimating  the change of the function $\varphi(c) $ with a slight change in $c$. First, note that for these
equations we assume $ r = O( pk^3n^{-1} )$ and therefore $ c = pk \pm o(1)$. In particular, we may assume that
$c$ is large and $ \varphi(c) = 1 - O( pk^3n^{-2})$. So, we may take the exponents in the expressions in (\ref{equ: phi ratio over k}) and (\ref{equ: phi ratio over r}) to be $ n-k-r$ (instead of $n-k-r-1$). Next, note that $\frac{d}{dc} \log(\varphi(c)) = c e^{-c}$, hence for $\Delta c$ at most a constant we have
\begin{equation*}
    \frac{\varphi(c+\Delta c)}{\varphi(c)} = \exp\left\{O(\Delta c \cdot 
    c e^{-c})\right\} = \exp \left\{O(\Delta c \cdot pk^3n^{-2})\right\} .
\end{equation*}
Setting $ \Delta c = p $ gives then gives the estimates in (\ref{equ: phi ratio over k}) and (\ref{equ: phi ratio over r}).

\end{proof}

We conclude this subsection with two simple Lemmas that allow us to pass between
different forms of our approximations for $ \varphi(c)^n$.
\begin{lem} \label{lem:cs}
Suppose $ c_1 = c_1(n)$ and $ c_2= c_2(n)$ satisfy
\[ c_1 , c_2 = 2\left( \log(np) - \log\log( np ) + \log(e/2) + o(1) \right) \ \ \
\text{ and } \ \ \ c_1-c_2 = o( (\log n)^{-3} n p^2). \] 
If $  t(n) \le n$ then we have
\[  \varphi(c_1)^{t(n)} \sim \varphi( c_2)^{t(n)}. \]
\end{lem}
\begin{proof}
Let $ f(c) = (c+1) e^{-c}$.
Note that $\frac{d f}{dc} =  - c e^{-c} $. Thus,
\[ t(n)| f(c_1) - f(c_2)| 
\le n ( c_1 - c_2) \cdot O \left( \frac{ (\log n)^3}{ (np)^2}   \right) = o(1),
\]
as desired.
\end{proof}
\noindent
For any probability $q$ define
\[ c_q = c_q(k,r) = q(k-r).\]
Note that using this notation we have $ c(m,k,r) = c_p(k,r)$.
\begin{lem} \label{lem:p to c}
If $ q = p(1 \pm o( (\log n)^{-1}))  $, $k$ satisfies (\ref{equ: k used in this range}), $ r = o(1/p)$ and $ t(n) \le n $ then we have
\[ \left( 1 - (1-q)^{k-r} - q (k-r) (1- q)^{k-r-1} \right)^{t(n)} \sim \varphi( c_q(k,r))^{t(n)}.   \]
\end{lem}
\begin{proof}
We have
\begin{align*} 
1 - (1-q)^{k-r} & - q (k-r) (1- q)^{k-r-1} = \exp\left\{ -\left(\frac{ c_q}{ 1-q} + 1\right) (1- q)^{k-r} + O\left( \left( c_q (1 - q)^{k-r} \right)^2  \right)   \right\} \\
& = \exp\left\{ -\left(c_q+ 1\right) (1- q)^{k-r} + O \left( q c_q ( 1 -q)^{k-r} \right) + O\left( \left( \frac{ (\log n)^3}{ (np)^2} \right)^2  \right)   \right\} \\
& = \exp\left\{ -\left(c_q+ 1\right) \exp\left\{ - q(k-r) + O( q^2k) \right\} + O \left(  \frac{ (\log n)^3 }{ n^2 p} \right) + O\left(  \frac{ (\log n)^6}{ n^4p^4}   \right)   \right\} \\
& = \exp\left\{ -\left(c_q+ 1\right) e^{-c_q} + O\left( \frac{ (\log n)^3}{ p^2n^2} \cdot p^2 k  \right) + O \left(  \frac{ (\log n)^3 }{ n^2 p} \right) + O\left(  \frac{ (\log n)^6}{ n^4p^4}   \right)   \right\}.
\end{align*}
The desired estimate follows (using the assumption $ p \ge n^{-3/4 + \epsilon}$).
\end{proof}

\subsection{Estimates for $k_0$ and $r_0(k)$}

Recall that we restrict our attention to $k$ that satisfy (\ref{equ: k used in this range}) and that we define $ r_0(k)$ to
the value of $r \leq r_1$ that maximizes $ \EE[ X_{k,r}]$. Furthermore $k_0$ is defined the be the smallest 
$k$ for which $ \EE[ X_{k, r_0(k)}] \le n^{\epsilon}$. In this section we estimates these two quantities.

We begin by noting two estimates for ratios of the variable $X'$ that follow
immediately from Lemma~\ref{lem: many ratios}. If $ k$ satisfies (\ref{equ: k used in this range}) and $ r = O( (\log n)^3/(np^2))$ then we have
\begin{align}
       \frac{X'(m,k+1,r)}{X'(m,k,r)} &\sim \frac{k}{e^2 n} \label{equ: X' ratio over k}, \ \text{ and} \\
       \frac{X'(m,k,r+1)}{X'(m,k,r)}  & \sim \frac{k^3 p}{2e^2 n (r+1)}.  \label{equ: X' ratio over r} 
       \end{align}
The approximation for $r_0$ can be quickly deduced by using (\ref{equ: X' ratio over r}):
\begin{lem} \label{lem: r0 approx}
    If $k$ satisfies (\ref{equ: k used in this range}) then $r_0(k) \sim \frac{k^3}{2e^2n}$.
\end{lem}

\begin{proof}
    Equation (\ref{equ: X' ratio over r}) implies that there is some function $f=f(n)$ such that
    \[ \frac{X'(m,k,r+1)}{X'(m,k,r)}  = f \cdot\frac{k^3 p}{2e^2 n (r+1)} \]
    and $ f \sim 1$. If $r+1 < \frac{f k^3p}{2e^2 n}$ then $X'$ is increasing in $r$, and if $ r+1> \frac{f k^3p}{2e^2 n}$ then $X'$ is decreasing in $r$. Hence, over $r \leq r_1$, $X'$ is at a maximum for some $r_0 \sim \frac{k^3p}{2e^2n}$.
\end{proof}

Note that with estimate in hand we can conclude that if $k$ satisfies (\ref{equ: k used in this range}) and $ r = o(1/p)$ then we have
\begin{equation}
\label{eq:standard varphi}
  ( \varphi( c(m,k,r)))^{n-k-r} = e^{-2r_0(1+ o(1))}.
\end{equation}

Now we certify that $k_0$ satisfies (\ref{equ: k used in this range}) and approximate the additive difference between $k_0$ and $k_V$.
\begin{lem} \label{lem: k_v minus k_0}
    $k_V - k_0 \sim \frac{k_V^2}{e^2 n}$; in particular, $k_0$ satisfies (\ref{equ: k used in this range}) and
    $k_V - k_0 = \omega(1)$.
\end{lem}
\begin{proof}
First, we use (\ref{equ: X' ratio over k}) and (\ref{equ: X' ratio over r}) 
along with Lemma~\ref{lem: r0 approx} to relate $ X'(m,k_0, r_0)$ and $ X'(m,k_V,0)$.
\begin{align*}
    \frac{X'(m,k_0,r_0)}{X'(m,k_V,0)} &= \left(\frac{X'(m,k_0,r_0)}{X'(m,k_0,0)} \right) \left(\frac{X'(m,k_0,0)}{X'(m,k_V,0)}\right) \\&=
    \left(\prod_{r=0}^{r_0-1} \left(\frac{X'(m,k_0,r+1)}{X'(m,k_0,r)} \right) \right)\left(\prod_{k=k_0}^{k_V - 1} \left(\frac{X'(m,k,0)}{X'(m,k+1,0)} \right)\right) \\&=
    \left(\frac{((1+o(1))r_0)^{r_0}}{r_0!}\right) \exp\left\{(\ln(n/k) + 2 +o(1))(k_V - k_0)\right\} \\&=
    \exp\{(1+o(1)) (r_0 + (k_V p /2) (k_V - k_0) )\}.
\end{align*}
Applying (\ref{eq:standard varphi}) we have
\begin{equation*}
    \frac{X'(m,k_0,r_0)}{N(k_V,0) U(m,k_V,0)} = \exp\{(1+o(1)) (-r_0 + (k_V p /2) (k_V - k_0) )\}.
\end{equation*}
Observe that we have $  n^{-1} < X'(m, k_0, r_0) < 2n^{\epsilon} $ by the definition of $ k_0$ and (\ref{equ: X' ratio over k}). Note that
$ N(k_V,0) U(m,k_V,0)$ is the expected number of independent sets of size $ k_V$. So we
have $  n^{-1} < N(k_V,0) U(m,k_V,0) \leq n^{\epsilon} $ by the definition of $ k_V$, (\ref{equ: N ratio over k}) and (\ref{equ: U ratio over k}). 
Therefore,
\[ \left| r_0 - k_V p ( k_V - k_0)/2 \right| = O( \log n).\]
Noting that $ r_0 = \omega( \log n)$, it follows that we have
\begin{equation*}
    k_V - k_0 \sim \frac{2 r_0}{k_V p} \sim \frac{k_V^2}{e^2 n}.
\end{equation*}
Recalling Lemma~\ref{lem: k vanilla}, we see that this implies that $ k_V - k_0 = \omega(1)$, $ k_V - k_0 = o(1/p)$, and 
$k_0$ satisfies (\ref{equ: k used in this range}).    
\end{proof}

\section{First moment: Proof of Theorem~\ref{thm:techy2}}

Recall that we set $r_1 = 8 \log(np)^3/(n p^2)$. Our goal in the Section is to prove that if $ k$ satisfies (\ref{equ: k used in this range}) then we have
\[ \sum_{r=0}^{ r_1 } \EE[ X_{k, r_0(k)}] = \Theta( \sqrt{ r_0(k)} \cdot \EE[X_{k,r_0(k)}]). \]
We begin by making an observation regarding the log-concavity in $r$ of the sequence $ X'(m,k,r)$ for
small $r$.
\begin{lem}
If $k$ satisfies (\ref{equ: k used in this range}) and $r \leq r_1$, then
    \begin{equation}
        \frac{X'(m,k,r+2) \cdot X'(m,k,r)}{(X'(m,k,r+1))^2} = 1 - \frac{1}{r+2} + o(r^{-1}). \label{equ: X' double ratio over r}
    \end{equation}
\end{lem}
\begin{proof}
We consider the three factors that comprise $X'$ in turn.
By two applications of (\ref{equ: N ratio over r}) we have
\begin{align*}
    \frac{N(k,r+2) \cdot N(k,r)}{(N(k,r+1))^2} &= \frac{(n-k-r-1)(k-r-1)(r+1)}{(n-k-r)(k-r)(r+2)} \\&= 
    1 - \frac{1}{r+2} - O(n^{-1} + k^{-1}) \\&=
    1 - \frac{1}{r+2} - o(r^{-1}).
    \end{align*}
    Similarly, two applications of (\ref{equ: U ratio over r}) 
    give
\begin{align*}    
    \frac{U(m,k,r+2) \cdot U(m,k,r)}{(U(m,k,r+1))^2} &= (1 \pm O(kp^2 + k^3 p n^{-2})) e^{3p} \\&=
    1 \pm O(kp^2 + k^3 p n^{-2}) \\&=
    1 + o(r^{-1}).
\end{align*}
It remains to consider the analogous ratios of the function $ \varphi(c)^{n-k-r}$. Unfortunately the error term in (\ref{equ: phi ratio over r}) is too large for us to use (\ref{equ: phi ratio over r}) directly, so we must resort to consider the second derivative. Set $ f =\log(\varphi(c))$. We have $\frac{d^2 f}{d c^2} = -(c-1) e^{-c}$, hence setting $x = (k-r)p,$ we have
\begin{align*}
    & \frac{\varphi( c(m,k,r+2))^{n-k-r-2}  \cdot \varphi(c(m,k,r))^{n-k-r}}{\left[ \varphi(c(m,k,r+1))^{n-k-r-1} \right]^2} \\
    & \hskip1.5cm =  \frac{\varphi( c(m,k,r+1))^2}{\varphi( c(m,k,r+2))^2} \cdot \exp \left\{ (n-k-r)[ f( x) + f(x-2p) - 2 f(x-p)]\right\} \\& \hskip1.5cm =
    \exp\{O( pk^3 n^{-2}) - O(p^2 \cdot n((k-r)p-1)e^{-(k-r)p})\} \\& \hskip1.5cm =
    1 + o(r^{-1}).
\end{align*}
The Lemma follows from taking the product of these three estimates.
\end{proof}

For ease of notation write $ r_0$ for $ r_0(k)$ for the remainder of this Section.
Now we write $ X'(m,k,r)$ in terms of $ X'(m,k,r_0)$ using a telescoping product. Note that by definition of $r_0$, we have 
\begin{equation*}
    \frac{X'(m,k,r_0)}{X'(m,k,r_0-1)} \geq 1 \geq \frac{X'(m,k,r_0 + 1)}{X'(m,k,r_0)},
\end{equation*}
hence by (\ref{equ: X' double ratio over r}) we have
\begin{equation}
    \frac{X'(m,k,r_0 + 1)}{X'(m,k,r_0)} = 1 - O(r_0^{-1}).\label{equ: X' ratio at r-OPT}
\end{equation}
Define $r_{low} := r_0 - r_0^{3/4}$ and $r_{high} := r_0 + r_0^{3/4}$. For any $r \in [r_0, r_{high}]$, by (\ref{equ: X' double ratio over r}) and (\ref{equ: X' ratio at r-OPT}) we have
\begin{align}
    \frac{X'(m,k,r)}{X'(m,k,r_0)} &= \prod_{r'=r_0}^{r-1} \frac{X'(m,k,r'+1)}{X'(m,k,r')} \nonumber \\&=
    \left(\prod_{r'=r_0+1}^{r-1} \left( \prod_{r'' = r_0}^{r'-1} \frac{X'(m,k,r''+2) \cdot X'(m,k,r'')}{(X'(m,k,r''+1))^2} \right)\right) \nonumber \\ & \qquad \times \left(\frac{X'(m,k,r_0+1)}{X'(m,k,r_0)}\right)^{r-r_0} \nonumber \\&=
    \exp\left\{-(1+o(1))\left(\binom{r-r_0}{2}\frac{1}{r_0}\right) - O(r_0^{-1}(r-r_0))\right\} \nonumber \\&=
    \exp\left\{-(1+o(1)) \frac{(r-r_0)^2}{2r_0} + o(1)\right\}.  
\end{align}
With similar calculations, one can also prove that the above also holds with $r \in {[r_{low}, r_0]}$, i.e:
\begin{align}
    \frac{X'(m,k,r)}{X'(m,k,r_0)} = \exp\left\{ - (1+o(1))\frac{(r-r_0)^2}{2r_0}+o(1)\right\} \qquad \text{for all $r \in [r_{low}, r_{high}].$}\label{equ: quadratic in exponent}
\end{align}

Equation (\ref{equ: quadratic in exponent}) establishes that $
\sum_{r=r_{low}}^{r_{high}} X'(m,k,r)$ approximates a Riemann sum of a dilated normal distribution with the domain size of a higher order than the standard deviation, therefore
\begin{equation*}
    \frac{1}{\sqrt{2 \pi r_0}}\sum_{r=r_{low}}^{r_{high}} X'(m,k,r) \sim  X'(m,k,r_0).
\end{equation*}

Finally, we need to show that the values of $r$ outside of $[r_{low}, r_{high}]$ contribute a negligible amount 
to the sum of $X'(m,k,r)$ over all $r \leq r_{1}$. Equation (\ref{equ: X' double ratio over r}) implies that $ X'(m,k,r) $ is increasing for $r < r_0$ and decreasing for $r > r_0$. Therefore, by (\ref{equ: quadratic in exponent}), we have
\begin{equation*}
  0 \le r < r_{low} \ \ \ \Rightarrow \ \ \  X'(m,k,r) \leq X'(m,k,r_{low}) \leq e^{-r_0^{1/2}/3}X(m,k,r_0),
\end{equation*}
and the analogous bound holds for $ r_{high} < r < r_1$.
Since $r_1 e^{-r_0^{1/2}/3} = o(1)$, the sum of $X(m,k,r)$ over all $r \in [0,r_1] \backslash [r_{low}, r_{high}]$ is $o(X'(m,k,r_0))$.

\section{Second Moment: Proof of Theorem~\ref{thm:techy1}}

In this Section we consider $k,r$ such that $k$  satisfies (\ref{equ: k used in this range}), $r = O( ( \log n)^3 p^{-2}n^{-1})$ and we 
have $ \EE[ X_{k,r}] = \omega(1) $ and show that $ \Var( X_{k, r}) = o ( \EE[ X_{k, r}]^2)$. We emphasize that we restrict our attention
to $ n^{5/4 + \epsilon} < m< n^{4/3 + \epsilon}$.

Given a fixed pair $A, M_A$ where $A$ is a set of $k$ vertices and $ M_A \subset \binom{A}{2}$ is a matching consisting of $r$ edges let
$ \cX_A$ be the event that the pair $A, M_A$ is counted by $ X_{k,r}$.
We have 
\[ \frac{ \Var[ X_{k,r}]}{ \EE[ X_{k,r}]^2}
= \frac{1}{ \EE[ X_{k,r}]^2} \sum_{A, M_A; B, M_B} \PP( \cX_A \wedge \cX_B) - \PP( \cX_A) \PP( \cX_B) .
\]
As is common for second moment calculations of this kind, we break this sum down depending on the intersection patterns of $A, M_A$ and $ B, M_B$. Define
\begin{equation} \label{eq:h'''}
h_i  :=  \frac{1}{ \EE[ X_{k,r}]^2} \sum_{ \stackrel{ A, M_A; B ,M_B}{ |A \cap B|= i} } \PP ( \cX_A \wedge \cX_B) -  \PP( \cX_A) \PP( \cX_B).
\end{equation}
We break the variance calculations into two main cases.

\subsection{ $i \le 1/(2p)$ }
In this case we need to use the full covariance term. Let $ A, M_A$ be fixed. We have
\[ 
h_i  =  \frac{1}{ | {\mathbb N} |} \sum_{ B ,M_B :|A \cap B|= i}  \frac{ \PP ( \cX_A \wedge \cX_B) }{ \PP( \cX_A) \PP( \cX_B)  } - 1, \]
where $ {\mathbb N}$ is the collection of all pairs $ B, M_B$.
Thus, the main work in this Section is in bounding $ \PP ( \cX_A \wedge \cX_B
) $ relative to $ \PP( \cX_A)^2$. 

We begin with our estimate for $ \PP( \cX_A)$. 
Let $ \cU_A$ be the event that all edges within $A$ appear
in accordance with the event $ \cX_A$.
Recall that $ \Phi $ is the probability that all vertices outside of $A$ have the required
two neighbors in $A \setminus ( \cup_{e \in M_A} e)$ when we condition on the event $ \cU_A$. 
So, we have 
\[ \PP( \cX_A ) = \PP( \cU_A) \Phi = \frac{ \binom{ \binom{n}{2} - \binom{k+r}{2}+ \binom{2r}{2} - r }{m-r}}{ \binom{ \binom{n}{2}}{m}}\Phi. \] 

In order to write an expression for an upper bound of $ \PP( \cX_A \wedge \cX_B)$, we need to 
introduce some parameters.
Given fixed pairs $ A, M_A$ and $B, M_B$ where $A, B$ are sets of $k$
vertices and $ M_A \subset \binom{A}{2}, M_B \subset \binom{B}{2} $ are matchings
that consist of $r$ edges, define
\begin{gather*}
 i := |A \cap B| \ \ \ \ \ \ \ell := |M_A \cap M_B| \ \ \ \ \ \  s := \left|  \left( \cup_{e \in M_A} e \right)  \cap \left( \cup_{e \in M_B} e \right) \right| \\
 a := \left| \left( \cup_{e \in M_A} e \right) \cap B \right| \ \ \ \text{ and }  \ \ \ b := \left| \left(  \cup_{e \in M_b} e \right) \cap A \right|.
\end{gather*}
In analogy with $\Phi$, we let $\Psi$  be the probability that all vertices outside
$A$ have the required two neighbors in $A \setminus ( \cup_{ e \in M_A}e )$ and all vertices outside $B$ have the required two neighbors in $B \setminus( \cup_{ e \in M_B}e ) $ when we
condition on the event  $ \cU_A \wedge \cU_B$. We have 
\[ \PP( \cX_A \wedge \cX_B) = \PP( \cU_A \wedge \cU_B) \Psi=  
\frac{  \binom{ \binom{n}{2} - 2\binom{k+r}{2} + \binom{ i}{2} + 2\binom{2r}{2} - \binom{s}{2} - 2r + \ell}{m - 2r + \ell} }
  { \binom{ \binom{n}{2}}{m}} \Psi.  \]
We note in passing that a precise expression for $ \Psi$ would depend on all of the parameters that we introduce here, but the estimate that we use is a function of only $ n,m,i,a$ and $b$. We suppress any dependence on these parameters in $\Psi$ in order to reduce clutter. (E.g. We could write $ \Psi_{n,m,i,a,b}$ in the place of $ \Psi$, but we refrain from doing this.) We also note that 
the parameter $s$ does not have a meaningful impact on the calculations that follows.

We seek to bound
\[ \frac{ \PP( \cX_A \wedge \cX_B )}{ \PP (\cX_A)^2} -1 = \frac{ \PP(\cU_A \wedge \cU_B)}{ \PP( \cU_A)^2} \cdot \frac{ \Psi}{ \Phi^2} -1.\]
So, we require an estimate for the ratio $ \Psi/ \Phi^2$. Note that if we condition on the event $ \cU_A \wedge \cU_B$, which fixes most of the edges in $ \binom{A}{2} \cup \binom{B}{2}$, then the 
event $ \cX_A \wedge \cX_B$ includes requirements on the edges between $ A \setminus B$ and $ B \setminus A$. In particular, if we consider the bipartite graph induced by these two sets of vertices then there is a minimum degree requirement on the vertices. This introduces delicate dependencies as the edges of this bipartite graph play a role in satisfying the minimum requirements for both $A$ and $B$. We resolve this issue by apply a variant on Janson's inequality that was recently introduced by
Bohman, Warnke and Zhu. In Section~\ref{sec:Psi} we use this variant to establish the following estimate.
\begin{lem} \label{lem:Psi}
Consider $k$ that satisfies (\ref{equ: k used in this range}) and
let $A, M_A$ and $ B,M_B$ be pairs in $ {\mathbb N}$. We set
$i = |A \cap B|$,  $ a = |B \cap ( \cup_{e \in M_A}e )|$ and $ b = |A \cap( \cup_{e \in M_B} e)|$. If $ i = O(1/p)$ then we have 
\[ \Psi \le (1 +o(1)) \Phi^2 \exp \left\{ O \left( i^2 \cdot \frac{ p^3 k^3}{n} \right) + O \left( (a+b) \cdot \frac{ k^3 p^2 }{n} \right) 
 \right\}.\]
\end{lem}

We bound the variance by summing the contributions for specific
choices of $i$ and $\ell$. For ease of notation set
\begin{align*}
M_1 & = \binom{k+r}{2} - \binom{2r}{2} + r  \ \ \ \text { and }\\
M_2 & = 2\binom{k+r}{2} - \binom{ i}{2} - 2\binom{2r}{2} + \binom{s}{2} + 2r - \ell.
    \end{align*}
Let $ h_{i,\ell} $ be the contribution to the variance coming from pairs
$ A, M_A; B, M_B$ such that $ |A \cap B| = i$ and $ |M_A \cap M_B| = \ell$. We have $ h_i = \sum_{\ell = 0}^r h_{i,\ell}$ where
\begin{equation}  
\begin{split} 
\label{eq:hiell}
h_{i, \ell} 
& = \frac{1}{ | {\mathbb N}| } \sum_{B, M_B}  \left(     
\frac{ \binom{ \binom{n}{2} - M_2}{m-2r+\ell} \binom{ \binom{n}{2}}{m} }{\binom{ \binom{n}{2} - M_1}{m-r}^2} \cdot \frac{ \Psi}{ \Phi^2}
- 1 \right),
    \end{split}                                       
 \end{equation}  
and the sum is over all pairs $B,M_B$ such that $|A \cap B|= i$ and $ |M_A \cap M_B |=\ell$.

Before we proceed to establish bounds on $h_{i,\ell}$ for particular choices of
the parameters $i$ and $\ell$, we carefully estimate $ \PP( \cU_A \wedge \cU_B) \PP( \cU_B)^{-2} $. We have 
\begin{equation*}
    \begin{split}
        & \frac{ \PP ( \cU_A \wedge \cU_B) }{ \PP( \cU_A) \PP ( \cU_B)}  = \frac{\binom{ \binom{n}{2} - M_2
        }{m - 2r + \ell} \binom{ \binom{n}{2}}{m}}{ \binom{ \binom{n}{2} - M_1 
        }{m-r} ^2 } 
        = \frac{ [m-r]_{r - \ell}}{ [m]_r} \cdot
        \frac{ \left[  \binom{n}{2} - M_2 
        \right]_{m - 2r + \ell}  \left[ \binom{n}{2} \right]_{m}}{ \left[ \binom{n}{2} - M_1  
        \right]_{m-r}^2} \\
    & \hskip5mm =  \frac{ [m-r]_{r - \ell}}{ [m]_r} \cdot \frac{  \left[ \binom{n}{2} - m +r \right]_{r}}{\left[ \binom{n}{2} - M_2 - m + 2r - \ell
    \right]_{r - \ell}}
    \cdot
        \frac{ \left[  \binom{n}{2} - M_2 
        \right]_{m - r}  \left[ \binom{n}{2} \right]_{m-r}}{ \left[ \binom{n}{2} - M_1
        \right]_{m-r}^2} \\
        & \hskip5mm =  \left[ \frac{ \binom{n}{2}}{m} \right]^\ell
        \cdot \left[\prod_{x=0}^{\ell-1} \frac{1}{1 - x/m} \left(1 - \frac{ m -r - x}{ \binom{n}{2}} \right) \right]\\
        & \hskip15mm \cdot \left[\prod_{x=0}^{r - \ell-1} \left( 1  - \frac{r-\ell}{ m - \ell -x}\right) \left( 1+
        \frac{ M_2 - r  }{ \binom{n}{2} - M_2 - m +2r - \ell -x } \right) 
        \right]\cdot
        \frac{ \left[  \binom{n}{2} - M_2 
        \right]_{m - r}  \left[ \binom{n}{2} \right]_{m-r}}{ \left[ \binom{n}{2} - M_1
        \right]_{m-r}^2} \\
&  \hskip5mm \le \left[ \frac{ \binom{n}{2}}{m} \right]^\ell \cdot \frac{ \left[  \binom{n}{2} - M_2 
        \right]_{m - r}  \left[ \binom{n}{2} \right]_{m-r}}{ \left[ \binom{n}{2} - M_1
        \right]_{m-r}^2} \cdot \exp \left\{ O(\ell^2/m) + O( r k^2/n^2 )  \right\} \\
        & \hskip5mm \le \left[ \frac{ \binom{n}{2}}{m} \right]^\ell \cdot \frac{ \left[  \binom{n}{2} - M_2 
        \right]_{m - r}  \left[ \binom{n}{2} \right]_{m-r}}{ \left[ \binom{n}{2} - M_1
        \right]_{m-r}^2} \cdot  \left( 1 + \frac{1}{ n^{ \Omega( \epsilon)}} \right).
    \end{split}
\end{equation*}
Note that we use $k \le n^{3/4- \epsilon}$. Next, we observe that for $ 0 \leq x < m-r$ we have 
\begin{multline*}
\Xi_x:= \frac{ \left(  \binom{n}{2} -  M_2 - x \right)  \left( \binom{n}{2} -x \right)}{ \left( \binom{n}{2} - M_1 - x \right)^2}
= 1 +
 \frac{  \left( \binom{i}{2} - \binom{s}{2} + \ell \right)  \left( \binom{n}{2} - x \right) - \left( \binom{k+r}{2} - \binom{2r}{2} + r \right)^2  } 
   { \left( \binom{n}{2} - \binom{k+r}{2} + \binom{2r}{2} - r - x \right)^2 } \\
   = 1 + \frac{ \binom{i}{2} \binom{n}{2} - \binom{k+r}{2}^2 + 2 \binom{2r}{2} \binom{k+r}{2} -x \binom{i}{2} - \left( \binom{s}{2} - \ell \right) \left( \binom{n}{2} - x \right) - 2r \binom{k+r}{2} - \left( \binom{2r}{2} - r \right)^2}{\left( \binom{n}{2} - \binom{k+r}{2} + \binom{2r}{2} - r - x \right)^2 } \\
   \le 1 + \frac{ \binom{i}{2} \binom{n}{2} - \binom{k+r}{2}^2  + O( p^2 k^8/n^2) }{ n^4/ 5}
=  1 + \frac{ 5\left( \binom{i}{2} \binom{n}{2} - \binom{k+r}{2}^2 \right)}{ n^4} 
+ O( p^2 k^8/n^6). 
\end{multline*}

As our bound does not depend on $x$ we set
\begin{equation}
    \label{eq:Xi}
\Xi:= \max_x \Xi_x \le  1 + \frac{ 5\left( \binom{i}{2} \binom{n}{2} - \binom{k+r}{2}^2 \right)}{ n^4} 
+ O( p^2 k^8/n^6).
\end{equation}
We have established
\begin{equation}
\label{eq: Us estimate}
\frac{ \PP ( \cU_A \wedge \cU_B) }{ \PP( \cU_A) \PP ( \cU_B)} 
\le \left( 1 + n^{- \Omega( \epsilon)} \right)
p^{-\ell} \Xi^{m-r}.  
\end{equation}

Next we turn to the term coming from the number of ways to choose $A, M_A$ 
and $ B, M_B$.
Recall that we define $N = N(k,r) = |{\mathbb N}|$ to be the number of choices of $A,M_A$. Let
$N_{i, \ell}$ be the number of choices for $ B, M_B$ such that $ |A \cap B|=i$ and
$M_A \cap M_B = \ell$. We have
\begin{align*}
    N &= \binom{n}{k+r} \frac{ [k+r]_{2r}}{ r! 2^r }, \\
    N_{i, \ell}& \le \binom{r}{\ell} \binom{k+r-2 \ell}{ i - 2 \ell} \binom{n-k-r}{k+r-i} \cdot \frac{ [k+r-2\ell]_{2r- 2\ell} }{ (r - \ell)! 2^{r-\ell}}.
    \end{align*}
Note that we have
\begin{multline}
\label{eq:Niell}
\frac{N_{i,\ell}}{N} \le \frac{ \binom{k+r}{i} \binom{n-k-r}{k+r-i}}{\binom{n}{k+r}} \cdot\binom{r}{\ell} 
\frac{ [i]_{2 \ell}} {[k+r]_{2 \ell}} \cdot \frac{[r]_\ell 2^\ell}{ [k+r]_{2\ell}} \\
\le \frac{ \binom{k+r}{i} \binom{n-k-r}{k+r-i}}{\binom{n}{k+r}} \left( \frac{ 2e r^2 i^2}{ \ell (k+r)^4} \right)^{\ell} \left(\prod_{x=0}^{2 \ell-1} \frac{1}{ 1 -x/(k+r)} \right)\\ \le \frac{ \binom{k+r}{i} \binom{n-k-r}{k+r-i}}{\binom{n}{k+r}} \left( \frac{ 2e r^2 i^2}{ \ell (k+r)^4} \right)^{\ell} \exp \left\{ O( \ell^2/k ) \right\}.
\end{multline}

We now have all of the necessary estimates in hand, and we are prepared to 
proceed to bound the sum of $ h_{i, \ell}$, which is 
defined in (\ref{eq:hiell}). We break the calculation in two subcases.

\noindent
{\bf Case 1.1} $ i \le  e^2 \frac{(k+r)^2}{n} $.\\[-6mm]
\begin{quote}
    We begin with careful estimates for $ \Xi^{m-r}$. Set $ \tk = k+r$. It turns out that $ \Xi$ is less than 1 for
    $i$ less than $ \tk^2/n$, and $ \Xi$ is small for other values of $i$ in the range treated in this case. Define $ f$ by setting $ i = (1+f)(\tk^2/n)$.    
    \begin{claim}
    \label{cl:Xi}
    \begin{enumerate}
\item If $ -1 \leq f \leq 0$  then \[ \Xi^{m-r} < 1- \Omega\left( \frac{ p k^2}{  n} \right). \]
\item If $ 0 < f < e^2 -1  $
then
\[ \Xi^{m-r} \le \exp \left\{ \frac{6 p k^4 f}{n^2}\right\}. \]
        \end{enumerate}        
    \end{claim}
    \begin{proof}
    
    Considering both parts for now, we have $ i = (1+f)\tk^2/n$ where $ -1 \leq f < e^2-1$. Here we have
    \begin{align*} 
    4 \left[ \binom{i}{2} \binom{ n}{2} - \binom{\tk}{2}^2 \right]
    &= \left( \frac{(1+f)\tk^2}{n}\right) \left( \frac{(1+f)\tk^2}{n} -1 \right)n(n-1) - \tk^2 (\tk-1)^2 \\
    &\le (f^2 + 2f)\tk^4 + (-1-f)\tk^2 n + O( \tk^3).  
    \end{align*}

    Now we consider each part of the Claim in turn: if $-1 \leq f \leq 0$, then the last expression will be $-\Omega(\tk^2 n)$, since both $(f^2 + 2f)$ and $(-1-f)$ are non-positive and not equal to 0 at the same time, and since $\tk^4 \gg \tk^2 n$. Hence, by (\ref{eq:Xi}) and the assumption $ k \le n^{3/4 - \epsilon}$ (which will allow the $O(p^2 \tk^8 / n^6)$ term in (\ref{eq:Xi}) to be dominated by $\tk^2 n / n^4$), the first part of the Claim holds. If $0 < f < e^2 - 1$ then the $(-1-f) \tk^2 n$ term will dominate the $O(\tk^3)$ as well as the $O(p^2 \tk^8 / n^6)$ in (\ref{eq:Xi}) (like in the previous case), and the first summand is bounded above by $(e^2 + 1)f \tk^4$, therefore
\begin{equation*}
\Xi^{m-r} \le \exp\left\{ \frac{5m}{4 n^4} \cdot (e^2+1)f \tk^4 \right\} \\
\le \exp \left\{ \frac{6 p k^4 f}{n^2}\right\}.
\end{equation*}

\end{proof}
With these estimates for $ \Xi$ in hand, we are ready to turn to bounding 
$ h_{i, \ell}$. In light of the dependence on $ a +b$ in the bound on $ \Psi/ \Phi^2$ given in Lemma~\ref{lem:Psi}, we treat pairs $ B, M_B$ such that $ a+b$ is large
separately. Let $ \cB = \cB_{i,\ell} $ be the set of pairs $ B, M_B$ such that $ |A \cap B| =i$
and  $ |M_A \cap M_B| = \ell$. We partition $ \cB$ into two parts. Set
\[ \vartheta =(\ell +1) n /( k ( \log n)^4)  \]
and let $ \cB' = \cB'_{i, \ell}$ be the
set of pairs $B,M_B$ in $ \cB$ such that $ a+b \le \vartheta $ and set $ \cB'' = \cB''_{i,\ell} = \cB_{i, \ell} \setminus \cB'_{i,\ell}$. We partition $ \cB''$ further based on $ a +b$; in particular, for each $y > \vartheta$ we set $ \cB_y''$ to be the collection of pairs in $ \cB''$ with $a+b=y$.

We use the Chernoff bound to show that $ |\cB_y''|$ is small. We view $A, M_A$ as
fixed, and choose $ B,M_B$ in a series of choices that we view as random
choices for the purpose of counting. First the 
set $B$ is chosen as a random set of cardinality $ \tk$ and
then the vertex set of $ M_B$ is chosen at random from $B$. We have
\[  \EE[a] = \EE\left[ \left| \left( \cup_{e \in M_A}e \right) \cap B \right| \right] 
= \frac{ 2r \tk}{n} \ll \frac{n}{k (\log n)^4 } \le y. \]
Thus, Theorem~\ref{thm:hyperchernoff} implies $ \PP (a \ge y/2) 
\le e^{- \Omega(y)}$. 
Once $B$ is fixed such that $ |A \cap B|=i$, we choose the vertices of $M_B$ at random. In this selection, the expected size of $b$ is $ 2ri/\tk \ll n/( k (\log n)^4)$. It again follows from Theorem~\ref{thm:hyperchernoff} that we have $ \PP (b \ge y/2) 
\le e^{- \Omega(y)}$.
Applying Lemma~\ref{lem:Psi}, Claim~\ref{cl:Xi}, and (\ref{eq: Us estimate}),
we have
\begin{multline*}
\frac{1}{N} \sum_{ B,M_B \in \cB''} \left(  \frac{ \PP( \cU_A \wedge \cU_B)}{ \PP( \cU_B)^2} \cdot \frac{ \Psi}{ \Phi^2} - 1 \right) \\
\le \sum_{ y > \vartheta} O(1) p^{- \ell}  \exp \left\{  O \left( \frac{pk^4}{n^2} \right) + y \cdot O\left( \frac{ k^3p^2}{n} \right) -\Omega(y)  \right\} = \exp\left\{  - \Omega \left( \frac{n}{ k (\log n)^4}  \right) \right\}. 
\end{multline*}
Note that we use the fact that pairs $ B, M_B$ in $ \cB''$ have $ \ell \le (a+b)(k (\log n)^4)/n$ to absorb the $p^{-\ell}$ into the $\Omega(y)$.

For the remainder of this case we restrict our attention to $ \cB'$. Note that Lemma~\ref{lem:Psi} implies that we may now assume
\begin{equation} \label{eq:Psi for small i}
 \Psi \le (1+o(1))\exp \left\{  O\left( \frac{ 1+ \ell}{ ( \log n)^2 }\right) \right\} \Phi^2.
 \end{equation}
We begin with $ \ell =0$. Set 
\[ i_1 = \tk^2/n + n/((\log n) p k^2) \ \ \  \text{ and }  \ \ \
i_2 = e^2 \tk^2/n.  \]
Applying (\ref{eq: Us estimate}) - (\ref{eq:Psi for small i}) and Claim~\ref{cl:Xi} we have 
\[ \sum_{i=0}^{ i_1 } 
\frac{1}{N} \sum_{ B,M_B \in \cB'_{i, 0}} \left(  \frac{ \PP( \cU_A \wedge \cU_B)}{ \PP( \cU_B)^2} \cdot \frac{ \Psi}{ \Phi^2} - 1 \right) 
\le \sum_{i=0}^{ i_1 }
\frac{ \binom{\tk}{i } \binom{n-\tk}{ \tk - i}  }{ \binom{n}{\tk}}  \left(\frac{7 }{ \log n} +o(1) \right) = o(1).\]
For $ i_1 < i \leq i_2 $  (if any such $i$ exist) we apply the Chernoff bound. Note that we have
\[ i_1 \ge  \tk^2/n + k ( \log n)/ n^{1/2}. \]
Recalling that we define $f$ by setting
$ i = (1+f)\tk^2/n$, an application of Theorem~\ref{thm:hyperchernoff}
gives
\[ \frac{ \binom{\tk}{i } \binom{n-\tk}{ \tk - i}  }{ \binom{n}{\tk}} \le 
\exp \left\{ - \Omega \left( f^2 \tk^2/n   \right)   \right\} 
\le \exp \left\{ - \Omega( (\log n)^2 ) \right\}.\]
So, we have 
\begin{multline}
\sum_{i=i_1}^{ i_2 }  
\frac{1}{N} \sum_{ B,M_B \in \cB'_{i,0}} \left(  \frac{ \PP( \cU_A \wedge \cU_B)}{ \PP( \cU_B)^2} \cdot \frac{ \Psi}{ \Phi^2} - 1 \right)  \\
\label{eq:i2}
\le \sum_{i=i_1}^{ i_2 } \exp \left\{o(1)+ O \left( \frac{f pk^4}{n^2}  \right)   - \Omega \left(  \frac{ f^2 k^2 }{n}  \right)   \right\}  = \sum_{i=i_1}^{ i_2 } \exp \left\{- \Omega(( \log n)^2)\right\}. 
\end{multline}

For larger values of $\ell$ our estimates follow the same argument, taking into account the dependence on $\ell$ that appears in 
(\ref{eq: Us estimate}) - (\ref{eq:Psi for small i}).
Note that 
\[ \frac{ r^2 i^2}{ k^4 p } = O \left( \frac{p^2 k^6}{n^2} \cdot \frac{ k^4}{ n^2} \cdot \frac{1}{k^4 p} \right)   = O \left( \frac{ p k^6}{ n^4 }\right) \le \frac{1}{n^{1/4 + 4\epsilon}}.\]
It follows that for $ i \le i_1$, $ \ell \ge 1 $ and $ B,M_B \in \cB'_{i, \ell}$ we have
\[ \frac{ N_{i,\ell}}{N} \cdot  \frac{ \PP( \cU_A \wedge \cU_B)}{ \PP( \cU_B)^2}  \cdot \frac{ \Psi }{ \Phi^2}   \le  
 \frac{ \binom{\tk}{i } \binom{n-\tk}{ \tk - i}  }{ \binom{n}{\tk}}  \left( \frac{1}{ \ell n^{1/4 + \epsilon}} \right)^\ell \exp \left\{ O \left( \frac{( \log n) p k^3 }{n^{3/2}}  \right) \right\} . 
\]
Therefore,
\[ \sum_{i=0}^{ i_1 }  \sum_{\ell=1}^r
\frac{1}{N} \sum_{ B,M_B \in \cB'_{i, \ell}} \frac{ \PP( \cU_A \wedge \cU_B)}{ \PP( \cU_B)^2} \cdot \frac{ \Psi}{ \Phi^2} 
\le  \sum_{\ell = 1}^{ r} \left( \frac{1}{ \ell n^{1/4}} \right)^\ell = O \left( \frac{1}{ n^{1/4}}  \right).   \]
For the $ \ell > 0$ and $i$ in the range $ i_1 < i <i_2$ the argument is even simpler as incorporating 
the impact of $ \ell > 0$ into (\ref{eq:i2}) decreases
an estimate which is already sufficiently small.
\end{quote}

\noindent {\bf Case 1.2.} $  e^2 (k+r)^2/n \le i \le 1/(2p) $  
\begin{quote}
    In this case the $ \binom{i}{2} \binom{n}{2}$ term is no longer balanced by the $\binom{k+r}{2}^2$ term in the numerator of the fraction in the bound on $ \Xi$ given in (\ref{eq:Xi}) (although the $-\binom{k+r}{2}^2$ term still absorbs the $O(p^2 k^8/n^6)$). So, again applying Lemma~\ref{lem:Psi} and (\ref{eq: Us estimate}), we have
\begin{multline*}
    h_{i,0} \le \frac{ \binom{ \tk}{i} \binom{n-\tk}{\tk-i} }{ \binom{n}{\tk}} 
    \cdot e^{ O \left( i^2 \cdot \frac{ p^3 k^3}{n}    \right) + O \left( (a+b) \frac{ k^3p^2}{n} \right) +o(1)} \cdot e^{ 2m \frac{ i^2}{ n^2} } \le \left( \frac{\tk e}{i}  \right)^i \left( \frac{\tk}{n} \right)^i e^{ pi^2 + o(i)} \\
    =    
    \left[ \frac{ \tk^2 e^{1 + o(1)} }{i n} \cdot e^{p i} \right]^i  \le \left[ e^{-1+o(1)} \right]^i. 
    \end{multline*}

For larger values of $ \ell$ we once again apply (\ref{eq: Us estimate}) and (\ref{eq:Niell}). 
We have
\[\left( \frac{ 2e^{1 +o(1)} r^2 i^2}{ \ell \tk^4 p} \right)^{\ell}  
= \left(  \frac{ \Theta(1) p k^2 i^2}{ \ell n^2} \right)^{\ell}.  \]
If $ \ell < i / ( \log n)^2$ then this expression is $ e^{ o( i) }$. On the other hand, if $ \ell > i/ (\log n)^2$ then this expression is less than $  \left(( \log n)^3 k^2 / n^2 \right)^\ell $. In either case the contribution to variance from the sum over all of these values of $i$ and $ \ell$ is negligible.

\end{quote}

\subsection{ $i > 1/(2p)$}

We bound $h_i$ as defined in (\ref{eq:h'''}). Recall that 
we set $ \tk = k+r$. For simplicity, let $j = \tilde{k}-i$. In order to highlight the role of the constant $1/2$ that defines this case we set $ C = 1/2$ and assume $ i > C/p$. In order to bound $ h_i$ we will sum over the possible sizes of the intersection of $ M_A$ and $M_B$. 
Define
\[ r_O = | M_A \setminus M_B| = |M_B \setminus M_A|.\]
(Of course $ r_0 = r - \ell$ where $ \ell = |M_A \cap M_B|$ was an important parameter in the previous subsection). Let $ {\mathbb N}_{i,r_O}$ be
the collection of pairs $ A, M_A; B, M_B$ such that $ |A \cap B| = i$ and $ |M_A \setminus M_B| = r_O$. We begin by bounding the size of $ {\mathbb N}_{i, r_O}$.

\begin{claim} \label{claim: n_i_r}
    \[ | {\mathbb N}_{i, r_O}| \leq |\mathbb{N}|\left(\frac{4e^2 n \tilde{k}}{j^2}\right)^j (12r)^{r_O}.\]
\end{claim}

\begin{center}
\begin{figure}
\begin{tikzpicture}

\node at (0,8){\large Figure 1: Potential structures in $ M_A \cup M_B$ when $A, M_A$ and $B, M_B$ are counted by $X_{k,r}$. };

\node at (0,-8){\LARGE Key};

\node at (-1.5,-9){Edges in $M_A \backslash M_B$};
\node at (-1.5,-9.5){Edges in $M_B \backslash M_A$};
\node at (-1.38,-10){Edges in $M_A \cap M_B$};

\draw[dotted,  red, thick] (1,-9) -- (2.9,-9);
\draw[dashed, blue] (1,-9.5) -- (2.9,-9.5);
\draw (1,-10) -- (2.9,-10);

\draw[thick] (-3.5,-7.2) -- (3.5,-7.2) (3.5,-7.2) -- (3.5,-10.6) (3.5,-10.6) -- (-3.5,-10.6) (-3.5,-10.6) -- (-3.5,-7.2);

\draw (-1.5,0) circle [radius = 6];
\draw (1.5,0) circle [radius = 6];

\node at (-2, 6.7) {\Huge{$A$}};
\node at (2, 6.7) {\Huge{$B$}};

\node (SS1) at (-6.3, 1.5){};
\node (SS2) at (-6.5, 0){};
\node (SS3) at (-6.7, -0.7){};
\node (SS4) at (-6.3, -2.0){};

\draw[dotted,  red, thick] ($(SS1)$) -- ($(SS2)$);
\draw[dotted,  red, thick] ($(SS3)$) -- ($(SS4)$);

\node (TT1) at (6.6, 0.7){};
\node (TT2) at (6.6, -0.7){};

\draw[dashed, blue] ($(TT1)$) -- ($(TT2)$);

\node (S1) at (-4.4,4.2){};

\node (S2) at (-5.1,2.8){};
\node (ST2) at (-3.05,2){};
\node (S3) at (-5.45,1.5){};
\node (ST3) at (-3.33,1){};
\node (S4) at (-5.6,0){};
\node (ST4) at (-3.4,0){};
\node (S5) at (-5.45,-1.5){};
\node (ST5) at (-3.33,-1){};
\node (S6) at (-5.1,-2.8){};
\node (ST6) at (-3.05,-2){};
\node (S7) at (-4.5,-4){};
\node (ST7) at (-2.5,-3){};

\draw[dotted,  red, thick] ($(S2)$) -- ($(ST2)$);
\draw[dotted,  red, thick] ($(S3)$) -- ($(ST3)$);
\draw[dotted,  red, thick] ($(S4)$) -- ($(ST4)$);
\draw[dotted,  red, thick] ($(S5)$) -- ($(ST5)$);
\draw[dotted,  red, thick] ($(S6)$) -- ($(ST6)$);
\draw[dotted,  red, thick] ($(S7)$) -- ($(ST7)$);

\node (T1) at (3.7,4.7){};

\node (T2) at (4.5,4){};
\node (TS2) at (2.5,3){};
\node (T3) at (5.1,2.8){};
\node (TS3) at (3.05,2){};
\node (T4) at (5.45,1.5){};
\node (TS4) at (3.33,1){};
\node (T5) at (5.6,0){};
\node (TS5) at (3.4,0){};
\node (T6) at (5.45,-1.5){};
\node (TS6) at (3.33,-1){};
\node (T7) at (5.1,-2.8){};
\node (TS7) at (3.05,-2){};
\node (T8) at (4.5,-4){};
\node (TS8) at (2.5,-3){};

\draw[dashed, blue] ($(T2)$) -- ($(TS2)$);
\draw[dashed, blue] ($(T3)$) -- ($(TS3)$);
\draw[dashed, blue] ($(T4)$) -- ($(TS4)$);
\draw[dashed, blue] ($(T5)$) -- ($(TS5)$);
\draw[dashed, blue] ($(T6)$) -- ($(TS6)$);
\draw[dashed, blue] ($(T7)$) -- ($(TS7)$);
\draw[dashed, blue] ($(T8)$) -- ($(TS8)$);

\node (I1) at (-0.3, 5){};

\node (I2) at (-1.5,3.5){};
\node (I3) at (-0.25,2.5){};
\node (I4) at (0.75,4){};

\node (I5) at (-0.8,1.2){};
\node (I6) at (0.2,1.3){};
\node (I7) at (0.4,0.1){};
\node (I8) at (-0.85,0){};

\node (I9) at (1.9,-1){};
\node (I10) at (1.5,-2.1){};

\node (I11) at (-1.5,-2.4){};
\node (I12) at (0.5,-1.9){};
\node (I13) at (-1,-3){};
\node (I14) at (0.8,-3.4){};
\node (I15) at (-1.7,-3.4){};
\node (I16) at (-0.5,-4.5){};
\node (I17) at (0.4,-4.6){};
\node (I18) at (1.5,-3.8){};

\draw[dotted, thick, red] ($(S1)$) -- ($(I1)$);
\draw[dashed, blue] ($(I1)$) -- ($(T1)$);

\draw[dashed, blue] ($(ST2)$) -- ($(I2)$);
\draw[dotted, thick, red] ($(I2)$) -- ($(I3)$);
\draw[dashed, blue] ($(I3)$) -- ($(I4)$);
\draw[dotted, thick, red] ($(I4)$) -- ($(TS2)$);

\draw[dashed, blue] ($(I5)$) -- ($(I6)$);
\draw[dotted, thick, red] ($(I6)$) -- ($(I7)$);
\draw[dashed, blue] ($(I7)$) -- ($(I8)$);
\draw[dotted, thick, red] ($(I8)$) -- ($(I5)$);

\draw[dashed, blue] ($(ST5)$) -- ($(ST6)$);

\draw[dotted, thick, red] ($(TS6)$) -- ($(I9)$);
\draw[dashed, blue] ($(I9)$) -- ($(I10)$);
\draw[dotted, thick, red] ($(I10)$) -- ($(TS8)$);

\draw ($(I11)$) -- ($(I12)$);
\draw ($(I13)$) -- ($(I14)$);
\draw ($(I15)$) -- ($(I16)$);
\draw ($(I17)$) -- ($(I18)$);

\end{tikzpicture}

\vskip5mm

In this example we have $ \ell= 4, r_0 =15, r_A = 3$ and $ r_B = 4$. Note that the graph with edge set 
$ M_A \cup M_B$ consists of isolated edges and an alternating cycle in $ A \cap B$ and alternating paths that begin
and end in $A \triangle B$.

\end{figure}
\end{center}
\begin{proof}

There are \[ \binom{n}{\Tilde{k}} \binom{n-\Tilde{k}}{j} \binom{\Tilde{k}}{j}\] ways to choose $A$ and $B$. Now fix sets $A$ and $B$.

We now consider the ways to choose the matchings $M_A$ and $M_B$ that 
have $|M_A \backslash M_B| = |M_B \backslash M_A| = r_O$ and
are compatible with the event $ \cX_A \wedge \cX_B$. (N.b., the $O$ in $ r_O$ standing for ``outer"; i.e., not in the intersection). These conditions impose a number of conditions on the matching $M_A$ and $M_B$. Figure 1 depicts the different structures that $M_A$ and $M_B$ can form. There may edges in $M_A \cap M_B$; these are disjoint from every other matching edge from $M_A$ and/or $M_B$. There may be edges in $M_A \backslash M_B$ disjoint from any $M_B$ edges; such edges have at least one vertex in $A \backslash B$, and vice versa for $M_B \backslash M_A$. All other edges of $M_A$ and $M_B$ appear in alternating paths or cycles; the cycles are completely in the intersection, and the paths have all vertices in $A \cap B$ except the endpoints, which {\it cannot} be in $A \cap B$.

Now we bound the number of ways to choose $M_A$ and $M_B$, given the positions of $A$ and $B$ (hence, $i$ and $j = \tilde{k} - i$ are known and fixed). To begin, we choose $M_A \triangle M_B$ - using the convention that $\triangle$ denotes symmetric difference - by the following process.

\begin{enumerate}
    \item Order the vertices of $ A\backslash B$ and of $B \backslash A$.

    \item Process the vertices $ v \in A \backslash B$ in order. If $v$ is not yet part of an edge in $M_A$ choose whether or not it is part of an $M_A$ edge and if so pick its neighbor $u$ (restricted to $A$, of course). If $ u \in A \cap B$ then choose whether or not this vertex is incident to an $M_B$ edge and if so choose its neighbor $w$ in $B$. If $ w \in A \cap B$ then it must also be in an edge in $M_A$; choose this neighbor in $A$ and continue this process until a vertex outside of $A \cap B$ is chosen.
    
    \item Perform the same procedure for each vertex in $B \backslash A$ (with the roles of $A$ and $B$ reversed).

    \item If there are still more $m_A \backslash m_B$ and $m_B \backslash m_A$ edges to choose, they must appear in alternating cycles. As long as there are more of these edges to choose: first, pick a vertex in $ v \in A \cap B$, then pick its neighbor in $u$ in an edge of $ M_A$ from $ A\cap B$ and then choose the neighbor of $u$ in an edges of $M_B$ and iteratively continue this process. At each even step (after the second step) choose whether or not to complete the cycle.
    
    \end{enumerate}

Throughout this procedure, we make series of two types of choices: we make some ``yes/no" choices (either starting a matching edge from $A \triangle B$, designating an edge from $M_A \triangle M_B$ to be isolated or not, or completing a cycle), and some choices regarding the next vertex in a path or cycle. Note that for each of the second type of choice there are at most $\tilde{k}$ possibilities. The procedure requires at most $2j + 2r_O$ choices of the first type. Indeed, such choices are made at particular vertices $A \triangle B$ to specify whether or not to include the vertex in a matching edge and after we place specific edges of $ M_A \triangle M_B$ to determine whether or not we continue an alternating path or cycle. The number of choices of the second type is {\it exactly} $2r_O$ as we make such a choice to complete each edge of $M_A \triangle M_B$ (except, in the case of cycles, the first choice of the ``second type" starts the cycle, and the last edge is chosen to complete the cycle by a choice of the first type, but the total number of second-type choices is still $2r_O$). Hence, the number of ways to set $M_A \triangle M_B$ is at most $4^j (2 \tilde{k})^{2 r_O}$. 

Now consider $M_A \cap M_B$. As these edge are contained in $A$, 
the number of ways to set these edges is at most 
$$ \frac{ [\tilde{k}]_{2(r - r_O)}}{(2^{r - r_O} (r - r_O)!)} \leq \frac{ [\tilde{k}]_{2r}}{(2^{r} r!)} \cdot  \left( \frac{3r}{\tilde{k}^2}\right)^{r_O}$$ 
(since $r_O \leq r \ll \tilde{k}$). 

Putting everything together we have (letting anything to the power of $j$ equal 1 if $j = 0$)
\[ | {\mathbb N}_{i, r_O} | \le
\binom{n}{\tilde{k}}\binom{n-\tilde{k}}{j}\binom{\tilde{k}}{j} \cdot 4^j (2 \tk)^{2 r_O} \cdot 
\frac{ [\tilde{k}]_{2r}}{(2^{r} r!)} \cdot  \left( \frac{3r}{\tilde{k}^2}\right)^{r_O}
\leq |\mathbb{N}|\left(\frac{4e^2 n \tilde{k}}{j^2}\right)^j (12r)^{r_O}. \]

\end{proof}

Now we bound $ \PP( \cX_A \wedge \cX_B)$. We achieve this bound by revealing $ G_{n,m}$ in a series of rounds, in some situations we reveal potential edges one at a time. Define
$\mathcal{E}$ to be the event that 
\begin{itemize}
    \item All edges of $M_A$ are edges, and all pairs of vertices in $A$ with at least one vertex not in $M_A$ are {\it not} edges
    \item All vertices not in $A \cup B$ have at least two neighbors in $A  \backslash ( \cup_{e \in M_A} e)$.
\end{itemize}
The event $\mathcal{E}$ is the same as the event $ \cX_A $ {\it except} that we ignore the minimum degree condition of vertices in $B \backslash A$. We first observe as much as $G_{n,m}$ as is necessary to verify $\mathcal{E}$. We do this so that we observe as few edges as possible; that is, for each vertex $v \notin A \cup B$ we reveal potential edges join $v$ and $ A \setminus ( \cup_{e \in M_A}e )$ one at at a time and we stop as soon as we observe that $v$ has two neighbors in this set
(note that we may come back to $v$ for verifying $ \cX_B$; we take care of this later in the proof). {\it Then} (if $\mathcal{E}$ happens) we look at the rest of the possible edge sites necessary to determine
whether or not the event $ \cX_A \wedge \cX_B$ occurs, again minimizing observations of edges (for verifying the condition for vertices outside $A$ or $B$ having edges to the inside, we stop once we observe two of these for each outer vertex into each set). Note that throughout the {\it entire} process, we observe at most $4n + 2r < 5n$ edges (the $2r$ coming from the matchings, the $4n$ for the minimum degree condition for vertices outside $A$ and/or $B$) and we establish the status of at most $2 n \tilde{k}$ vertex pairs (this is simply an upper bound on the edges that have at least one vertex in either $A$ or $B$). Therefore, whenever we reveal the status of a pair of vertices, the probability that an edge appears will be at least 
$$ p' = \frac{m - 5n}{\binom{n}{2}}$$ and at most 
$$ p'' = \frac{m}{\binom{n}{2} - 2n\tilde{k}}.$$ Note that $ p' \sim p$ and $ p'' \sim p$ and
\begin{equation} \label{equ: p' to p for big i}
    (1-p')^{\tk} \leq (1+o(1))(1-p)^{\tk}.
\end{equation}

We now bound the probability that each step of the process results in an outcome in $ \cX_A \wedge \cX_B$. 
Note that the only difference between $\mathcal{E}$ and $\cX_{A} $ is the minimum degree condition for vertices in 
$B \backslash A$. Recalling Lemma~\ref{lem:p to c}, we observe that the probability that 
the vertices in $ B \setminus A$ satisfy this minimum degree condition is at least 
\[ (1 - ( 1 - p')^{k-r} - (k-r)p'(1-p')^{k-r-1})^{j} \sim \varphi( p'(k-r))^j \ge \left(1 - \delta\right)^j,\]
for some $\delta>0$.
Therefore,
$$\PP[\mathcal{E}] \leq (1+2\delta)^j\PP(\cX_A).$$
Conditioning on the event $\mathcal{E}$, the event $ \cX_A \wedge \cX_B $ occurs if we have the following:
\begin{enumerate}
    \item All edges of $M_B \backslash M_A$ are actually edges. This occurs with probability at most $ (p'')^{r_O}$.
    \item All pairs of vertices in $\binom{B}{2} \backslash \left( \binom{A}{2} \cup \binom{ \cup_{e \in M_B} e }{2} \right)$ 
    must not be edges. Consider vertex pairs that are contained in $ \cup_{ e \in M_B} e$ but are not contained in $A$. One vertex
    of such a pair is in $B\setminus A$ and therefore in one of the $r_O$ outer edges of $ M_B$. Thus, the number of such pairs is at most $ 2 r_O \cdot 2r$. The probability that none of these pairs appear is at most
    $$(1 - p')^{(\tilde{k} - j)j + \binom{j}{2} - 4 r r_O} \leq  2^{r_O}(1 - p')^{(\tilde{k} - (j+1)/2)j}.$$
    \item Each vertex in $A \backslash B$ must have at least  two neighbors in $B \backslash ( \cup_{e \in M_B} e)$ and each vertex in $B \backslash A$ must have at least two neighbors in $A \backslash ( \cup_{e \in M_A} e)$ without violating our placements of $M_A$ and $M_B$ (specifically, the non-edge constraint outside of the matching vertex sets). We must now define two new variables: let $R_A$ be the set of edges $m_A$ which ``cross" into $A \cap B$ (i.e. these edges have one vertex in $A \setminus B$ and one vertex in $A \cap B$)  and do not intersect any edge in $M_B$, and let $r_A = |R_A|$. Define $ R_B, r_B$ analogously (see Figure 1).

    We bound the probability that the minimum degree conditions are satisfied as follows. We start with the vertices in $ B \setminus A$ that are in edges in $R_B$. Each such vertex $v$ has one neighbor in $A \setminus ( \cup_{e \in M_A}e)$ via the edge in $R_B$ containing $v$. The vertex $v$ must have an additional neighbor in $ A \setminus ( \cup_{e \in M_A} e)$. There are two possible choices for such a neighbor: a vertex in $A \backslash B$ or a vertex in $ A \cap B$ that is in an edge in $ R_B$; any other choice would violate the non-edge constraint inside $B$ or connect to a vertex that is an edge in $M_A$. As $ r_B \le j$ the number of possible neighbors is hence at must $2j$. So for this first part of the scheme our upper bound on the probability is $(2jp'')^{r_B}$ (or simply 1 if $j = 0$). Next, consider a vertex $ u \in A \backslash B$. The vertex $u$ must have two neighbors in $B \setminus ( \cup_{e \in M_B} e)$ (again, without violating non-edge constraints of $A$) {\it except} the vertices in $ \cup_{e \in R_A} e \backslash B$, which must have only one additional neighbor (the matching edge counts for the edge requirement in this case). Note that the edges from the previous step do not count, since they are incident with $M_B$. Hence an upper bound for the probability of this part of the scheme happening is $(2jp'')^{2j - r_A}$. (Note that there are some vertices in $ A \setminus ( B \cup ( \cup_{e \in R_A} e))$ that have $ r_A + j \le 2j$ potential neighbors. These are vertices that are in edges in $M_A$ that are not included in $R_A$. These vertices could potentially be adjacent to the vertices in $ ( \cup_{e \in R_A}e) \cap B$.)

    We have shown that $(2jp'')^{2j - r_A + r_B}$ is an upper bound for the probability of each vertex in $A \triangle B$ satisfying the minimum degree condition. However, by switching the roles of $A$ and $B$, we see that $(2jp'')^{2j - r_B + r_A}$ is also an upper bound on the probability of this event. Taking the smaller bound, we see that the probability of this event is at most $(2jp'')^{2j}$. Note that we could have $j> 1/(2p'')$. In this situation we simply bound this probability by 1. 

    Finally, each vertex outside of $A \cup B$ must have at least two neighbors in $B\backslash (\cup_{e \in M_B} e)$. We just upper bound this probability by 1.
\end{enumerate}
Putting everything together, we have
\begin{equation*}
    \PP[ \cX_A \wedge \cX_B] \leq (1 + 2\delta)^j \PP(\cX_A) \cdot (p'')^{r_O} \cdot  2^{r_O}(1-p')^{(\tilde{k} - (j+1)/2)j} \cdot \min\{(2jp'')^{2j},1\} 
\end{equation*}
Recalling (\ref{eq:h'''}), we have
\begin{align*}
    \sum_{i = C/p}^{\tilde{k}} h_i  & = \frac{1}{ \EE[ X_{k,r}]^2}  \sum_{i=C/p}^{\tk} \sum_{ \stackrel{ A, M_A; B ,M_B}{ |A \cap B|= i} } 
    \PP ( \cX_B \wedge \cX_A) -  \PP( \cX_A) \PP( \cX_B) \\
    &\leq \frac{1}{\EE[X_{k,r}]}\sum_{j = 0}^{\tilde{k} - C/p}\sum_{r_O = 0}^{r} \left( \frac{4e^2 n \tilde{k}}{j^2} \right)^{j} (12r)^{r_O} \cdot (2p'')^{r_O} (1 + 2\delta)^j \min\{(2jp'')^{2j},1\} (1 - p')^{(\tilde{k} - (j+1)/2)j}  \\&=
    \frac{1}{\EE[X_{k,r}]}\sum_{j = 0}^{\tilde{k} - C/p}\sum_{r_O = 0}^r\left(4e^2(1 + 2\delta) n \tilde{k} (1 - p')^{(\tilde{k} - (j+1)/2)} \cdot \min\{4p''^2,j^{-2}\} \right)^{j}\left(24rp''\right)^{r_O} \\&\leq
    \frac{1}{\EE[X_{k,r}]}\sum_{j = 0}^{\tilde{k} - C/p}\sum_{r_O = 0}^r\left(4e^2(1 + 3\delta) n \tilde{k} (1 - p)^{(\tilde{k} - j/2)} \cdot \min\{5p^2,j^{-2}\} \right)^{j}\left(25rp\right)^{r_O}  \qquad \text{by (\ref{equ: p' to p for big i})}\\&\leq
    \frac{1}{\EE[X_{k,r}]}\sum_{j = 0}^{\tilde{k} - C/p}\sum_{r_O = 0}^r\left(4e^2(1 + 4\delta) \tilde{k}^2 (1-p)^{(\tilde{k} - j)/2} \cdot \min\{5p^2,j^{-2}\}  \right)^{j}\left(25rp\right)^{r_O} \qquad \text{by (\ref{equ: k pre-approx})} \\&\leq
    \frac{2}{\EE[X_{k,r}]}\sum_{j = 0}^{\tilde{k} - C/p}\left(4e(1 + 4\delta) \tilde{k}^2 (1-p)^{(\tilde{k} - j)/2} \cdot \min\{5p^2,j^{-2}\}  \right)^{j} \qquad \qquad \text{since $rp = o(1)$}.
\end{align*}
It remains to observe that the expression raised to the power of $j$ is bounded by some constant less than 1. For $j < (1 - \delta)\tilde{k}$ , we have $(1-p)^{(\tilde{k} - j)/2} = n^{-\Omega(1)}$; then, bounding $\min\{5p^2, j^{-2}\}$ by $5p^2$ we see that the expression has the form $\Theta(\ln^2(n) n^{-\Omega(1)})$, which is $o(1)$.

For $j \geq (1 - \delta)\tilde{k}$, then the expression is of the form $\Theta( (1-p)^{(\tilde{k} - j)/2})$ (as $j^{-2} < 5p^2$), which is less than 1 so long as $ i = \tilde{k} - j = i$ is larger than a specific constant multiple of $1/p$. From the expression above such constant would be bigger than 1/2, but we now remark how the $4e$ term (which contributes a $(4e)^j$ multiplier) can be brought down to $(1 + c\delta)$ (for some universal constant $c$ and small enough  $\delta$) by a tightening of a couple of bounds earlier in the proof: first, there is a factor of $4^j$ which appears in the bound of $|\mathbb{N}_{i,r_O}|$. This just comes from the counting of ``yes/no" choices regarding whether each vertex from $A \backslash B$ and $B \backslash A$ is part of a matching edge in $A$ or $B$; however, this can be bounded by a term on the order of $\binom{j}{r}^2$, and since $r \ll j$ in this case, the $4^j$ can replaced by a $(1 + \delta)^j$. Secondly, one factor of $e^j$ comes from the $\binom{\tk}{j}$ factor in the bound of $\mathbb{N}_{i,r_O}$ (at the end of the proof of Claim \ref{claim: n_i_r}). For small $j$ this factor of $e^j$ is essential, but since $j \geq (1 - \delta)\tk$  here, $\binom{\tk}{j}$ can also be bounded by $(1 + c_1 \delta)^j$ for universal constant $c_1$ (by considering the slope of the entropy function). Finally, there is a $\tk^2$ and a $j^{-2}$ inside the parenthesis which when multiplied is bounded by $1 + 3\delta$. Hence, $\sum_{i = C/p}^{\tk} h_i$ is indeed $o(1)$ for any positive constant $C$.

\section{Poisson Paradigm} \label{sec:poisson}

This Section contains the proofs of Lemmas~\ref{lem:Phi}~and~\ref{lem:Psi}, which give our estimates for $ \Phi$ and $ \Psi$, respectively. Recall that $\Phi$ gives the probability that
all vertices outside a pair $ A, M_A$ that might counted by $ X_{k,r}$ have the required two neighbors in $A \setminus ( \cup_{e \in M_A} e)$, and
$ \Psi $ gives the analogous probability that pairs $ A, M_A$ and $ B,M_B$ both satisfy this minimum degree condition.

Both of these proofs make use of the Poisson paradigm. We achieve a tight estimate for $ \Phi$ by counting the number of bipartite graphs with vertex parts $ V \setminus A$ and $ A \setminus ( \cup_{e \in M_A}e )$ with a fixed number of edges and minimum degree 2 for all vertices $ V \setminus A$. We achieve this counting by sampling degrees from a truncated Poisson distribution, following methods introduced by Pittel and Wormald to count graphs with minimum degree 2 \cite{pittelwormald}.
We establish an upper bound on $ \Psi$ by passing to $ G_{n,p}$, as set forth in Lemma~\ref{lem: passing to Gnp}. The main difficulty here is dealing with the minimum degree condition for vertices in $ A \setminus B$ and $ B \setminus A$ as the edges between these two sets have an impact on the minimum degree condition for vertices in both sets. We handle this situation using a variation on Janson's inequality recently introduced by Bohman, Warnke and Zhu in the context of establishing 2-point concentration of the domination number of $ G_{n,p}$ \cite{bwz}.

\subsection{An asymptotic estimate for $ \Phi$: Proof of Lemma~\ref{lem:Phi}}

Note that Lemma~\ref{lem:p to c} implies that we can pass
between the expression for $\Phi$ in terms of $p$ and the expression
in terms of $ \varphi(c_p)$ where $ c_p = p(k-r) = c(m,k,r)$. Let $K$ be a fixed set of $k+r$ vertices and let $M \subset \binom{K}{2}$ be a matching consisting of $r$ edges.

We begin with the upper bound on $\Phi$. 
Here we pass from $ G_{n,m}$ to $ G_{n,p'}$ 
with $p' = p + (\log n) \max\{ p k^2/n^2, p^{1/2}/n\}$, as discussed in Section~\ref{sec:passthru}. Let $ \cE$ be the event that every vertex outside of $K$ as at least two neighbors in $ K \setminus( \cup_{e \in M} e )$. Note that $ \Phi = \PP_m( \cE \mid \cU)$. Lemma~\ref{lem:p to c} implies that we have $ \PP_{p'}( \cE) \sim \varphi( c_{p'})^{n-k-r} $. Next, note
that we have
\[ c_{p'} - c_p = O \left( \frac{ (\log n) p k^3}{ n^2 } \right) + O \left( \frac{ ( \log n) p^{1/2}k}{n}  \right) \ll \frac{ np^2}{ (\log n)^3}.\]
Thus, applying Lemmas~\ref{lem: passing to Gnp}~and~\ref{lem:cs} we have
\[  \Phi = \PP_m( \cE \mid \cU) \le(1 +o(1)) \PP_{p'}( \cE) \sim \varphi(c_{p'})^{n-k-r} \sim \varphi( c_p)^{n-k-r}. \]

It remains to establish the lower bound on $ \Phi$. Here we apply the Poisson paradigm. We condition on the number of 
edges between $ K \setminus ( \cup_{e \in M} e)$ and $ V \setminus K$ and consider a balls in bins model. For ease of notation set
\[ \beta = n-k-r \ \ \ \text{ and }  \ \ \   \gamma = k-r \ \ \ \text{ and } \ \ \ A = (n-k-r)(k-r),\]
and let $ \kappa$ be the number of edges between $ V \setminus K$ and $K \setminus ( \cup_{e \in M} e)$. We view the random bipartite graph 
between $ V \setminus K$  and $K \setminus ( \cup_{e \in M} e)$ as a $\beta \times \gamma$, $0-1$ matrix chosen uniformly at random
from all such matrices with exactly $\kappa$ 1's. 
Let $ f(\beta, \gamma, \kappa)$ be the probability that this random matrix has at least two 1's in every row. Then we have
\[ \Phi = \sum_{ \kappa = 0}^{ m-r} \frac{ \binom{A}{ \kappa} \binom{ \binom{n}{2} - M_1 -A}{ m-r- \kappa} }{  \binom{ \binom{n}{2} - M_1}{ m - r}}  f(\beta, \gamma, \kappa ), \]
where $ M_1 = \binom{k+r}{2} - \binom{2r}{2} + r$.

\begin{thm} \label{thm: enumeration} Let $ \delta >0$ be fixed. If $\kappa, \beta, \gamma \to \infty$ such that $\frac{1}{2} \beta \log \beta < \kappa < \frac{1}{\delta} \beta \log \beta$ and $\gamma < \beta < \gamma^{2 - \delta}$, then
    \begin{align*}
        f(\beta, \gamma, \kappa) &\geq (1+o(1)) \left(1 - (c+1) e^{-c} \pm e^{-2c + O(\log(c))}\right)^{\beta} 
    \end{align*}
where $ c = \kappa/ \beta$ is the expected number of 1's in each row of our matrix.
\end{thm}
\noindent
The proof of Theorem~\ref{thm: enumeration} is given below. 

We are now ready to complete the proof of the lower bound in Lemma~\ref{lem:Phi}.
Set 
$$ \kappa_{\rm low} = (m-r) \cdot \frac{A}{\binom{n}{2} - M_1} - (\log n)n^{1/2}  .$$  
Note that $ \kappa_{\rm low} = ( 1 +o(1)) nkp= (2 +o(1)) \log(np) n \ge (1/2 + \ep) n ( \log n)$.
Note
that $ f( \beta, \gamma, \kappa)$ is increasing in $ \kappa$ (by a simple 
coupling argument). Thus we have
\[ \Phi \ge \sum_{ \kappa = \kappa_{\rm low} }^{m-r}\frac{ \binom{A}{ \kappa} \binom{ \binom{n}{2} - M_1 -A}{ m-r- \kappa} }{  \binom{ \binom{n}{2} - M_1}{ m - r}} f(\beta, \gamma, \kappa_{\rm low}) 
= f(\beta, \gamma, \kappa_{\rm low}) \left( 1 - \sum_{ \kappa =0}^{ \kappa_{\rm low}-1 } \frac{ \binom{A}{ \kappa} \binom{ \binom{n}{2} - M_1 -A}{ m-r- \kappa} }{  \binom{ \binom{n}{2} - M_1}{ m - r}} \right). \]
Now note that the final sum is simply the probability of a particular deviation 
for a hypergeometric random variable. Note that
\[ \EE[ \kappa] = \frac{ (m-r)( n-k-r)(k-r)}{ \binom{n}{2} - M_1} = (1 +o(1)) pnk. \]
Thus, by Theorem~\ref{thm:hyperchernoff}, the final sum is 
bounded above by
\[ \exp \left\{  - \Omega \left( \frac{ (\log n)^2n}{ nkp } \right) \right\} = \exp \left\{ - \Omega( \log n) \right\}. \]
We have established
$ \Phi \ge (1 - o(1)) f( \beta, \gamma, \kappa_{\rm low})$.
Set $ c_{\rm low} = \kappa_{\rm low}/\beta$.
Now we note that $ p \ge n^{-3/4 + \epsilon}$ implies 
\[ e^{-2c_{\rm low}} = (1+o(1)) \left(\frac{ \log (np)}{ np} \cdot \frac{2}{e} \right)^4 = O(n^{-1-3\ep}),\]
and therefore
\[ \left(1 - (c_{\rm low}+1) e^{-c_{\rm low}} + e^{-2c_{\rm low} + O(\log(c_{\rm low}))}\right)^{\beta} \sim \varphi( c_{\rm low}). \]
Finally, we note that $ |c_{\rm low} - c_p| = O \left( ( \log n) n^{-1/2}\right)$,
and $ p \ge n^{-3/4 + \epsilon}$ implies that this quantity is $ o(( \log n)^{-3} n p^2)$. Thus, by Lemma~\ref{lem:cs} we have $ \varphi( c_{\rm low})^{\beta} \sim \varphi( c_p)^{\beta}$. Putting it all together, using Theorem \ref{thm: enumeration},
\begin{multline*}
    \Phi \ge (1 - o(1)) f( \beta, \gamma, \kappa_{\rm low})
\ge ( 1 -o(1)) \left(1 - (c_{\rm low}+1) e^{-c_{\rm low}} + e^{-2c_{\rm low} + O(\log(c_{\rm low}))}\right)^{\beta} \\
\sim \varphi( c_{\rm low})^\beta \sim \varphi( c_p)^{\beta}, 
\end{multline*}
as desired.

\subsubsection{Proof of Theorem \ref{thm: enumeration}}

For positive integer variables $\beta, \gamma, \kappa$ with $\kappa \leq  \beta \gamma$, we choose a $\beta \times \gamma$ - sized matrix with entries in $\{0,1\}$ uniformly at random from the collection of all such matrices that have exactly $\kappa$ 1's. We are interested in the probability $f(\beta, \gamma, \kappa)$ of every row having at least two 1's.
We work directly with the total number of valid matricies $C(\beta, \gamma, \kappa) := \binom{\beta \gamma}{\kappa} f(\beta, \gamma, \kappa)$. We establish a lower bound on $ C(\beta, \gamma, \kappa)$ by working with random
matrices in which the number of 1's in each row is given by an independent truncated Poisson random variable with minimum degree 2.

\subsubsection{Preliminaries}

In this Subsection we make some preliminary definitions and state some
facts from \cite{pittelwormald}.

To begin, define, for any degree-sequence vector $ \vd = (d_1,\dots,d_{\beta} )$,
$$\eta( \vd ) := \frac{1}{\kappa}\sum_{i=1}^{\beta} d_i(d_i-1).$$
In order to motivate this definition, let $ g(\vd)$ be the number of ways to choose a 0-1 matrix of size $\beta \times \gamma$ such that row $i$ has $d_i$ 1's. Setting
$ \kappa= d_1 + \dots + d_{\beta} $ and letting $ \| \vd \|_{\infty}$ be the maximum of the $d_i$'s we have
\begin{align}
g( \vd) = \prod_{i=1}^{\beta}\binom{\gamma}{d_i}&= \left(\prod_{i=1}^{\beta}[\gamma]_{d_i}\right)\left(\prod_{i=1}^{\beta}\frac{1}{d_i!}\right) \nonumber \\&=
\gamma^{\kappa}\left(\prod_{i=1}^{\beta} \prod_{j=0}^{d_i-1} \left(1 - \frac{j}{\gamma}\right)\right)\left(\prod_{i=1}^{\beta}\frac{1}{d_i!}\right) \nonumber \\&=
\left(1 + O(\kappa \| \vd \|_\infty^2\gamma^{-2})\right) \gamma^{\kappa} e^{-\kappa \eta(\vd)/2\gamma} \left(\prod_{i=1}^{\beta}\frac{1}{d_i!}\right). \label{equ: g formula}
\end{align}

Next we introduce and recall some properties of truncated Poisson random variables.
Let $\chi = \chi(\lambda)$ be the truncated Poisson random variable with parameter $\lambda>0$ and a lower bound 2: for all $j \in \NN $,
\begin{equation}
    \Pr[\chi = j] = \begin{cases}
        \frac{\lambda^j}{j! (e^\lambda - \lambda - 1)}, &j \geq 2 \\
        0, &j < 2.
    \end{cases} \label{equ: chi def}
\end{equation}
Next, recall that $ c = \kappa/\beta$, and for $c >2$,
define $\lambda_c$ to be the positive value of $\lambda$ for which
$$ \EE[\chi(\lambda)]= \frac{\lambda (e^\lambda - 1)}{e^\lambda - \lambda - 1} = c$$ 
(Lemma \ref{lem: lambda and eta asymptotics} below gives a proof of the existence and uniqueness of $\lambda_c$). Finally, define
$$\overline{\eta}_c := \frac{\lambda_c e^{\lambda_c}}{e^{\lambda_c} - 1},$$
and note that
\begin{equation}
    \EE[\chi(\lambda_c)(\chi(\lambda_c)-1)] = c \overline{\eta}_c. \label{equ: eta c motivation}
\end{equation}
Hence, if the components of a $\beta$-length vector $\vd$ consist of independent copies of $\chi(\lambda_c)$, then $\EE[\eta( \vd)] = \overline{\eta}_c$; this motivates our definition of $\overline{\eta}_c$. 

Now we show $\lambda_c$ is well-defined, as well as give asymptotic formulas for $\lambda_c$ and $\overline{\eta}_c$.

\begin{lem} \label{lem: lambda and eta asymptotics}
If $c > 2$ then $\lambda_c$ is well-defined and 
\begin{align}
    \lambda_c &= c - c^2 e^{-c} \pm  e^{-2c + O(\log(c))} \label{equ: lambda asymptotics} \\
    \overline{\eta}_c &= c \pm e^{-c + O(\log(c))} \label{equ: eta asymptotics}
\end{align}

\end{lem}

\begin{proof}

Note that $\EE[\chi(\lambda)]$ is strictly increasing in $\lambda$ (see the first paragraph of the proof of Lemma 1 in \cite{pittelwormald}, and compare (5) in \cite{pittelwormald} with (\ref{equ: chi def}) above for notation differences) and that $ \lim_{ \lambda \to 0} \EE[\chi( \lambda)] = 2$; hence, there exists a unique $\lambda_c$ for which
\begin{equation}
    \EE[\chi(\lambda_c)] = \frac{\lambda_c(e^{\lambda_c} -1)}{e^{\lambda_c} - \lambda_c - 1} = \lambda_c \left(1 + \frac{\lambda_c}{e^{\lambda_c} - \lambda_c - 1}\right) = c. \label{equ: lambda c and c}
\end{equation}
Note that $\EE[\chi(\lambda)]$ is unbounded for $\lambda$ large and that we
can restrict our attention to
$c$ large due to $O(\log(c))$ term in the final exponent of (\ref{equ: lambda asymptotics}); hence, by (\ref{equ: lambda c and c}), $\lambda_c$ must be large as well.
With $\lambda$ large, one can show that 
\begin{equation}
    \frac{d \EE[\chi(\lambda)]}{d\lambda} = 1 + o_{\lambda}(1). \label{equ: expectation of chi derivative}
\end{equation}
Next, note that we have
\begin{align}
    \EE[\chi(c - c^2 e^{-c})] &= (c - c^2 e^{-c}) \left(1 + \frac{c - c^2 e^{-c}}{e^{c - c^2 e^{-c}} - (c - c^2 e^{-c}) - 1}\right) \nonumber \\&=
    (c - c^2 e^{-c})\left(1 + (c - c^2 e^{-c})e^{-c}(1 + O(c^2e^{-c})\right) \nonumber \\&=
    c \pm e^{-2c + O(\log(c))}. \label{equ: almost to lambda c}
\end{align}
Finally, combining (\ref{equ: expectation of chi derivative}) and (\ref{equ: almost to lambda c}) gives (\ref{equ: lambda asymptotics}).

We now use (\ref{equ: lambda asymptotics}) and the definition of $\overline{\eta}_c$ to verify (\ref{equ: eta asymptotics}):
\begin{align*}
    \overline{\eta}_c = \frac{\lambda_c e^{\lambda_c}}{e^{\lambda_c} - 1} &= \lambda_c \left(1 + \frac{1}{e^{\lambda_c} - 1}\right) \\&=
    (c - e^{-c + O(\log(c))} )\left(1 + \frac{1}{\exp\{c - e^{-c + O(\log(c))}\} - 1}\right) \\&=
    c \pm e^{-c + O(\log(c))}.
\end{align*}
\end{proof}

In our application of Lemma~\ref{lem: lambda and eta asymptotics}, we only need (\ref{equ: eta asymptotics}) 
to verify that $\overline{\eta}_c$ converges to $c$ sufficiently quickly as $c \to \infty$. Note that 
(\ref{equ: lambda asymptotics}) can be used to verify that
\begin{equation}
    (e^{\lambda_c} - \lambda_c - 1) (c e^{-1} \lambda_c^{-1})^c = 1 - (c+1)e^{-c} \pm e^{-2c + O(\log(c))}. \label{equ: Poisson asymptotics}
\end{equation}
Finally, we state two large deviation estimates for truncated Poisson random variables
established by Pittel and Wormald:
\begin{lem}[Pittel and Wormald \cite{pittelwormald}] \label{lem:pwA} 
Let $ \vchi $ be a vector of length $ \beta$ consisting of independent truncated Poisson random variables with parameter $ \lambda_c$ where $ c = \kappa/\beta$ and minimum value $2$. If $ \kappa = O( \beta \log( \beta))$ then we have
\[ \PP \left( | \eta( \vchi ) - \EE( \eta( \vchi ))| \ge \beta^{1/2} \kappa^{-1} ( \log \beta)^8 \right) 
\le \exp\left\{ - \Omega(( \log \beta)^3) \right\}.\]
\end{lem}
\begin{proof}
This is equation (33) in \cite{pittelwormald}, with $2m$ replaced by $\kappa$, $n$ replaced by $\beta$, $Y$ replaced by $\chi$, $k$ set equal to 2, and $S$ set equal to $ \eta( \vchi)/2$  (see (28) in \cite{pittelwormald} for the definition of $S$). Note that this requires the condition $\kappa = \Theta( \beta (\log \beta)).$
\end{proof}
\begin{lem}[Pittel and Wormald \cite{pittelwormald}]
If  $\chi $ is a truncated Poisson random variables with parameter 
$\lambda$ and minimum value $2$ then 
\begin{equation}
    \PP[\chi \geq j_0] = O(e^{-j_0/2})  \ \ \text{ for \ $j_0 > 2e\lambda.$} \label{equ: upper tail}
\end{equation}
\end{lem}
\begin{proof}
This is (27) in \cite{pittelwormald}. It follows from the simple observations that we have
\[ \sum_{ j \ge 2}  \frac{ \lambda^j}{ j!( e^\lambda - \lambda -1)} =1 \]
and the ratios of consecutive terms in this sum are at most $e^{-1}$ for $j \ge j_0/2$.
\end{proof}

\subsubsection{Enumerating by passing to independent random degrees}

Given a degree sequence vector $ \vd= (d_1, \dots, d_\beta)$, define $$U({\vd}) := g({\vd}) \left(\prod_{i=1}^{\beta}d_i!\right) = \prod_{i=1}^{\beta}[\gamma]_{d_i},$$
and note that (\ref{equ: g formula}) implies
\begin{equation}
    U( \vd) = \left(1 + O(\kappa \| \vd \|_\infty^2\gamma^{-2})\right) \gamma^{\kappa} e^{-\kappa \eta(\vd)/2\gamma}. \label{equ: U formula}
\end{equation}

We appeal to generating functions, following Pittel and Wormald. Note that we have
\begin{align*}
    C(\beta, \gamma, \kappa) &= \sum_{ \substack{ d_1, \dots, d_\beta \ge 2 \\ d_1 + \dots + d_{\beta} = \kappa}} \prod_{i=1}^{\beta}\binom{\gamma}{d_i} \nonumber =
    [x^{\kappa}]\left( \sum_{d \geq 2}  \binom{\gamma}{d} x^{d}\right)^{\beta 
    }.
\end{align*}

Let $\chi_1, \dots, \chi_{\beta}$ be independent copies of $\chi(\lambda)$, and let $\vec{\chi} := ( \chi_{1},\dots,\chi_{\beta})$. For any $ \lambda >0$ we have
\begin{align}
    C(\beta, \gamma, \kappa) &= \frac{(e^\lambda - \lambda - 1)^{\beta}}{\lambda^{\kappa}} \cdot [x^{\kappa}] \left( \frac{1}{(e^\lambda - \lambda - 1)}\sum_{d \geq 2}  \binom{\gamma}{d} (\lambda x)^{d}\right)^{\beta 
    }\nonumber  \\&=
    \frac{(e^\lambda - \lambda - 1)^{\beta}}{\lambda^{\kappa}} \EE\left[U(\vec{\chi}) {\bf 1}_{\left\{\sum_{j=1}^{\beta} \chi_j = \kappa\right\}}\right] \nonumber  \\&=
    \frac{(e^\lambda - \lambda - 1)^{\beta}}{\lambda^{\kappa}} \EE\left[U(\vec{\chi})  \mathrel{\Big|} \sum_{j=1}^{\beta} \chi_j = \kappa\right] \PP\left[\sum_{j=1}^{\beta} \chi_j = \kappa\right]. \label{equ: CR generate complete}
\end{align}
Thus computing the last two terms in (\ref{equ: CR generate complete}) is the main task here.
\begin{thm} \label{thm:U details}
    Let $ \delta >0$. Suppose that $\kappa, \beta, \gamma \to \infty$ such that $\frac{1}{2} \beta \log \beta < \kappa < \frac{1}{\delta} \beta \log \beta$ and $\gamma < \beta < \gamma^{2-\delta}$. Set $c = \kappa/\beta$ and choose $ \lambda_c$ so that $ \EE[ \chi( \lambda_c) ] = c $. If $ \chi_1, \dots, \chi_\beta$ are i.i.d. copies of $ \chi( \lambda_c)$ then we have
    \begin{enumerate}[label=(\alph*)]
        \item 
        \begin{equation*}
        \PP\left[\sum_{j=1}^{\beta} \chi_j = \kappa\right] = \frac{1 + O(\kappa^{-1})}{\sqrt{2 \pi \kappa (1 + \overline{\eta}_c - c)}}
        \end{equation*}

        \item     
        \begin{equation*}
            \EE\left[U(\vec{\chi})  \mathrel{\Big|}  \sum_{j=1}^{\beta} \chi_j = \kappa  \right] \geq (1 + O(\beta^{1/2} (\log^8 \beta) \gamma^{-1}))\gamma^{\kappa}e^{-\kappa \overline{\eta}_c/2\gamma}
        \end{equation*}

    \end{enumerate}
    \end{thm}
\noindent The proof Theorem~\ref{thm:U details} is given in the following subsection.

Setting $\lambda = \lambda_c$ in (\ref{equ: CR generate complete}) and applying Theorem~\ref{thm:U details} we have
    \begin{align*}
        C(\beta, \gamma, \kappa) &\geq (1 + O(\beta^{1/2} (\log \beta)^8 \gamma^{-1})) \frac{\gamma^{\kappa }(e^{\lambda_c} - \lambda_c - 1)^{\beta}}{(\lambda_c)^{\kappa}e^{\kappa \overline{\eta}_c/2\gamma}\sqrt{2\pi \beta c (1 + \overline{\eta}_c - c)}}. 
    \end{align*}
Now we recall a standard bound on binomial coefficients (E.g. see Lemma 21.1 in \cite{alanmichal}). As $ \kappa = O( \beta ( \log \beta)) = o\left( (\beta \gamma)^{3/4} \right)$, we have 
\begin{multline*}
\binom{ \beta \gamma}{ \kappa } = \left(1 + O \left( \kappa^4 ( \beta \gamma)^{-3} \right) + O \left( \kappa ( \beta \gamma)^{-1} \right )\right) \frac{ (\beta \gamma)^\kappa}{ \kappa!} \exp \left\{ -\frac{ \kappa^2}{ 2 \beta \gamma} - \frac{ \kappa^3}{ 6 \beta^2 \gamma^2} \right\} \\ = (1 - O(1/\kappa) - O( \kappa^3 \beta^{-2} \gamma^{-2})) (2 \pi \kappa)^{-1/2} \left( \frac{ e \beta \gamma}{ \kappa}   \right)^\kappa\exp \left\{ -\frac{ \kappa^2}{ 2 \beta \gamma} \right\}
\end{multline*}
So, with the aid of (\ref{equ: eta asymptotics}), (\ref{equ: Poisson asymptotics}), and 
recalling that $c = \kappa/\beta$ and our bounds on $\kappa$ (which imply $e^{-c} \le \beta^{-1/2}$) and $\gamma$ with respect to $\beta$, we have
\begin{align*}
    f(\beta,\gamma,\kappa) &= \binom{\beta \gamma}{\kappa}^{-1} C(\beta, \gamma, \kappa) \\&\geq
    (1 + O(\beta^{1/2}(\log \beta)^8 \gamma^{-1} + \kappa^{-1} + \kappa^3 (\beta \gamma)^{-2})) \\
    & \qquad \times \left(\frac{e\beta\gamma }{\kappa e^{\kappa/(2\beta\gamma)}}\right)^{-\kappa}(2 \pi \kappa)^{1/2} \left(\frac{\gamma^{\kappa }(e^{\lambda_c} - \lambda_c - 1)^{\beta}}{(\lambda_c)^{\kappa}e^{\kappa \overline{\eta}_c/2\gamma}\sqrt{2\pi \beta c (1 + \overline{\eta}_c - c)}}\right) \\&=
    (1+O(\beta^{1/2}(\log \beta)^8 \gamma^{-1}))\frac{c^{\kappa} \left( e^{\lambda_c} - 1 - \lambda_c \right)^{\beta}}{e^{\kappa}\lambda_c^\kappa} \left(e^{\kappa(c - \overline{\eta}_c)/(2\gamma)}(1 + \overline{\eta}_c - c)^{-1/2}\right) \\&=
    (1+O(\beta^{- \delta/4})) \left((e^{\lambda_c} - \lambda_c - 1) (c e^{-1} \lambda_c^{-1})^c\right)^{\beta} \\&=
    (1+o(1)) \left(1 - (c+1) e^{-c}  \pm e^{-2c + O(\log(c))}\right)^{\beta},
\end{align*}
as desired.

\subsubsection{Proof of Theorem~\ref{thm:U details}}
\begin{proof}
Statement~(a) is identical to statement~(a) of Theorem~4 of \cite{pittelwormald} (with $2m$ replaced by $\kappa$, $n$ replaced by $\beta$, $Y$ replaced by $\chi$, and $k$ replaced by $2$). Note that Theorem~4 of \cite{pittelwormald} requires $\kappa = O(\beta \log \beta)$, which is a condition that we impose here.

We now turn to the proof of Statement~(b). 
    By Statement~(a) and Lemma~\ref{lem:pwA} we have
    \begin{align*}
        \PP\left[(\kappa/2\gamma)|\eta(\vec{\chi}) - \EE[\eta(\vec{\chi})]| \geq \gamma^{-1} \beta^{1/2 } (\log \beta)^8 \mathrel{\Big|} \sum_{j=1}^{\beta} \chi_j = \kappa\right] \leq 
        \exp\{-\Omega((\log \beta)^3)\}
    \end{align*}
Next we bound the probability that any copy of $\chi$ is too large: by (\ref{equ: upper tail}), and since $ (\log \beta)^3 \gg c$ for large $\beta$, we have 

$$\PP\left[||\vec{\chi}||_{\infty} > ( \log \beta)^3 \mathrel{\Big|} \sum_{j=1}^{\beta} \chi_j = \kappa\right] = \exp\{-\Omega( ( \log \beta)^3 )\}.$$
Therefore
    \begin{align*}
        \EE\Big[U(\vec{\chi}) & \mathrel{\Big|}  \sum_{j=1}^{\beta} \chi_j = \kappa \Big] \\ & \geq \EE\left[U(\vec{\chi}) {\bf 1}_{||\vec{\chi}||_{\infty} \leq \log^3\beta}   \mathrel{\Big|} \sum_{j=1}^{\beta} \chi_j = \kappa\right] \\ & =
        (1 + O(\beta (\log \beta)^7 \gamma^{-2})) \EE\left[\gamma^{\kappa}e^{-\kappa \eta(\vec{\chi})/2\gamma} {\bf 1}_{||\vec{\chi}||_{\infty}\leq \log^3\beta}  \mathrel{\Big|} \sum_{j=1}^{\beta} \chi_j = \kappa\right] \\ & =
        (1+O(\beta (\log \beta)^7 \gamma^{-2})) \gamma^{\kappa} e^{-\kappa \overline{\eta}_c/2 \gamma} \EE\left[ e^{-(\kappa/2\gamma) (\eta(\vec{\chi}) - \EE[\eta(\vec{\chi})])} {\bf 1}_{||\vec{\chi}||_{\infty}\leq \log^3\beta} \mathrel{\Big|} \sum_{j=1}^{\beta} \chi_j = \kappa\right] \\ & \geq 
        (1+O(\beta (\log \beta)^7 \gamma^{-2})) \gamma^{\kappa} e^{-\kappa \overline{\eta}_c/2 \gamma} e^{ - \gamma^{-1} \beta^{1/2} ( \log \beta)^8 } 
        \left( 1 - e^{ -\Omega( (\log \beta)^3)} \right) \\ & =
        \Big(1 + O(\beta (\log \beta)^7 \gamma^{-2}) + \beta^{1/2} (\log \beta)^8 \gamma^{-1} + \exp\{-\Omega((\log \beta)^3 )\})\Big) \gamma^{\kappa} e^{-\kappa \overline{\eta}_c/2 \gamma} \\ & =
        (1 + O(\beta^{1/2} (\log \beta)^8 \gamma^{-1}) ) \gamma^{\kappa} e^{-\kappa \overline{\eta}_c/2 \gamma}.
    \end{align*}

\end{proof}

\subsection{An upper bound on $ \Psi$: Proof of Lemma~\ref{lem:Psi}.} \label{sec:Psi}

Let $A, M_A$ and $B, M_B$ be fixed where $A$ and $B$  are sets of 
$ \tk = k+r$ vertices and $M_A \subset \binom{A}{2}, M_B \subset \binom{B}{2}$
are matchings such that $|M_A| = |M_B| = r$. Let $ I = A \cap B $
and set
\[ i := |I| = |A \cap B|.\]
We will also need to keep track of 
\[ a = \left| B \cap \left( \cup_{e \in M_A} e \right) \right| \ \ \ \text{ and } 
\ \ \ b = \left| A \cap \left( \cup_{e \in M_B} e \right) \right|. 
\]
Let $ \cF$ be the event that every vertex outside of $A$ has at
least two neighbors in $A \setminus ( \cup_{e \in M_A} e)$ and every vertex outside of $B$ has at least
two neighbors in $B \setminus ( \cup_{e \in M_B}e)$. Our goal in this Section is to bound
$ \Psi = \PP_m( \cF \mid \cU_A \wedge \cU_B)$.
We work in $ G_{n,p'}$, with 
$p' = p + (\log n) \max\{  p k^2/ n^2, p^{1/2}/n \}$ 
as set forth in Lemma~\ref{lem: passing to Gnp}. We have $ \Psi \le (1 +o(1))\PP_{p'}(\cF \mid \cU_A \wedge \cU_B)$. 
We begin by revealing the edges between
$I$ and the complement of $A \cup B$. Let $ \Gamma =\Gamma(I)$ be the number of edges joining $I$ and $ V \setminus (A \cup B)$. Note
that we have
\[
\mu_\Gamma := \EE_{p'}[\Gamma \mid \cU_A \wedge \cU_B] = p'i(n-2 \tk +i). \]
Let $ N_1$ be the number of vertices in $ V \setminus A \cup B$ that have 
exactly one neighbor in $ I$, and let $ N_2$ be the number of vertices in $ V \setminus (A \cup B)$ that have at least two neighbors in $I$.
Note that we have the simple bound $ N_1 \le \Gamma$. Set
\[ J = N_1 -p'i(n- 2\tk +i). \]
By the Chernoff bound we have
\begin{align}
\label{eq:J bound low}
j < npi \ \ \ \Rightarrow \ \ \ 
\PP_{p'}( J \ge j )
\le \exp\left\{ - \frac{ j^2}{ 3npi}  \right\} \\
\label{eq:J bound high}
j \ge npi \ \ \ \Rightarrow \ \ \ 
\PP_{p'}( J \ge j )
\le \exp\left\{ - \frac{ j}{ 3}  \right\}.
\end{align}
Furthermore, setting $ \xi = \max\{  4 n i^2 p^2, ( \log n) p k^3/n \} $, we have 
\begin{equation}
\label{eq:N2 bound}
\PP_{p'}(N_2 \ge \xi ) \le \binom{n}{\xi} \binom{i}{2}^{\xi} (p')^{2\xi}
\le \left( \frac{ne i^2 p'^2}{ 2 \xi}   \right)^{\xi} \le \exp\{ - \Omega( (\log n) p k^3/n ) \} = o( \Phi^2). 
    \end{equation}
Let $ \cB$ be the event $ \{ N_2 \ge \xi \}$.

Now consider the
edges between $ A \Delta B$ and $ V \setminus ( A \cup B )$.
Let $ \cF_A$ be the event that all vertices outside of $ A \cup B \cup N_1 \cup N_2$ have
at least two neighbors in $A \setminus ( I \cup( \cup_{e \in M_A} e))$ and all vertices 
in $N_1$ have at least one vertex in this set. We define $ \cF_B$ analogously. For ease of
notation, set $ i_a = i + 2r -a$. We have 
\begin{equation*}
    \begin{split}
\PP_{p'}( \cF_A & \mid N_1, N_2 ) \\
& \le \left( 1 - (1-p')^{\tk-i_a} - (\tk-i_a)p'(1 - p')^{\tk-i_a-1} \right)^{ n -2\tk+i - N_1 - N_2} \left( 1 - (1-p')^{\tk-i_a}\right)^{ N_1} \\
& \le \exp \left\{ - (1-p')^{\tk-i_a} \left[ ( n -2\tk+i - N_1 - N_2 )( 1 + (\tk-i_a)p'(1 -p')^{-1}) + N_1 \right]   \right\} \\
& \le  \exp \left\{ - (1-p')^{\tk-i_a} \left[ (\tk-i_a)p'( n -2\tk+i - N_1)  + ( n -2\tk+i) \right] + O( N_2pk^3/n^2)   \right\} \\
& = \exp \left\{ - (1-p')^{\tk-2r-i} \left[ (\tk - 2r-i)p'( n -2\tk+i - N_1)  + ( n -2\tk+i) \right] \right. \\
& \hskip2cm  \left. +O (a\cdot k^3p^2/n) + O( N_2pk^3/n^2)   \right\}.
\end{split}
\end{equation*}
The analogous bound holds for $ \cF_B$. 

Next let $ \cF_C$ be the event that every vertex in $A \setminus I$ has at least two neighbors in 
$ {B \setminus ( \cup_{e \in M_B}e)}$ and every vertex in $B \setminus I$ has at least two neighbors in 
$ A \setminus ( \cup_{e \in M_A}e)$. 
It follows from Lemma~\ref{lem:withLutz} (stated and proved below - see (\ref{eq:bipartite conclusion})), letting $A \backslash B$ and $B \backslash A$ be the partite sets and $T = {((\cup_{e \in M_A}e) \backslash B) \cup ((\cup_{e \in M_B}e) \backslash A)}$, that we have 
\begin{equation*}
\begin{split}
\PP_{p'}( \cF_C) & \le  \exp\left\{ - \left( 2(\tk-2r-i) +a+b\right)\left[ (\tk-i) p' (1-p')^{\tk-i-1} + (1-p')^{\tk-i} \right] + O( (\log n) p^3 k^8/ n^4) \right\} \\   
& \le \exp\left\{ - 2 (1-p')^{\tk-i} (\tk - 2r -i)\left[ (\tk-i) p' + 1 \right] + O( (\log n) p^3 k^8/ n^4) \right\}. 
\end{split}
\end{equation*}

Putting these bounds together, observing that $ \cF_A, \cF_B$ and $ \cF_C$ are independent in $ G_{n,p'}$, we have 
\begin{equation*}
    \begin{split}
\PP_{p'} &( \cF_A \wedge \cF_B \wedge \cF_C \mid N_1, N_2 ) \\
& \le \exp \left\{ - 2(1-p')^{\tk-2r-i} \left[ (\tk- 2r-i)p'( n -\tk - N_1 ) + ( n -\tk ) \right] \right\} \\
& \hskip2cm \times \exp \left\{ O ((a+b)\cdot k^3p^2/n)+ O( N_2 \cdot pk^3/n^2)  +O( p^3k^7/n^3)  \right\}  \\
& \le   \exp \left\{ -  2 (n-k-r)(1-p')^{k-r} \cdot \Upsilon 
 +O( p^3k^7/n^3)  +O ((a+b)\cdot k^3p^2/n) + O( N_2 \cdot pk^3/n^2)  
\right\} 
\end{split}
\end{equation*}
where 
\begin{equation*}
\begin{split}
\Upsilon &= p'(k-r) \left( 1 - \frac{i}{k-r} \right) (1 - p')^{-i} \left[ 1 - \frac{N_1}{n-k-r}\right] + (1-p')^{-i} \\
& \ge p'(k-r) \left( 1 - \frac{i}{k-r} \right) (1 + p'i) \left[ 1 - \frac{ p'i(n-2\tk+i) + J}{n-k-r}\right] + (1+p'i)  \\
& \ge p'(k-r) \left( 1 - \frac{i}{k-r} \right) (1 + p'i) \left[ 1 - p'i - \frac{J}{n-k-r}\right] + 1+p'i   \\
& \ge  \left[ p'(k-r) +1 \right] + p'(k-r) \left[- \frac{i}{k-r} - p'^2i^2- \frac{ J (1 +p'i)}{n-k-r}\right]  +p'i \\
& = \left[ p'(k-r) +1 \right] - p'(k-r) \left[p'^2i^2 + \frac{ J (1 +p'i)}{n-k-r}\right]. 
\end{split}
\end{equation*}
Inserting this bound on $ \Upsilon$ into the expression above, recalling the bound on $ (1-p')^k/(1-p)^k$ given in
Lemma~\ref{lem: passing to Gnp}, and using Lemma \ref{lem:Phi}, we have 
\begin{equation*}
\begin{split}
  \PP_{p'}  & ( \cF_A \wedge \cF_B \wedge \cF_C \mid N_1, N_2 ) \\
& \le \exp \left\{ -2(n-k-r) ( 1 - p)^{k-r} (p(k-r) +1) + O \left( p^{5/2}k^5/n^2  \right)    \right\} \\
& \hskip1cm  \times
\exp \left\{  O \left( \frac{ p k^3}{n} \cdot \left[ p^2i^2 + \frac{J}{ n} \right] \right) +O( p^3k^7/n^3)  +O ((a+b)\cdot k^3p^2/n) + O( N_2 \cdot pk^3/n^2)  \right\} \\
& \le \Phi^2 \exp \left\{  O \left( \frac{ p k^3}{n} \cdot \left[ p^2i^2 + \frac{J}{ n} \right] \right) +O ((a+b)\cdot k^3p^2/n) + O( N_2 \cdot pk^3/n^2)  + o(1) \right\}. 
\end{split}
\end{equation*}

Set $ j_0 = n^2 / ( (\log n) pk^3) $ and $ j_1= nip$.
Recalling Lemma~\ref{lem: passing to Gnp} we have $ \Psi$ is bounded above by $ (1 +o(1))$ times
\begin{multline*}
\PP_{p'}( \cF_A \wedge \cF_B \wedge \cF_C ) 
 \le \PP_{p'}( \cB) 
+ \PP_{p'} ( \cF_A \wedge \cF_B \wedge \cF_C \mid \overline{ \cB} \wedge \{ J \le j_0\})\\
 \hskip2.5cm + \sum_{j \ge j_0}  \PP_{p'} ( \cF_A \wedge \cF_B \wedge \cF_C \mid \overline{ \cB} \wedge \{ J = j\}) \PP_{p'} (J \ge j).
\end{multline*}
Recall that (\ref{eq:N2 bound}) gives a bound on $\PP_{p'}( \cB)$.  Note that if $ N_2 < \xi$ then the corresponding term is absorbed into other
error terms. For $ J \le j_0$ we observe that the corresponding term is $ o(1) $. With these observations in hand, and recalling (\ref{eq:J bound low})~and~(\ref{eq:J bound high}), we have 
\begin{equation*}
\begin{split}
\Psi 
& \le o( \Phi^2) + \Phi^2 \exp\left\{ O\left(i^2 \cdot p^3k^3/n   \right)  +O ((a+b)\cdot k^3p^2/n)    + o(1)  \right\} \\
& \hskip3cm \times \left[ 1 + \sum_{j = j_0}^{j_1} \exp \left\{  O \left( j \cdot \frac{p k^3}{n^2} \right) - \Omega \left( j^2 \cdot \frac{1}{nip}   \right) \right\} \right. \\
& \hskip6cm \left.  + \sum_{j \ge j_1} \exp \left\{  O \left( j \cdot \frac{p k^3}{n^2} \right) - \Omega \left( j \right) \right\} 
\right].
\end{split}
\end{equation*}
The desired bound follows (using $ k < n^{3/4 - \epsilon}$ and $i = O(1/p)$ to conclude that $ j_0/ (nip) \gg p k^3/n^2$).

\subsubsection{Poisson approximation in a bipartite random graph}

In this Section we prove a Lemma that allows us to bound the 
probability of the event $ \cF_C$ defined above.
In order to make a general presentation, we consider
a random bipartite graph $ G_{N,N,p}$ on vertex set $V$. 
We should think of $p$ as being the $p'$ defined above and 
$N$ as being approximately equal to $k$.
We are looking to get a Poisson-type upper bound on the probability that all 
vertices of $ G_{N,N,p}$ outside some exceptional set $T$  
have degree at least 2.  Set $ t = |T|$. In our application of this Lemma, the vertex set
$ V $ has bipartition $ A \setminus B, B \setminus A$ and $T = V \cap \left[ (\cup_{e \in M_A}e) \cup(\cup_{e \in M_B} e) \right]$. For each vertex $v$ let $ J_v$ be 
the indicator random variable for the event that the degree of $v$ is 0 or 1. Furthermore, let
\[ X = \sum_{v \in V \setminus T} J_v.  \]
Let
\begin{align}
\mu & = \EE[X]  =  (2N - t) ( (1-p)^N + Np ( 1-p)^{N-1} ) \ \  \text{ and }  \label{eq:mu} \\
\sigma & =   pN  ( (1-p)^{N-1} + (N-1)p ( 1-p)^{N-2} )  \label{eq:sigma}
\end{align}
We are looking to establish a bound of the form
$ \PP(X = 0) \le e^{ -\mu }$. We achieve this following an argument recently
introduced by Bohman, Warnke and Zhu \cite{bwz}. This argument is an adaptation of
the original proof of Janson's inequality.
\begin{lem}
\label{lem:withLutz}
Let $V$ be the vertex set of the random bipartite graph $ G_{N,N,p}$.
If $ T\subset V$ is a set of $t$ vertices, and
$X$ is the number of vertices in $ V \setminus T$ that have degree at most 1, then
we have
    \[ \PP( X =0 ) \le \exp\left\{ - \mu  + \mu \sigma \log(1 + 1/\sigma) \right\}, \]
where $ \mu$ and $\sigma$ are defined in (\ref{eq:mu}) and (\ref{eq:sigma}), respectively.
\end{lem}
\noindent
Note that in our application of this Lemma we have $ i = O(1/p)$ and
therefore
\[
\mu = O \left( \frac{ p k^4}{ n^2} \right) \ \ \ \text{ and } \ \ \ 
\sigma = O \left( \frac{ p^2 k^4}{ n^2} \right).
\]
So, it follows that in our application we have
\begin{equation}
\label{eq:bipartite conclusion}
    \PP( X =0) \le \exp \left\{ - \mu + O\left( \frac{ (\log n) p^3 k^8}{ n^4}   \right) \right\}.
\end{equation}

\begin{proof}
For ease of notation set $ V' = V \setminus T$.
We make use of the moment generating function.
Define $ f(\lambda) = \EE[ e^{-\lambda X} ]$. Note that we have 
$ \PP( X = 0) \le f( \lambda)$ for all $ \lambda$. Therefore, as $ f(0)=1$, we have
\[ - \log( \PP(X=0)) \ge  - \log ( f( \lambda)) = - \int_0^\lambda \left( \log (f(s)) \right)' ds = -\int_{0}^\lambda \frac{ f'(s)}{ f(s)}ds.        \]
Note that if we let $ \Omega $ be our probability space (i.e. the collection of all bipartite graphs with the binomial measure) then we
have
\begin{multline*}
  f'(s) = \frac{d}{ds} \sum_{ \omega \in \Omega} \PP( \omega) e^{-s X(\omega)} = - \sum_{ \omega \in \Omega} \PP( \omega) X( \omega)  e^{-s X(\omega)} \\
  = - \sum_{ \omega \in \Omega} \PP( \omega)   e^{-s X(\omega)}   \sum_{v \in V'} J_v( \omega) =
 - \sum_{ v \in V'} \sum_{ \omega \in \Omega}   \PP( \omega)   e^{-s X(\omega)}  J_v( \omega) = - \sum_{v \in V'} \EE[ J_v e^{-sX}]
\end{multline*} 
Therefore,
\[ - \log ( \PP( X=0)) \ge \int_{0}^\lambda \sum_{v \in V'} \frac{ \EE[ J_v e^{-sX}] }{ \EE [e^{-sX}]} ds.\]
So, we would like to find a
lower bound on
\begin{equation}
\label{eq:goal}
\sum_{v \in V'} \frac{ \EE[ J_v e^{-sX}] }{ \EE[ e^{-sX}]}.
\end{equation}

Now we define an auxiliary collection of random variables. For distinct 
vertices $u,v$ let $ J_u^v$ be the indicator random variable for the
event that $u$ has degree 0 or 1 when we ignore the status of the potential edge $uv$. 
If $u$ and $v$ are in the same part in the bipartition then we simply set $ J_u^v = J_u$.
Define 
\[ X_v = \sum_{u \in V': u \neq v} J_u^v. \]
Note that we have
\[ \EE[ J_v e^{-sX}] = \PP( J_v =1) \EE[ e^{-sX} \mid J_v =1].   \]
Furthermore if we condition on the event that $ J_v =1$ then $ X = X_v +1$ in almost all cases. The only exception is the situation in which $ v $ has a unique neighbor $v' \in V'$ 
and $ J_{v'}^v =1$ while $ J_{v'}=0$; in other words, $v$ has unique neighbor $v' \in V'$ and this neighbor has degree 2 in the graph. Note that in this special case we have $ X = X_v$. It follows that when we condition on $ J_v=1$ we have $ X \le X_v +1$ and, as the random variables $J_v$ and $ X_v$ are independent, we have
\[
\EE[ e^{-sX} \mid J_v=1 ] \ge  \EE[ e^{-s(X_v + 1)} \mid J_v=1 ] = e^{-s} \EE[ e^{-sX_v} \mid J_v=1 ] 
= e^{-s} \EE[ e^{-sX_v}]
\]

In order to get a lower bound on the expression in (\ref{eq:goal}) we now look to bound $ \EE[e^{sX}]$ in
terms of $ \EE[e^{sX_v} ]$. To this end, define
\[ Y_v = \sum_{u \in \Gamma(v) \setminus T} J_u^v,\]
where $ \Gamma(v)$ is the neighborhood of $ v$ in $ G_{N,N,p}$.
Note that we have 
\[ X \ge J_v + X_v - Y_v \ge X_v - Y_v.\]
Therefore
\[ \EE[ e^{-sX}] \le \EE[ e^{-sX_v + s Y_v}] = \EE\left[ \EE[ e^{-sX_v} e^{sY_v} \mid \Gamma_v ] \right].\] 
Now we note that conditioning on the event $ \Gamma_v = S $ yields a product space 
on the $N(N-1)$ potential edges that do not include $v$. In this space $ e^{-sX_v}$ is increasing 
while $ e^{s Y_v}$ is decreasing. Thus, we can apply FKG to conclude
\[ \EE[ e^{-sX_v} e^{sY_v} \mid \Gamma_v ]  \le \EE[ e^{-sX_v} \mid \Gamma_v ] \EE[ e^{ sY_v} \mid \Gamma_v] = \EE[ e^{-sX_v}] \EE[ e^{ sY_v} \mid \Gamma_v]. \]
and therefore
\[ \EE[ e^{-sX} ] \le \EE[ \EE[ e^{-sX_v} ] \EE[ e^{sY_v} \mid \Gamma_v]] =\EE[ e^{-sX_v} ]  \EE [\EE[ e^{sY_v} \mid \Gamma_v]] = \EE[ e^{-sX_v} ] \EE[ e^{sY_v} ].  \]
Now, for $u,v$ in different parts of the bipartite graph, set 
\[ \Pi = \PP( J_u^v =1) =  (1-p)^{N-1} + (N-1)p (1-p)^{N-2}.\]
As the events $ \left\{ \{ J_u^v =1\} : u \in \Gamma_v \right\}$ are independent, we have
\[ \EE[ e^{ s Y_v} \mid \Gamma_v =S ]= \left( 1 - \Pi + e^{s} \Pi \right)^{ |S \setminus T|}  = ( 1 + \Pi( e^s -1))^{|S \setminus T|}.\]
Therefore, letting $V_2$ be the part in the bipartition that does not contain $v$ and setting $ V_2' = V_2 \setminus T$, we have 
\begin{multline*}
\EE[ e^{s Y_v}] = \sum_{S \subset V_2 } \PP( \Gamma_v = S) \EE[ e^{sY_v} \mid \Gamma_v = S]
= \sum_{ k=0}^{ |V_2'|} \binom{|V_2'|}{k} p^k (1-p)^{|V_2'|-k} ( 1 + \Pi( e^s -1))^{k} \\
= ( 1 - p + p ( 1 + \Pi( e^s -1)))^{|V_2'|} = ( 1 + p \Pi ( e^s -1))^{|V_2'|}
\le \exp \left\{ p \Pi |V_2'| ( e^s -1) \right\}.
\end{multline*}

Putting everything together we have
\[ -\log( \PP(X=0)) \ge \int_{0}^\lambda \mu \exp\left \{ -s - p \Pi N ( e^s -1) \right\} ds.\]
In order to estimate the integral we use the bound $ e^{-x} \ge 1 -x $ and
recall $ \sigma = p N \Pi$ to write
\begin{multline*}
\int_{0}^\lambda \mu \exp\left \{ -s - p \Pi N ( e^s -1) \right\} ds
\ge \mu \int_{0}^\lambda  e^{-s}  \left( 1 - \sigma (e^s -1) \right) ds \\
= \mu \int_0^\lambda - \sigma + e^{-s}( 1 + \sigma) ds
= \mu \left[  - \sigma \lambda - ( e^{-\lambda}- 1)( 1 + \sigma ) \right].
\end{multline*}
This expression is maximized at $ \lambda = \log ( 1 + 1/\sigma)$, and we have
\[-\log( \PP(X=0)) \ge \mu \left[ 1 - \sigma \log( 1 + 1/\sigma) \right]. \]
\end{proof}

\section{Proof of Theorem~\ref{thm:main}} \label{sec:WXYZ}

We begin by proving Theorem~\ref{thm:main} for $m > n^{4/3+ \epsilon}$. Here we will apply a result recently established by the authors.
\begin{thm}[Bohman, Hofstad \cite{bhindy}]
\label{thm:oldy}
    Let $ \epsilon>0$. If $ n^{-2/3 + \epsilon} < p \le 1$ then $ \alpha(G_{n,p})$ is concentrated on 2 values.
\end{thm}
We also use the following Lemma, the proof of which is left to the reader.
 \begin{lem} \label{lem: setup for m large}
    Let $Q_n$ be a sequence of random variables that take values in the integers. The sequence $Q_n$ is {\bf not} concentrated on 2 values if and only if there exists some sequence $n_1 <n_2 <\dots$, some $\varepsilon > 0$, and a sequence  $k_1, k_2, \dots $ such that
\[ \PP[Q_{n_i} < k_i] \geq \varepsilon \ \ \ \ \text{ and } \ \ \ \ \PP[Q_{n_i} > k_i] \geq \varepsilon. \]
for all $i \in \ZZ^+.$
\end{lem}

\begin{thm}
    If $ n^{4/3 +\ep} < m \leq \binom{n}{2}$ then $\alpha(G_{n,m})$ is two-point concentrated.
\end{thm}
\begin{proof}
Assume for the sake of contradiction that there exists some $m = m(n) > n^{4/3 + \ep}$ such that $\alpha(G_{n,m})$ is not two-point concentrated. 
By Lemma~\ref{lem: setup for m large}, there exists an increasing sequence $n_1,n_2,\dots$, some $\varepsilon > 0$, and some  $k_i$ such that
\begin{equation*}
    \PP[\alpha(G_{n_i,m(n_i)}) < k_i] >\varepsilon \qquad \text{and} \qquad \PP[\alpha(G_{n_i,m(n_i)}) > k_i] >\varepsilon \ \ \ \ \forall i \in \ZZ^+.
\end{equation*}
Now choose a probability $ p_i $ so $ m(n_i)$ is the median of $ e( G_{n_i,p_i})$ (where $ e( G)$ is the number of edges in the graph $G$). Since the independence number is a decreasing function with respect to edges, then
    \[
        \Pr[\alpha(G_{n_i,p_i}) < k_i] \geq 
        \Pr[\alpha(G_{n_i,p_i}) < k_i | e(G_{n_i,p_i}) \geq m(n_i)] \Pr[e(G_{n_i,p_i}) \geq m( n_i)] \geq
        \varepsilon/2.
    \]
    Similarly, $\Pr[\alpha(G_{n_i,p_i}) > k_i] \geq \varepsilon/2$. Therefore, by Lemma \ref{lem: setup for m large}, $\alpha(G_{n,p})$ is not two-point concentrated, which contradicts Theorem~\ref{thm:oldy}.
\end{proof}

In the remainder of this Section we assume $ n^{5/4+ \epsilon} < m < n^{4/3 + \epsilon}$, and prove the following two Lemmas. Theorem~\ref{thm:main} for $m$ in this range follows immediately. These Lemmas are also used in the proof of Theorem~\ref{thm:Gnp}, which appears in the next section.
\begin{lem} \label{lem:lower}
    Let $ \epsilon>0$ and $ n^{5/4+ \epsilon} < m < n^{4/3 + \epsilon}$.
    If $ r \le r_1$ and $k$ satisfies (\ref{equ: k used in this range}) and $\EE[ X_{k,r}] = \omega(1)$ 
    then $ \alpha( G_{n,m}) \ge k$ whp.
\end{lem}
\begin{lem} \label{lem:upper}
    Let $ \epsilon>0$ and $ n^{5/4+ \epsilon} < m < n^{4/3 + \epsilon}$.
    If $k$ satisfies (\ref{equ: k used in this range}) and $ \EE[ X_{k, r_0(k)}] < n^{ 2\epsilon}$ then $ \alpha( G_{n,m}) \le k
    $ whp.
\end{lem}
\noindent
Theorem~\ref{thm:main} follows immediately from Lemmas~\ref{lem:lower}~and~\ref{lem:upper}. We state this as a separate Lemma for use in the proof of Theorem~\ref{thm:Gnp} below. 
Recall that $k_0$ is
defined to be the smallest $k$ such that $ \EE[X_{k,r}] < n^{\epsilon}$ for all $r \leq r_1$. 
\begin{lem} \label{lem:picky}
Let $ \epsilon>0$ and $ n^{5/4+ \epsilon} < m < n^{4/3 + \epsilon}$. We have
$ \alpha(G_{n,m}) \in \{ k_0-1, k_0 \}$ whp.
\end{lem}
\begin{proof}
Note that we have $ n^{\epsilon} < \EE[X_{k_0-1, r_0(k_0-1)}]$ by the definition of $k_0$. So Lemma~\ref{lem:lower} implies $ \alpha( G_{n,m} ) \ge k_0 -1$ whp.
As we have $ \EE[X_{k_0,  r_0(k_0)}] < n^{\epsilon} < n^{2\epsilon}$, Lemma~\ref{lem:upper} implies $ \alpha( G_{n,m}) \le k_0 $ whp. 
\end{proof}

We now turn to the proof of Lemmas~\ref{lem:lower}~and~\ref{lem:upper}.
Recall that $Y_{k,r}$ counts the number of extended independent sets while Theorems~\ref{thm:techy1}~and~\ref{thm:techy2}, which are the main technical result in this work, addresses
the closely related variable $ X_{k,r}$. In the remainder of this Section we translate these results
for $ X_{k,r}$ to $ Y_{k,r}$ and prove Lemmas~\ref{lem:lower}~and~\ref{lem:upper}.

In order to facilitate the comparison of $ Y_{k,r}$ and $ X_{k,r}$, we introduce three
more random variables, i.e. we consider a collection of five closely
related random variables. These variables address the two differences between extended independent
sets and the pairs $ K,M$ counted by $ X_{k,r} $: the appearance of cycles within $K$ and allowing
vertices outside of $K$ to satisfy the minimum degree 2 requirement in a connected component in $G[K]$. Indeed
we define one random variable for each combination of allowing or prohibiting cycles in the interior, as well as allowing for or prohibiting neighbors of the same interior component to satisfy the 2 neighbors requirement for vertices outside of $K$. We also introduce a variable that does not depend at all of the edges that are not contained in $K$.
\begin{define} \label{def: wxyz} Let $U_{k,r} / W_{k,r}/X_{k,r}/Y_{k,r}/Z_{k,r}$ be the random variable that counts the number of pairs $K, M$ where $K$ is a set of $k+r$ vertices in $G=G_{n,m}$ and $M$ is a perfect matching consisting of $r$ edges 
contained within $K$ and
\begin{enumerate}
    \item $G[K]$ contains the matching $M$, and every edge in $G[K]$ has both vertices in $ \cup_{e \in M} e$.
    \item For variables $W_{k,r},X_{k,r},Y_{k,r},Z_{k,r}$, and for every vertex $v \not\in K$, the set $\mathcal{N}(v) \cap K$ contains: 
    \begin{itemize}
        \item ($W_{k,r}$ and $X_{k,r}$ only) At least two isolated vertices of $G[K]$.
        \item ($Y_{k,r}$ and $Z_{k,r}$ only) At least two isolated vertices of $G[K]$ OR at least two vertices of the same connected component of $G[K]$.
    \end{itemize}
    \item ($W_{k,r}$ and $Y_{k,r}$ only) $K$ contains no cycles.
\end{enumerate}
\end{define}
Below is a table that displays the differences between the five variables: 
\begin{center}
\begin{tabular}{ c | c | c |}
    & No cycles & Possible cycles \\
    & in interior & in interior \\ 
   \hline
   2-neighbors requirement & & \\
   must come from & $W_{k,r}$ & $X_{k,r}$ \\
   ``isolated" interior vertices & & \\ \hline
   2-neighbors requirement & & \\
   {\it could} come from & $Y_{k,r}$ & $Z_{k,r}$ \\ 
   common interior component & & \\ \hline
   & & \\
   No 2-neighbors requirement &  & $U_{k,r}$ \\
    & & \\ \hline
\end{tabular}
\end{center}
 
 \noindent
 Observe that $W_{k,r} \leq X_{k,r} \leq Z_{k,r}$ and $W_{k,r} \leq Y_{k,r} \leq Z_{k,r}$ and $Z_{k,r} \leq U_{k,r}$; in words, the further up and left on the table you are, the more strict the requirements become for a pair $ K,M$ that the corresponding variable counts. We begin by proving two Lemmas that bound the relationships among variables $W, X, Y, $ and $Z$. 
 \begin{lem} 
\label{lem:W} If $ r \le r_1$ and $k$ satisfies (\ref{equ: k used in this range}) then
     \[\EE[W_{k,r}] = (1-o(1))\EE[X_{k,r}]. \]
 \end{lem}

 \begin{lem}  \label{lem:Z} If $ r \le r_1$ and $k$ satisfies (\ref{equ: k used in this range}) then
     \[ \EE[ Z_{k,r} ] = (1 +o(1)) \EE[ X_{k,r}].\]
 \end{lem}

\begin{proof}[Proof of Lemma~\ref{lem:W}]
Let a set $K$ of $ k+r $ vertices and a matching $M$ in $K$ with $r$ edges be fixed.
Let $\cU$ 
be the event that $ G[K]$ contains the matching $M$ and
no edge that intersects $ K \setminus ( \cup_{e \in M} e)$.
Let $\mathcal{X}$ be the event that the pair $K, M$ is counted by the variable $X_{k,r}$; 
note that $\PP[\mathcal{X} ] =  \PP[ \cU ] \Phi$. 

Now, the pair $K,M$ is also counted by $ W_{k,r}$ unless a cycle appears within $K$.
To bound the probability of this event,
we bound the probability that any specific cycle inside $K$ appears 
conditioned on $\mathcal{X}$ and then apply the union bound. Note that we can
determine whether or not the event $\mathcal{X}$ holds by viewing $ O(nk) = o(n^2)$ pairs of vertices 
(to determine whether or not each is an edge or not). Under such
conditioning, the probability that a particular edge appears in $ G_{n,m}$ is at
most $ (1 + o(1))p$ (recall $ p = m \binom{n}{2}^{-1}$). Then a specific cycle in $K$ with 
$\ell$ non-matching edges appears with probability at most $\left(p+o(p) \right)^{\ell} < (2p)^{\ell}$. 
The number of cycles in $K$ with $\ell$ non-matching edges can bounded above by $(4r)^{\ell}$ (choose, in order, the matching edges the cycle will ``visit", allowing for possible repeats, and for each non-matching edge in the cycle, determine which of the two vertices of the matching edges will be the start and end vertices of the non-matching cycle edge). Hence, by the union bound,
\begin{equation*}
    \PP[G[K] \text{ contains a cycle} | \mathcal{X}] \leq \sum_{\ell = 2}^{\infty} (8rp)^{\ell} = O \left( \left( \frac{ p^2 k^3}{n}  \right)^2 \right).
\end{equation*}
\end{proof}

\begin{proof} [Proof of Lemma~\ref{lem:Z}]
    Let $K$ be a fixed set of $k+r$ vertices and let $M$ be a matching of size $r$ in $K$.
    Let $\mathcal{Z}$ be the event that the pair $ K,M$ is counted by the variable $Z_{k,r}$. And again let
    $\mathcal{X}$ be the event that $ K,M$ is counted by the random variable $ X_{k,r}$. Note that 
    $\mathcal{X} \subset \mathcal{Z} \subset \mathcal{U}$. It suffices to show 
    $\Pr[\mathcal{Z}] \leq (1+o(1)) \Pr[\mathcal{X}]$. 

Let $S$ be some set of vertices outside $K$ such that $|S| = s$ (we will eventually enumerate over all $2^{n-k-r}$ possible choices of $S$). Define $\mathcal{S} = \cS_S$ to be the event that each vertex in $S$ has at most one neighbor in $K \setminus( \cup_{e \in M} e )$ and at least two neighbors in $ \cup_{e \in M} e$, while the rest of the vertices outside $K$ have at least two neighbors in $K \setminus ( \cup_{e \in M} e )$. 
Note that $ \cX = \cU \wedge \cS_\emptyset$ and
\[ \cZ \subsetneq \bigcup_{S \subseteq V \setminus K} \cU \wedge \cS_S. \] 
We now bound $\PP[\mathcal{U} \wedge \mathcal{S} ]$ by exposing edges in batches through the following process: 
    \begin{itemize}
        \item[0.] View vertex pairs in $ \left( \binom{K}{2} \setminus \binom{ \cup_{e \in M}e}{2} \right) \cup M$ to determine the event $\mathcal{U}$.

        \item[1.] View only as many vertex pairs as necessary between $V \backslash{(K \cup S)}$ and $K \setminus (\cup_{e \in M} e)$ to verify that all vertices in $V \setminus {(K \cup S)}$ have at least two neighbors in $K \setminus (\cup_{e \in M} e)$.

        \item[2.] View {\it all} pairs of vertices between $S$ and $K \setminus ( \cup_{ e\in M} e)$ to verify that the vertices in $S$ have at most one neighbor in $K \setminus( \cup_{e \in M} e)$.

        \item[3.] View vertex pairs between $S$ and $\cup_{e \in M}e$ to verify that each vertex in $S$ has at least two neighbors in $ \cup_{e \in M} e$.
    \end{itemize}

    We now consider the probability that these steps result in an outcome in the event $\mathcal{U} \wedge \cS$. For step 1, we pass to the binomial model as set forth in Lemma~\ref{lem: passing to Gnp}. In particular, we compare with $ G_{n,p'}$ with $ p' = p + (\log n) \max\{ p k^2/n^2 , p^{1/2}/n \}$. Let $ \cE$ be the event that every vertex outside of $K \cup S $ as at least two neighbors in $ K \setminus( \cup_{e \in M}e)$. 

We have $c_0 - c_{p'} = O( (\log n) ( p k^3/n^2 + p^{1/2}k/n ) \ll \frac{ np^2}{ (\log n)^3}$.
Thus, by 
Lemmas \ref{lem:cs} and \ref{lem:Phi} we have
\begin{equation}
\label{eq:step 1}
\PP( \cE \mid \cU) \le(1 +o(1)) \PP_{p'}( \cE) \sim \varphi(c_{p'})^{n-k-r - s} \sim \varphi( c_0)^{n-k-r-s} \sim \frac{ \Phi}{ \varphi( c_0)^s}. 
\end{equation}
    
    For step 2, we note that by the end of step 2 we reveal at most $O(n)$ edges and $O(nk)$ non-edges, hence at each pair of vertices we view the probability of viewing an edge is both $(1 - O(n/m))p = (1 - O(1/np))p$ and $(1+O(k/n))p$. Hence, the probability the second step is completed in accordance with the event $\cZ$ is 
    \begin{equation}
        \label{eq:step 2}
    \Big((1+o(1)) kp \big(1-(1- O(1/np))p\big)^{k-1}\Big)^{s} \le \Big( 2 k^3 p n^{-2} e^{-2}\Big)^s.
    \end{equation}    
    For the final step, we again note that we view $O(nk)$ pairs of vertices, so the probability of uncovering an edge for a given pair of vertices is at most $(1+O(k/n))p.$ As there are $2r$ vertices in $ \cup_{e \in M} e $, the
    probability every vertex in $S$ has the required neighbors in $ \cup_{e \in M} e$ is at most
    \[ \left( 3 r^2 p^2 \right)^s.  \]
    
    Putting this observation together with (\ref{eq:step 1}) and (\ref{eq:step 2}) we have
    \begin{equation}
        \PP[ \mathcal{S} | \mathcal{U}] = \frac{\Phi}{ \varphi(c_0)^s} \left(O(1) \cdot \frac{ k^3 p}{ n^{2}} \cdot r^2p^2 \right)^s = \Phi \left(O(1) \frac{ k^9 p^5}{ n^4} \right)^s \label{equ: bounding Z event}.
    \end{equation}
    We now complete the proof:
       \begin{multline*}
     \PP( \cZ \mid \cU)
      \le \sum_{S \subseteq V \setminus K} \PP( \cS_S \mid \cU) 
     \le  \Phi \sum_{s=0}^{n-k-r} \binom{n-k-r}{s} \left(O(1) \frac{ k^9 p^5}{ n^4} \right)^s \\
         \le \Phi \sum_{s=0}^{n-k-r}\left(O(1) \frac{ k^9p^5}{n^3} \right)^s 
          = ( 1 +o(1)) \Phi = (1 +o(1)) \PP( \cX \mid \cU).
    \end{multline*}
 \end{proof}

With Lemma~\ref{lem:W}~and~\ref{lem:Z} in hand we now turn to 
the proofs of Lemma~\ref{lem:lower}~and~\ref{lem:upper}.

\begin{proof}[Proof of Lemma~\ref{lem:lower}] 
Note that $ \EE[W_{k,r}] \to \infty$ by Lemma~\ref{lem:W}.
As $ W_{k,r} \le X_{k,r}$, and applying Theorem~\ref{thm:techy1} and Lemma~\ref{lem:W}, we have
\[ \EE[ W_{k,r}^2] \le \EE[X_{k,r}^2] \le ( 1 +o(1)) \EE[X_{k,r}]^2 \le (1 +o(1)) \EE[W_{k,r}]^2.\]
Therefore, $ \Var[W_{k,r}] = o( \EE[W_{k,r}]^2)$ and by the second moment method we have
$ W_{k,r}>0$ whp. As $ W_{k,r} \le Y_{k,r}$, we have $ Y_{k,r} >0$ whp. The Lemma now follows from Lemma~\ref{Lem:MNI connection}.
\end{proof}

\begin{proof}[Proof of Lemma~\ref{lem:upper}]
By
Lemma~\ref{Lem:MNI connection}, and the fact that whp $\alpha(G_{n,m}) \leq k_V$ by Lemma~\ref{lem: k_v markov} and Markov's Inequality, it suffices to show
\[ \sum_{\hk = k+1}^{k_V}\sum_{r=0}^{ \hk} \EE \left[ Y_{\hk,r} \right] = o(1). \]
For the remainder of this
Section we always work with $\hk$ between $k+1$ and $k_V$ inclusive and hence we are assuming
that $\hk$ satisfies (\ref{equ: k used in this range}). For ease of notation we 
set $ r_0 = r_0(k)$, and recall that we set $ r_1 = 8 ( \log np)^3/( n p^2)$, and $ r_1 \sim pk^3/n \sim 2e^2 r_0$.
We begin by noting that for any $r \le r_1 $, using (\ref{equ: X' asymptotic formula}) and iterated application of (\ref{equ: X' ratio over k}) gives
\[ \sum_{ \hk = k+1}^{k_V} \EE[ X_{\hk,r}] \le   \EE[ X_{k+1,r}] \sum_{ j=0}^{k_V - k -1} \left( \frac{k}{n} \right)^j = (1 +o(1))  \EE[ X_{k+1,r}] .\]
It then follows from 
Lemma~\ref{lem:Z}  and Theorem~\ref{thm:techy2} and (\ref{equ: X' asymptotic formula}) and (\ref{equ: X' ratio over k}) that we have
\begin{multline*}
\sum_{\hk = k+1}^{k_V}\sum_{r=0}^{ r_1}  \EE[ Y_{\hk,r}] \le \sum_{\hk = k+1}^{k_V} \sum_{r=0}^{r_1} \EE[ Z_{\hk,r}]
\le \sum_{\hk = k+1}^{k_V} \sum_{r=0}^{r_1} (1 +o(1)) \EE[ X_{\hk,r}] \\\leq
2 \sum_{r=0}^{r_1} \EE[ X_{k+1,r}] = O\left(  \sqrt{ r_0} \cdot \EE[ X_{k+1, r_0}] \right)
= O \left( \frac{ p^{1/2} k^{3/2}}{ n^{1/2}} \cdot \frac{ k}{ n}  \cdot \EE[ X_{k, r_0}]  \right) = o(1),
\end{multline*}
using the assumption $ m > n^{5/4 + \epsilon}$ (which implies $ k = \tilde{ O}( n^{3/4 - \epsilon})$).

Next we consider $ r$ such that $ r_1 < r < 1/p$. Here it suffices to bound $ \EE[ Y_{\hk,r}] $ using only bounds on $ \EE[U_{\hk,r}]$. We begin
by noting that we can infer a bound on $ \EE[ U_{k+1, r_0}]$ from the condition on $ \EE[ X_{k,r_0}]$. Indeed, we have 
$ \varphi(m,k,r_0)^{n-k-r_0} = \exp \left\{ - (1 +o(1)) npk \cdot ( k e^{-1} n^{-1} )^2 \right\} 
= e^{ - (2 +o(1)) r_0} $. Furthermore,
\[ \EE[ U_{k +1, r_0}] \varphi(m,k+1, r_0)^{n-k-1-r} \sim \EE[ X_{k +1, r_0}] < \EE[ X_{k, r_0}] \le n^{2\epsilon}. \]
Combining these two observations we have
\begin{equation*} 
\EE[ U_{k +1, r_0}] \le e^{ (2 +o(1)) r_0}.
\end{equation*}
Next we observe (\ref{equ: N ratio over r}) and (\ref{equ: U ratio over r}) imply that for any 
$ r $ we have
\[ \frac{ \EE[ U_{\hk, r+1}] }{ \EE[ U_{ \hk, r}]} = (1+o(1)) \frac{n \hk }{2r } \cdot e^{3rp} p \cdot \frac{ \hk^2}{e^2 n^2} = ( 1 +o(1)) \frac{ p\hk^3}{ 2 e^2 n} \cdot \frac{ e^{3pr}}{r}.  \]
It follows that  if $ r \ge r_0$ then $ \EE[U_{\hk, r+1}] / \EE[ U_{ \hk, r}] \le 1+o(1)$ and 
\begin{equation} 
\label{eq: U step r}
3r_0 < r < 1/p \ \ \ \Rightarrow \ \ \ \frac{ \EE[ U_{\hk, r+1}] }{ \EE[ U_{ \hk, r}]} \le \frac{1}{e}. 
\end{equation}
Therefore, 
\[ \EE[ U_{k+1, r_1} ] = \EE[ U_{k+1, r_0}] \prod_{ j = r_0}^{r_1-1} \frac{ \EE[ U_{k+1, r+1}] }{ \EE[ U_{ k+1, r}]} \le e^{ (2 +o(1)) r_0} \left( \frac{1}{e} \right)^{r_1 - 3r_0} \le e^{ - 5 r_0}. \]
Next, we observe that by (\ref{equ: N ratio over k}) and (\ref{equ: U ratio over k}) we have
\begin{equation} \label{eq: U step k}  \frac{ \EE[ U_{\hk+1,r_1}]}{ \EE[U_{\hk,r_1}]} \le \frac{ k}{en}.   \end{equation}
Iterative application of (\ref{eq: U step r}) and (\ref{eq: U step k}) then gives
\begin{equation*}
    \sum_{\hk = k+1}^{k_V}\sum_{r=r_1}^{1/p}  \EE[ Y_{\hk,r}] \leq \sum_{\hk = k+1}^{k_V}\sum_{r=r_1}^{1/p} \EE[ U_{\hk,r}] \leq 2 \sum_{\hk = k+1}^{k_V}  \EE[ U_{\hk,r_1}] = O \left( \EE[ U_{k +1, r_1}] \right) = o(1).
\end{equation*}

For $r > 1/p$ we once again do not make use of the degree 2 condition for vertices outside the set $K$. However, we do make use of the condition on edges within $ \cup_{e \in M} e$ in the definition of an extended independent set. In particular, we use the fact that this collection of edges induces no cycle. This implies that the number of edges in this set is at most $ 2r-1$. We have 
\begin{equation} \label{eq:Y upper large r} \EE[ Y_{\hk,r}] \le  \sum_{j=0}^{r-1} \binom{n}{\hk+r} \binom{ \binom{ \hk+r}{2}}{ r +j}  \frac{ \binom{ \binom{n}{2} - \binom{\hk+r}{2}}{m -r -j}}{ \binom{ \binom{n}{2}}{m}  }.   
\end{equation}
Next we note that we have
\begin{multline*}
  \frac{ \binom{ \binom{n}{2} - \binom{\hk+r}{2}}{m -r -j}}{ \binom{ \binom{n}{2}}{m}  } = \frac{ (m)_{r+j}}{ \left( \binom{ n}{2} \right)_{r+j}} 
\cdot \frac{ \left( \binom{n}{2} - \binom{ \hk+r}{2} \right)_{m-r-j}}{ \left( \binom{ n}{2} -r -j \right)_{m-r-j}} 
\le  \left( \frac{ m}{ \binom{ n}{2} } \right)^{ r+j} 
\cdot \left(  \frac{ \binom{n}{2} - \binom{ \hk+r}{2}} { \binom{ n}{2} -r -j} \right)^{m-r-j} \\
= p^{r+j} \left( 1 - \frac{ \binom{ \hk+r}{2}}{ \binom{n}{2}}  \right)^{m-r-j} 
\left( 1 + \frac{ r+j}{ \binom{n}{2} -r -j}  \right)^{m-r-j}  \\
\le p^{r+j} \exp\left\{ - (m-r-j) \frac{\binom{\hk+r}{2}}{ \binom{n}{2}} + 8 \frac{ r m }{ n^2} \right\} \le p^{r +j} \exp \left\{ - p \binom{\hk+r}{2} +o(r) \right\}.
\end{multline*}
We now regroup the terms.
We have
\begin{multline*}
\binom{n}{\hk+r} e^{-p\binom{\hk+r}{2}} \le \exp \left\{ (\hk+r) \left[ \log \left( \frac{ ne}{ \hk+r} \right) - \frac{ p(\hk+r -1)}{2}   \right]    \right\}  \\
= \exp \left\{(\hk+r)\left[ \log \left( \frac{\hk}{\hk+r} \right) - \frac{pr}{2} +o(1) \right]   \right\} = e^{ - \Omega\left( (\log n)  \cdot r \right)},
\end{multline*}
and
\[  \binom{ \binom{\hk+r}{2}}{ r+j} p^{r+j} \le \left( \frac{ (\hk+r)^2ep}{2(r+j)}   \right)^{r+j}  \le \left( 2e (\log n)^2 \right)^{r+j} = e^{ O \left( \log\log n \right) \cdot r }  \]
As there are only a polynomial number of terms in the sum in (\ref{eq:Y upper large r}) we
conclude $ \EE[ Y_{\hk,r}] \le e^{ - \Omega\left( (\log n) \right) \cdot r } $. It follows that we have
\[ \sum_{\hk = k+1}^{k_V}\sum_{r = 1/p}^{\hk} \EE[ Y_{\hk,r} ] = o(1),\]
as desired.
\end{proof}

\section{Regarding $ G_{n,p}$: Proof of Theorem~\ref{thm:Gnp}.} \label{sec:G(n,p)}

In this Section we prove Theorem~\ref{thm:Gnp} as a Corollary of our results on
$G_{n,m}$.
We begin by noting how $X'(m,k,r)$ changes as we change the number of edges.
\begin{lem} \label{lem: X' ratio over m} If $k$ satisfies (\ref{equ: k used in this range}) and $r \leq r_1$ we have
    \begin{equation}
       \frac{X'(m+1,k,r)}{X'(m,k,r)} = 1 - (1\pm o(1))\frac{k^2}{n^2}. \label{equ: X' ratio over m}
    \end{equation}
\end{lem}
\begin{proof}
Note that $ N(k,r)$ does not depend on $m$. We consider
the changes in $ U(m,k,r)$ and $ \varphi(c(m,k,r))^{n-k-r}$ that result from incrementing $m$ by 1 in turn.
We have
\begin{equation*}
        \frac{ \binom{\binom{n}{2} - M_1}{m+1-r}}{\binom{ \binom{n}{2} - M_1 }{m-r}} = \frac{ \binom{n}{2} - M_1 - m + r}{m+1-r} \qquad \text{and} \qquad \frac{\binom{\binom{n}{2}}{m}}{\binom{\binom{n}{2}}{m+1}} = \frac{m+1}{\binom{n}{2} - m},
\end{equation*}
hence, we have
\begin{align*}
    \frac{U(m+1,k,r)}{U(m,k,r)} &= \left(\frac{ \binom{n}{2} - M_1 - m + r}{\binom{n}{2} - m}\right)\left(\frac{m+1}{m+1-r}\right) \\&=
    1 - (1 \pm o(1)) \frac{k^2}{n^2} + O\left(\frac{r_1}{m}\right) \\&=
    1 - (1 \pm o(1)) \frac{k^2}{n^2}.
\end{align*}
Next recall that we have
$\frac{d}{dc} \log(\varphi(c)) = c e^{-c}$, hence for $c \gg 1$ and any $\Delta c$ at most a constant we have
\begin{equation*}
    \frac{\varphi(c+\Delta c)^{n-k-r}}{\varphi(c)^{n-k-r}} = \exp\left\{ O(\Delta c \cdot
    (n-k-r) c e^{-c})\right\}.
\end{equation*}
Noting that increasing $m$ by 1 results in an increase of $ (k-r)/n^2$ in $c$, it follows that
we have
\begin{equation*}
     \left( \frac{\varphi(c(m+1,k,r))}{\varphi(c(m,k,r))} \right)^{n-k-r} 
     = \exp\{O(k n^{-2} \cdot k^3 p n^{-1})\} = 1 + o(k^2 n^{-2}).
\end{equation*}
Multiplying these two estimates completes the proof of the Lemma.
\end{proof}

In order to apply our results for $ G_{n,m}$ to $G_{n,p}$ we need to
identify the values of $m$ at which $k_0$ changes from one value to the next.
For $ n^{ 2/3 - \ep} < k < n^{3/4-\ep}$
define $m_0(k)$ be the smallest $m$ for which $k_0(m) = k$. 
\begin{lem} \label{lem: m to k ratio}

If $ n^{2/3 - \ep} < k < n^{3/4 - \ep}$ then
    \begin{equation*}
        m_0(k) - m_0(k+1) \sim \frac{n^2 \log(n/k) }{k^2 }.
    \end{equation*}

\begin{proof}
Note that $r_0$ is a function of $m$. For the purpose of this proof we write $ r_0(m,k)$ 
in place of $ r_0(k)$. Set $m' = m_0(k), \hat{m} = m_0(k+1), r' = r_0(m', k)$, and $\hr = r_0(\hm, k+1)$. 
We begin by noting that (\ref{equ: X' ratio over m}) and (\ref{equ: X' asymptotic formula}) and the definitions of $k_0$ and $r_0$ imply
    \begin{align}
        \max_{r \leq r_1}\{X'(m', k, r)\} &\sim n^{\ep} \label{equ: n to the epsilon 1st} \\
        \max_{r \leq r_1} \{X'(\hm, k+1, r)\} &\sim n^{\ep} \label{equ: n to the epsilon 2nd},
    \end{align} 
and it then follows from (\ref{equ: X' ratio over k}) that we have
    \begin{equation}
        \max_{r \leq r_1}\{X'(\hm, k, r)\} \sim e^2 n^{1 +\ep} k^{-1}. \label{equ: n to the epsilon prime}
    \end{equation} 
(Here we should note that $r_1$ does depend on $m$; however, $r_1$ varies only slightly for the range of values of $m$ we consider, and the values of $X'$ with $r$ near $r_1$ are much lower than the maximum value of $X'$ by (\ref{equ: X' ratio over k}), hence we can essentially treat $r_1$ as fixed for the proof of this Lemma). Now let $\delta > 0$ be an arbitrary constant. By (\ref{equ: X' ratio over m}) and (\ref{equ: n to the epsilon prime}), if $m - \hm < (1 - \delta)\frac{n^2 \log(n/k)}{k^2}$ then 
\begin{equation*}
    \max_{r \leq r_1}\{X'(m,k,r)\} \gg n^{\ep}
\end{equation*}
and if $m - \hm > (1 + \delta)\frac{n^2 \log(n/k)}{k^2}$ then 
\begin{equation*}
    \max_{r \leq r_1}\{X'(m,k,r)\} \ll n^{\ep}.
\end{equation*}
Therefore, by (\ref{equ: n to the epsilon 1st}), the Lemma holds.
\end{proof}

\end{lem}

    Now we prove Theorem~\ref{thm:Gnp} using Lemma~\ref{lem: m to k ratio} together with Lemmas~\ref{lem:lower},~\ref{lem:upper},~and~\ref{lem:picky} and the distribution of the number of edges of $G_{n,p}$. We prove each statement in turn.

    We begin with the concentration statement in part (ii) of Theorem~\ref{thm:Gnp}.
As the number of edges of $G_{n,p}$ is simply a binomial random variable with 
variance $\sim n^2p/2$, there exists some $m_{low}$ and $m_{high}$ such that $m_{high} - m_{low} = \ell n p^{1/2} $ 
        and $e(G_{n,p}) \in (m_{low},m_{high})$ whp. Then by Lemma \ref{lem: m to k ratio}, and recalling (\ref{equ: k used in this range}), we have 
        \[ k_0(m_{low}) - k_0(m_{high}) = \Theta\left( \frac{ \ell n p^{1/2}} { n^2 p^2 / (\log n)}\right) = \Theta\left( \ell ( \log n) p^{-3/2} n^{-1}\right). \] 
        Therefore, applying Lemma~\ref{lem:picky} (and noting that $k_0$ as a function of $m$ is decreasing), we see that $\alpha(G_{n,p})$ is concentrated on an interval of the specified length.

Now we turn to the anti-concentration statement in part~(ii) of Theorem~\ref{thm:Gnp}. Set $m = \binom{n}{2}p$, and let $ \hk = k_0(m)$. Consider $k_{low} = \hk - \frac{C}{2} \log(n/ \hk) n^{-1} p^{-3/2}$, and $k_{high} = \hk + \frac{C}{2} \log(n/ \hk) n^{-1} p^{-3/2}$. It is enough to show that there exists some constant $\ep_1 > 0$ such that $\PP[\alpha(G_{n,p}) < k_{low}] > \ep_1$ and $\Pr[\alpha(G_{n,p}) > k_{high}] > \ep_1$. We establish the first inequality; the other follows by a similar argument. 
        Let $m_{high} = m_0(k_{low} - 2)$. By Lemma~\ref{lem: m to k ratio}, we have 
        \begin{equation*}
            m_{high} - m \sim \frac{ n^2 \log( n/ \hk)}{ \hk^2}(\hk - k_{low} + 2) \sim \frac{C}{8} n \sqrt{p}.
        \end{equation*}
        Note that the probability of the event $ \{ e(G_{n,p}) > m_{high} \} $ is bounded below by some constant, say $ 2\ep_1$.
        Since the independence number is a decreasing function with respect to edges, if $m \geq m_{high}$ then $\alpha(G_{n,m}) < k_{low} $ whp. Thus,
        \[ \PP( \alpha( G_{n,p}) < k_{low}) \ge \PP( e( G_{n,p}) > m_{high}) - o(1) \ge \ep_1. \]
         
    Finally, we prove part~(i) of Theorem~\ref{thm:Gnp}.
        We first recall that the authors established 2 point concentration of $ \alpha( G_{n,p})$ with $ p > n^{-2/3 + \epsilon}$ in a previous work \cite{bhindy}.
        So we may assume $ ( \log n)^{2/3} n^{-2/3} \ll p < n^{-2/3+ \epsilon}$. Write
        $p = \ell n^{-2/3} \ln(n)^{2/3}$ where $ \ell = \omega(1)$.          
        Set 
        $ \hk = k_0(m)$ and $ \hr = r_0(m, \hk)$.
        Again choose  $m_{low}$ and $m_{high}$ such that $m_{high} - m_{low} = \ell n p^{1/2}$ and $e(G_{n,p}) \in (m_{low}, m_{high})$ whp.  
        Note that Lemma~\ref{lem: X' ratio over m} implies
        \begin{multline*} 
        X'( m_{low}, \hk , \hr ) \le n^{\epsilon} \exp\left\{    \ell n p^{1/2} \cdot \Theta( ( \log n)^2 n^{-2} p^{-2} )  \right\} \\
        =  n^{\epsilon} \exp\left\{    \Theta( \ell ( \log n)^2 n^{-1} p^{-3/2} )  \right\} = n^{\epsilon} \exp\left\{    \Theta( \ell^{-1/2} (\log n) ) \right\}
        = n^{( 1 + o(1)) \epsilon}.
        \end{multline*}
        Lemma~\ref{lem:upper} (and recalling (\ref{equ: X' ratio over k})) then implies that we have $ \alpha( G_{n,m_{low}}) \le \hk$ whp, and it follows that $ \alpha(G_{n,p}) \le \hk $ whp. A similar calculation shows $ X'( m_{high}, \hk - 1, \hr ) \ge  n^{ (1 -o(1))\epsilon} $. So it follows from Lemma~\ref{lem:lower} (again recalling (\ref{equ: X' ratio over k})) that we have $ \alpha( G_{n,m_{high}}) \ge \hk-1$ whp. Thus $ \alpha( G_{n,p}) \ge \hk -1$ whp, as desired.

\noindent {\bf Acknowledgement.} The first author thanks Lutz Warnke and Emily Zhu for many useful conversations and Mehtaab Sawhney for pointing out connections with the work of Ding, Sly and Sun.

\bibliographystyle{plain}
\bibliography{Bib}
\end{document}